\documentclass[a4paper,oneside,english,10pt]{elsarticle}

\usepackage{amssymb}
\usepackage{amsmath}
\usepackage{pifont}

\usepackage{tikz}

\usepackage{tcolorbox}
\newcommand{\RR}{{\mathbb R}}

\newcommand{\x}{{\mathbf x}}
\newcommand{\y}{{\mathbf y}}
\newcommand{\n}{{\mathbf n}}
\newcommand{\p}{{\mathbf P}}

\newcommand{\Omegap}{\Omega^i}
\newcommand{\Omegam}{\Omega}

\newcommand{\R}{\mathbb R}

\newcommand{\bfx}{{\mathbf x}}
\newcommand{\bfy}{{\mathbf y}}
\newcommand{\bfn}{{\mathbf n}}
\newcommand{\bfv}{{\mathbf v}}
\newcommand{\bfPhi}{{\boldsymbol\Phi}}

\newcommand{\bfB}{{\mathbf B}}
\newcommand{\T}{\mathrm{T}}
\newcommand{\wh}[1]{\widehat{#1}}
\newcommand{\Cex}{\mathrm{C}(\bfx,\epsilon)}

\newcommand{\Cexq}{\text{C}_q(\x,\epsilon)}

\journal{Computers \& Mathematics with Applications}

\begin{document}

	\begin{frontmatter}
	
	\title{CQ-collocation BEMs for exterior wave propagation problems: the 3D Neumann case%
		\tnoteref{t1}}
	\tnotetext[t1]{This research was supported by the GNCS-INDAM 2025 program
		``High-order BEM-based numerical techniques for wave propagation problems''.}

	\author[1]{D. Collato}
	\ead{davide.collato@polito.it}
	\author[1]{S. Falletta\corref{cor1}}
	\ead{silvia.falletta@polito.it}
	\author[1]{G. Monegato}
	\ead{giovanni.monegato@formerfaculty.polito.it}
	\author[1]{L. Scuderi}
	\ead{letizia.scuderi@polito.it}
	\cortext[cor1]{Corresponding author}

	\affiliation[1]{organization={Dipartimento di Scienze Matematiche, Politecnico di Torino},
		addressline={C.so Duca degli Abruzzi, 24},
		postcode={10129},
		city={Torino},
		country={Italy}}

\begin{abstract}
In this paper, we consider three-dimensional exterior wave propagation problems with non-homogeneous Neumann boundary conditions. We reformulate the problem as a space-time boundary integral equation involving a hypersingular operator. To solve the resulting formulation, we combine Lubich's convolution quadrature with a shifted collocation boundary element discretization and investigate numerically the efficiency, stability, and convergence of the resulting numerical scheme. Numerical tests show that the proposed approach achieves an accuracy comparable to that of Lubich's convolution quadrature combined with a Galerkin boundary element method, and provide numerical evidence of the unconditional stability of the collocation scheme.

\end{abstract}


\end{frontmatter}

\section{Introduction}\label{sec:intro}

Let $\Omegap\subset \RR^3$ be a bounded, open domain with Lipschitz boundary $\Gamma$, and define the exterior domain $\Omegam = \RR^3 \setminus \overline{\Omegap}$, containing a homogeneous, compressible acoustic medium. We consider the following problem governed by the scalar wave equation with a Neumann boundary condition
\begin{align}
\label{pb:esterno}
\begin{cases}
\begin{aligned}
&\frac{\partial^2 u}{\partial t^2}(\x,t) - \Delta u(\x,t) = 0 && (\x,t) \in \Omegam \times (0,T], \\[4pt]
&u(\x,0) = 0, \quad \frac{\partial u}{\partial t}(\x,0) = 0 && \x \in  \Omegam,\\[4pt]
&\frac{\partial u}{\partial \n_\x}(\x,t) = g(\x,t) && (\x,t) \in  \Gamma \times (0,T],
\end{aligned}
\end{cases}
\end{align}
where $\n_\x$ denotes the outward unit normal to $\Gamma$ at $\x$. 

 Existence and uniqueness of the solution to the initial-boundary value
 Problem \eqref{pb:esterno} have been proved in \cite{BaHa} provided that \(g\) belonging to a suitable
 anisotropic Sobolev space. To report this result, we 
 consider the usual Sobolev space \(H^\alpha(D)\),
 \(\alpha \in [-k,k]\), $D\subset\R^3$ and \(k=1\) for the Lipschitz case considered here. Denoting the norm on this space by $\|\cdot\|_{H^\alpha(D)}$, the anisotropic Sobolev space \(H^r(\mathbb{R};H^\alpha(D))\) of order
 \(r \in \mathbb{R}\) is defined by
 $$
 H^r(\mathbb{R};H^\alpha(D))
 =
 \left\{
 f\hskip-0.1cm:\hskip-0.1cm D \times \mathbb{R} \to \mathbb{R}
 \;:\;
 \int_{\mathbb{R}}
 (1+|\omega|)^{2r}
 \,
 \|\mathcal{F}f(\cdot,\omega)\|_{H^\alpha(D)}^2
 \,d\omega
 <\infty
 \right\},
 $$
 where \(\mathcal{F}\) denotes the Fourier transform on \(\mathbb{R}\).
 This space can be restricted to a finite time interval by defining, as in \cite{Lubich2},
 \begin{equation*}
 	H_0^r(0,T;H^\alpha(D))
 	=
 	\left\{
 	f|_{D\times(0,T)}
 	\;:\;
 	f\in H^r(\mathbb{R};H^\alpha(D)),
 \ f\equiv0
 	\text{ on }
 	D\times(-\infty,0)
 	\right\}.
 	\label{eq:Hr0}
 \end{equation*}
 Then, for a  boundary datum
 $
 g\in H_0^{r+1}(0,T;H^{-1/2}(\Gamma)),
 $
 there exists a unique solution
 $
 u\in H_0^r(0,T;H^1(\Omega))
 $
 of Problem \eqref{pb:esterno},
 depending continuously on the boundary datum. 
 
 Note that for integer \(r\geq 0\), the condition on the datum implies that \(g\) and its first \(r\) time derivatives vanish at \(t=0\), so that for $r\geq 1$ the compatibility conditions at \(t=0\) are  automatically satisfied.

\

Exterior Neumann problems of the form \eqref{pb:esterno} arise in several wave propagation applications. A notable example is time-domain sound-hard scattering, where $\Omegap$ denotes a rigid obstacle and $u$ is the scattered acoustic field. In this case, the total field $u^{\text{tot}}$ is decomposed as
\(
u^{\text{tot}}(\x,t)=u^{\text{inc}}(\x,t)+u(\x,t),
\)
where $u^{\text{inc}}$ is the incident field. The rigidity of the obstacle imposes the Neumann boundary condition
\(\partial u^{\text{tot}}/\partial \n_\x=0\)
which, in terms of the scattered field, yields
\[
g(\x,t)=\frac{\partial u}{\partial \n_\x}
=-\frac{\partial u^{\text{inc}}}{\partial \n_\x}.
\]
For a broader context on such problems, we refer to \cite{ColtonKress2013}, \cite{Lax} and \cite{Moench1996}.

\

The aim of this paper is to propose a numerical strategy for solving Problem~\eqref{pb:esterno}, based on its reformulation as a space-time boundary integral equation, using the double-layer potential representation of the solution $u$. This formulation leads to an integral equation that involves a hypersingular operator acting on the unknown boundary density. For the discretization of the time integral, we adopt the Convolution Quadrature (CQ) rule introduced by Lubich \cite{Lubich1, Lubich2}, which provides a robust and stable tool in the numerical solution of the boundary integral equation associated with wave propagation problems. 
The Lubich CQ, combined with a spatial Galerkin boundary-element discretization, was introduced in \cite{Lubich2} for exterior 3D Dirichlet wave propagation problems and later studied in \cite{Chappell} for exterior 3D Neumann ones. For both problems, unconditional stability and convergence results were proved. Related works include the study of some fast algorithms for the implementation of the CQ-Galerkin method (see, for example, \cite{BaSau}, \cite{HaKrSa}), its extension to multistep and multistage based CQ (see, for example, \cite{Banjiai_RK}), and its application to other engineering problems such as, for instance, elastodynamic problems (see, for example, \cite{KieSch},\cite{FMS_bemelasto}).

Alongside Galerkin methods, collocation-based discretizations have attracted considerable attention owing to their reduced computational complexity, which allows their application in industrial contexts. In particular, collocation methods avoid the evaluation of double integrals, making them especially appealing for large-scale simulations. Nevertheless, their theoretical understanding of stability and convergence generally remains incomplete. 

In this paper, we propose a discretization strategy combining Lubich's convolution quadrature in time with a shifted collocation method in space. 
The shifted collocation approach is suited for Lipschitz surfaces, where the unit normal vector $\n_{\x}$ is not defined at edges and corners.
The hypersingular space-time integral equation is first regularized by means of integration by parts, leading to a strongly singular formulation. The resulting spatial integrals are then evaluated using tailored cubature formulas specifically designed for this setting.
Extensive numerical experiments are provided to support the validity and robustness of the proposed approach. In particular, in the absence of a formal proof, the numerical experiments indicate unconditional stability; indeed, apparent numerical instabilities, such as spurious oscillations, are found to disappear upon increasing the accuracy of the cubature rule used to evaluate the boundary integrals. This confirms the crucial role of cubature accuracy in guaranteeing stability. These findings extend and complement our previous works on 2D and 3D Dirichlet cases (\cite{FMS_Dir,FS_fascette,FMS_altavelocita}) and on 2D Neumann configurations (\cite{FMS_Neu}), where the Lubich-CQ approach was successfully applied with collocation-based spatial discretizations. In \cite{FS_fascette}, we also 
compared our approach applied to a 3D Dirichlet problem with other time-domain collocation boundary element formulations based on retarded potential operators, noting that the latter were conditionally stable under a CFL condition. By contrast, the Lubich-based scheme was shown to be stable for any choice of space-time discretization parameters.

This paper completes our investigation by addressing the application of the CQ-collocation framework to wave propagation problems in three-dimensional Neumann problems, providing a novel, robust and efficient boundary element method. Although the underlying idea is inspired by the two-dimensional CQ-shifted collocation method proposed in \cite{FMS_Neu}, the extension to three dimensions requires a substantially different numerical treatment. In particular, the use of surface meshes, the numerical computation of the boundary integrals, the treatment of strongly singular integrals, and the choice of shifted collocation points all involve new challenges specific to the three-dimensional setting.

\

The paper is organized as follows. In Section~\ref{sec:bem}, we derive the space--time boundary integral formulation associated with Problem~\eqref{pb:esterno}, together with its temporal and spatial discretization. Therein, the collocation approach is presented by detailing the selection of collocation points and the numerical cubature strategy.
Since no theoretical stability or convergence results are currently available for the 3D collocation method, in Section~\ref{sec:num_res} we present a numerical investigation of these properties through some benchmark problems. The analysis includes a comparison with the CQ-Galerkin method, showing that both approaches yield comparable accuracy and an optimal conditioning of the final linear system. Nonetheless, the collocation method achieves these results at a substantially lower computational cost, making it a more efficient and competitive alternative in practice. 
Finally, the main conclusions are summarized in Section~\ref{sec:conclusions}.

\section{The boundary element method}\label{sec:bem}
Let us assume that $\Gamma$ and the boundary data of Problem~\eqref{pb:esterno} satisfy the regularity and compatibility conditions described at the beginning of Section~\ref{sec:intro}. Under these assumptions, its solution can be represented by means of the following double layer potential
\begin{equation}
u(\x,t) = (D \varphi)(\x,t) := \int_\Gamma \int_0^t 
\frac{\partial}{\partial \n_\y} G(\x-\y, t-\tau)\,\varphi(\y,\tau)\, \text{d}\tau\,\text{d}\Gamma_\y ,
\ \x \in \Omegam, \; t \in (0,T],
\label{eq:double_layer}
\end{equation}
where 
\begin{equation}
G(\x,t) = \frac{1}{4\pi ||\x||}\,\delta\!\left(t - ||\x||\right)
\label{eq:fundamental_solution}
\end{equation}
denotes the fundamental solution of the wave equation,  $\delta$ being the Dirac delta distribution. The unknown density function $\varphi$ represents the jump of the solution $u$ across the boundary (see \cite[Proposition~2.1]{FMS_Neu} for the two-dimensional case; the corresponding extension to the three-dimensional case follows analogously).
To derive the boundary integral equation considered here, we follow \cite[Section 6.5]{Steinbach2008} and define the following operator
\begin{equation}
\begin{aligned}
(W \varphi)(\x,t) &:=  \n_\x \cdot \lim_{\tilde{\x}\in\Omega,\tilde{\x} \to \x\in\Gamma} 
\nabla_{\tilde\x} (D\varphi)(\tilde\x, t)\\
& = \n_\x \cdot \lim_{\tilde{\x}\in\Omega,\tilde{\x} \to \x\in\Gamma}
\nabla_{\tilde\x} 
\int_\Gamma  \int_0^t
\frac{\partial}{\partial \n_\y} G(\tilde\x-\y, t-\tau)\,\varphi(\y,\tau)\, \text{d}\tau\,\text{d}\Gamma_\y \\
& =  \lim_{\epsilon\to 0^+} 
\int\limits_{\substack{\mathbf{y}\in{\Gamma}:\\\|\mathbf{x}-\mathbf{y}\|>\epsilon}}  
\int_0^t \frac{\partial^2}{\partial \n_\x\partial \n_\y} G(\x-\y, t-\tau)\,\varphi(\y,\tau)\, \text{d}\tau\,\text{d}\Gamma_\y.
\end{aligned}
\label{eq:W_operator}
\end{equation}
For $\x\in\Gamma$, the limit as $\epsilon\to0^+$ of the above integral does not exist in the sense of the Cauchy principal value. Indeed, taking into account the expression of the fundamental solution given in \eqref{eq:fundamental_solution}, the dominant singular behaviour is determined by the term
$$
\frac{\partial^2}{\partial {\n}_{\mathbf{x}}\partial {\n}_{\mathbf{y}}}
\left(
\frac{1}{\|\mathbf{x}-\mathbf{y}\|}
\right) = 	\frac{
	\n_{\mathbf{x}}\cdot\n_{\mathbf{y}}
}
{\|\mathbf{x}-\mathbf{y}\|^3}
-
3
\frac{
	(\n_{\mathbf{x}}\cdot(\mathbf{x}-\mathbf{y}))
	(\n_{\mathbf{y}}\cdot(\mathbf{x}-\mathbf{y}))
}
{\|\mathbf{x}-\mathbf{y}\|^5},
$$
which is of order
$\mathcal{O}\left(\|\mathbf{x}-\mathbf{y}\|^{-3}\right)$ as
$\mathbf{y}\rightarrow\mathbf{x}$. 
Consequently, the singularity of the kernel in \eqref{eq:W_operator} is too strong for the corresponding integral to be interpreted as a Cauchy principal value integral. Therefore, the operator $W$ is referred to as a hypersingular boundary integral operator. To obtain an explicit representation of $W$, an appropriate regularization procedure is required. A common approach consists in subtracting the leading singularity which can then be treated by analytical techniques or special quadrature rules (see, e.g., \cite{Steinbach2008}).
 Here, we rely on the alternative approach based on integration by parts, which reduces the singularity of the kernel and leads to either a strongly singular formulation within the collocation approach (see \eqref{eq:Dkappa_tangential}) or a weakly singular formulation within the Galerkin method (see \eqref{eq:W_galerk}).

Finally, combining \eqref{eq:W_operator} with \eqref{eq:double_layer} and imposing the Neumann boundary condition $\partial u/\partial \n_\x = g$, we obtain the following boundary integral equation in the unknown density $\varphi$
\begin{equation}
(W \varphi)(\x,t) = g(\x,t),
\quad \x \in \Gamma, t \in (0,T],
\label{eq:BIE}
\end{equation}
where the hypersingular boundary integral operator $W$ is defined as in \eqref{eq:W_operator} and is formally written as
\begin{equation}
(W \varphi)(\x,t)
= 
\int_\Gamma \int_0^t \frac{\partial^2}{\partial \n_\x\partial \n_\y}  G(\x-\y, t-\tau)\varphi(\y,\tau)\,  \text{d}\tau\,\text{d}\Gamma_\y.
\label{eq:hyper_W}
\end{equation}

In the next sections, we describe the time and space discretization of \eqref{eq:hyper_W}, the former based on a Lubich-CQ and the latter on a shifted collocation method.

\subsection{Time discretization: convolution quadrature}
We first remark that the hypersingular integral in \eqref{eq:hyper_W} is regularized through an integration-by-parts procedure (see \eqref{eq:Dkappa_tangential}) before the space discretization is performed. Therefore, before introducing the time discretization, it is convenient to rewrite \eqref{eq:BIE} in the following form
\begin{equation}
\frac{\partial}{\partial \n_\x} 
\int_\Gamma \,\frac{\partial }{\partial \n_\y}  \left(\int_0^t G(\x-\y, t-\tau)\,\varphi(\y,\tau)  \, \text{d}\tau\right) \,\text{d}\Gamma_\y  = g(\x,t),
\, \x \in \Gamma, \; t \in (0,T].
\label{eq:hyper_W2}
\end{equation}
To approximate the time integral in (\ref{eq:hyper_W2}), we decompose the interval $[0,T]$ into $N$
subintervals of equal length $\Delta_t = T/N$ and collocate \eqref{eq:hyper_W2} at the discrete time levels
$t_n = n \Delta_t$, with $n = 0, \ldots, N$,
\begin{equation}
\frac{\partial}{\partial \n_\x} 
\int_\Gamma \frac{\partial }{\partial \n_\y}  \left(\int_0^{t_n} G(\x-\y, t_n-\tau)\,\varphi(\y,\tau)  \, \text{d}\tau\right) \,\text{d}\Gamma_\y  = g(\x,t_n),
\quad \x \in \Gamma.
\label{eq:hyper_W2_t}
\end{equation}
The time integrals  in \eqref{eq:hyper_W2_t} are then discretized by means of the Lubich
convolution quadrature rule associated with the second order BDF method (BDF$2$)
(see \cite{Lubich1})
\begin{equation}\label{E:omega}
\int_0^{t_n} G(\x-\y, t_n-\tau)\,\varphi(\y,\tau)  \, \text{d}\tau \approx  \sum_{j=0}^n
\omega_{n-j}(\Delta_t;r)\,\varphi(\y,t_j),
\quad r = \|\x - \y\|,
\end{equation}
where the coefficients $\omega_n$ are defined by
\begin{equation}\label{eq:coeffwn}
\omega_n (\Delta_t;r) =
\frac{1}{2\pi \imath}\int_{|z|=\rho}
K\!\left(r,\frac{\gamma(z)}{\Delta_t}\right) z^{-(n+1)}\,\text{d}z,\quad n = 0,\ldots,N.
\end{equation}
In \eqref{eq:coeffwn}, $K$ denotes the Laplace transform of the fundamental solution $G$, whose expression is 
\begin{equation}\label{eq:lapl_transf}
K(r,s) = \frac{1}{4\pi} \frac{e^{-rs}}{r}, 
\end{equation}
$\gamma(z) = \tfrac{3}{2} - 2z + \tfrac{1}{2}z^2$ is
the characteristic quotient of the BDF$2$ method and the parameter $\rho$ is chosen such that, for $\vert z\vert \leq \rho$, the corresponding $\gamma(z)$ lies
in the domain of analyticity of $K$ (see \cite{Lubich1}).  

Introducing the polar parametrization $z = \rho e^{\imath\vartheta}$, the integrals in \eqref{eq:coeffwn}
can be efficiently computed via a trapezoidal rule with $L$ equal steps of length $2\pi/L$ as follows
\begin{equation}\label{E:ot}
\omega_{n}(\Delta_t;r) \approx
\frac{\rho^{-n}}{L}\sum_{\ell=0}^{L-1}
K\!\left(r, \frac{\gamma\!\left(\rho e^{\imath 2\pi \ell / L}\right)}{\Delta_t}\right)
e^{-\imath n 2\pi \ell/L}.
\end{equation}
Denoting by $\varphi^{(j)}_{\Delta_t}(\x)\approx \varphi(\x,t_j)$ the approximation of the exact solution of \eqref{eq:hyper_W2} evaluated at the time instants $t_j$, and inserting \eqref{E:ot} in \eqref{eq:hyper_W2_t}, we obtain
\begin{align}
&\sum_{j=0}^n \frac{\rho^{-(n-j)}}{L}\sum_{\ell=0}^{L-1} \left(\frac{\partial}{\partial \n_\x} 
\int_\Gamma  \frac{\partial}{\partial \n_\y} 
K\!\left(r, \frac{\gamma\!\left(\rho e^{\imath 2\pi \ell / L}\right)}{\Delta_t}\right)\,\varphi^{(j)}_{\Delta_t}(\y) \text{d}\Gamma_\y\right)e^{-\imath (n-j) 2\pi \ell/L} \nonumber\\
& \qquad = g(\x,t_n), \quad n=0,\ldots,N, \ \x \in \Gamma. 
\label{eq:hyper_W2_t_bis}
\end{align}
		For the numerical experiments of Section \ref{sec:num_res}, in~\eqref{E:ot} we set \(L=N+1\) and, to compute the coefficients $\omega_n$ with accuracy $\varepsilon$, \(\rho=\varepsilon^{1/N}\). The choice of these parameters is theoretically justified by Chappell in~\cite[Theorem 5.1]{Chappell} for the Neumann Problem \eqref{pb:esterno} discretized by a Galerkin boundary element method.
	
	To investigate the effectiveness of these parameter values, we also considered the alternative choice \(L=2N\) and \(\rho=\varepsilon^{1/(2N)}\), proposed by Lubich in~\cite{Lubich1} for the Dirichlet problem combined with the Galerkin approach. The two parameter sets yield the same accuracy for our approach, the former choice being preferable from a computational point of view, as it requires the evaluation of the boundary element matrices for a smaller number of frequencies.
	Moreover, in the numerical examples, we set \( \varepsilon=10^{-10} \), which was found to be sufficient for the accuracy of the numerical results reported in Section~\ref{sec:num_res}.

\

We finally remark that a key advantage of the convolution quadrature method is that all the coefficients \( \omega_n \) can be computed simultaneously using the FFT, with a computational cost of \( \mathcal{O}(N \log N) \) flops.

\subsection{Space discretization: collocation method} 

To describe the space discretization of \eqref{eq:hyper_W2_t_bis}, we start by introducing a conforming decomposition $\mathcal{T}_h$ of the surface $\Gamma$ into planar triangular elements $\T\in\mathcal{T}_h$. This allows to define the piecewise triangular surface $\Gamma_\Delta = \displaystyle{\cup_{\T\in\mathcal{T}_h} \T}$, interpolating $\Gamma$ at the triangle vertices $\{\x_{i}, i=1,\dots,M\}$. We denote by $h$ the mesh size of $\Gamma_\Delta$, which is given by the maximum triangle diameter.

At each time instant $t_n$, the unknown function $\varphi_{\Delta_t}(\x,t_n)$ is then approximated by
\begin{equation*}
\varphi^{(n)}_{\Delta_t,h}(\x) := \sum_{i=1}^{M} \varphi_{i}^{(n)} \Psi_{i}(\x),\quad \x\in\Gamma_\Delta,
\end{equation*}
where $\{\Psi_{i}\}_{i=1}^{M}$ are the classical continuous piecewise linear basis functions associated with the chosen triangulation, and $\varphi_{i}^{(n)} = \varphi^{(n)}_{\Delta_t,h}(\x_i)\approx \varphi(\x_{i},t_n)$.

Therefore, the discrete collocation problem consists of finding
$\varphi^{(n)}_{i}$, for $i = 1,\ldots,M$ and $n=0,\ldots,N$ 
such that
\begin{align}\label{E:BIEc}
&\sum_{j=0}^n \frac{\rho^{-(n-j)}}{L} 
\sum_{\ell=0}^{L-1} 
\Bigg[
\sum_{i=1}^M \varphi^{(j)}_{i} 
\frac{\partial}{\partial \n_{\boldsymbol{\xi}_m}} 
\int_\Gamma
\frac{\partial}{\partial \n_\y}  
K\!\Big(r_m, \frac{\gamma\!\big(\rho e^{\frac{\imath 2\pi \ell}{L}}\big)}{\Delta_t}\Big)
\Psi_i(\y)\, \text{d}\Gamma_\y
\Bigg] e^{-\frac{\imath (n-j) 2\pi \ell}{L}}
\notag\\
& \qquad  = g(\boldsymbol{\xi}_m,t_n),
\end{align}
where $r_m = \|\boldsymbol{\xi}_m -\y\|$ for $m=1,\ldots,M$, and $\boldsymbol{\xi}_m$ are the collocation points defined in Section \ref{sec:coll}.
Before proceeding numerically, we apply an integration by parts strategy to the space integrals appearing in (\ref{eq:hyper_W2_t_bis}), which effectively weakens the hypersingularity and enables an easier numerical discretization via collocation techniques. 
To this aim, we observe that \eqref{eq:lapl_transf} coincides with the fundamental solution $U^*_\kappa(\x,\y)=\frac{1}{4\pi}\frac{e^{\imath \kappa r}}{r}$ of the Helmholtz equation with wavenumber \(\kappa\). Hence, by setting \(s=-\imath\kappa\), we apply Equation (8.8) in \cite{Steinbach2013}, derived by integration by parts for \(U^*_\kappa\) and valid for Lipschitz boundaries, which in our case takes the form
	\begin{align}
	\frac{\partial}{\partial \n_\x} \int_\Gamma  \frac{\partial}{\partial \n_\y} 
	K(r, s)\,v(\y) \, \mathrm{d}\Gamma_\y
	&= \int_\Gamma 
	\Big[ \big( \n_\x \times \nabla_\x K(r,s) \big) \cdot \big( \n_\y \times \nabla v(\y) \big) \notag\\
	&\quad - \, s^2\, (\n_\x \cdot \n_\y)\, K(r,s)\, v(\y) \Big] \, \mathrm{d}\Gamma_\y.
	\label{eq:Dkappa_tangential}
	\end{align}
We point out that the reformulation \eqref{eq:Dkappa_tangential} involves tangential rather than normal derivatives, and that the term $\n_\y \times \nabla v(\y)$, with $\y \in \Gamma$, represents the surface curl of the boundary function $v$.
Further, we highlight that the hyper-singularity of \eqref{eq:hyper_W} has been weakened. In particular, the first term in the right hand side of \eqref{eq:Dkappa_tangential} is of order $\mathcal{O}(\|\x-\y\|^{-2})$, and its integral is defined in the Cauchy principal value sense as $\y$ approaches $\x$, whereas the second term is of order $\mathcal{O}(\|\x-\y\|^{-1})$ and the corresponding integral is weakly singular. The numerical treatment of these singular integrals is detailed in Section \ref{sec:cub}.
 
 \
 
To define the final linear system associated with the whole space-time discretization, we introduce the matrices $\mathbf{A}^{(n)} \in \R^{M\times M}$, $n = 0,\ldots,N$, whose generic entry, by using \eqref{E:BIEc} and \eqref{eq:Dkappa_tangential}, and restricting the integration over the support of the basis functions $\Psi_i$, is given by

\begin{align}
\mathbf{A}^{(n)}_{m,i}
&= \frac{\rho^{-n}}{L}\sum_{\ell=0}^{L-1}
\Biggl(
\displaystyle\int_{\operatorname{supp}(\Psi_i)}
\Bigl[
\Bigl( \n_{\boldsymbol{\xi}_m} \times \nabla_\x
K\Bigl(r_m,\tfrac{\gamma(\rho e^{\frac{\imath2\pi\ell}{L}})}{\Delta_t}\Bigr) \Bigr)
\cdot
\bigl( \n_\y \times \nabla \Psi_i(\y) \bigr)
\notag\\[1ex]
&
-
\Bigl(\tfrac{\gamma(\rho e^{\frac{\imath2\pi\ell}{L}})}{\Delta_t}\Bigr)^2
(\n_{\boldsymbol{\xi}_m}\cdot\n_\y)\,
K\Bigl(r_m,\tfrac{\gamma(\rho e^{\frac{\imath2\pi\ell}{L}})}{\Delta_t}\Bigr)
\Psi_i(\y)
\Bigr]
\,\mathrm{d}\Gamma_\y
\Biggr)
e^{-\frac{\imath n2\pi\ell}{L}}.
\label{eq:matrix_entries}
\end{align}

 The collocation method consists of finding the  unknown vector $\pmb{\varphi}^{(j)}=(\varphi_{1}^{(j)},\ldots,\\ \varphi_{M}^{(j)})^T$,
$j=0,\ldots,N$, solution of the following Toeplitz block lower
triangular linear system
\begin{equation*}
\sum_{j=0}^n \mathbf{A}^{(n-j)}\pmb{\varphi}^{(j)}=\mathbf{g}^{(n)},\quad
n=0,\ldots , N,
\end{equation*}
with $\mathbf{g}^{(n)} = (g_{1}^{(n)},\ldots,g_{M}^{(n)})^T = (g(\boldsymbol{\xi}_1,t_n),\ldots,g(\boldsymbol{\xi}_M,t_n))^T$.
Since we assume that $g(x,t_0)\equiv 0$ (see Section \ref{sec:num_res}), and since the matrix $\mathbf{A}^{(0)}$ is nonsingular, we have $\pmb{\varphi}^{(0)}=\mathbf{0}$.
Hence, the solution of the linear system is obtained with a block forward substitution; in particular, for all instant $t_n$, we solve reduced linear systems of the type
\begin{equation}\label{E:SYS}
\mathbf{A}^{(0)}\pmb{\varphi}^{(n)}=\mathbf{g}^{(n)} - \sum_{j=0}^{n-1} \mathbf{A}^{(n-j)}\pmb{\varphi}^{(j)}, \quad n=1,\ldots , N,
\end{equation}
all having the same coefficient matrix $\mathbf{A}^{(0)}$.

Finally, once $\varphi^{(n)}_{\Delta_t,h}$ has been computed, we recover the approximation $u^{(n)}_{\Delta_t,h}(\p)$ to the solution of Problem~\eqref{pb:esterno} at the point $\p\in\Omegam$ and time $t_n$ by means of the discrete double layer potential representation \eqref{eq:double_layer}, obtained by applying Lubich-CQ to the double layer potential operator. Specifically,
\begin{align}
&u^{(n)}_{\Delta_t,h}(\p) = \nonumber\\
&\sum_{j=0}^n \frac{\rho^{-(n-j)}}{L}\sum_{\ell=0}^{L-1} \left(
\int_\Gamma  \frac{\partial}{\partial \n_\y} 
K\!\left(r_{_\p}, \frac{\gamma\!\left(\rho e^{\frac{\imath 2\pi \ell}{L}}\right)}{\Delta_t}\right)\,\varphi^{(j)}_{\Delta_t,h}(\y) \text{d}\Gamma_\y\right)e^{-\frac{\imath (n-j) 2\pi \ell}{L}}
\label{eq:u_rec}
\end{align}
for $n = 0,\ldots,N$ and $r_{_\p}=\|\p-\y\|$.

\subsubsection{Collocation points}\label{sec:coll}

Since we consider Lipschitz surfaces, where the unit normal vector $\n_{\x}$ may not be defined at edges and corner points, a standard nodal collocation approach cannot be employed. We therefore adopt a shifted collocation strategy. The shifted nodes are selected according to the procedure proposed in \cite{Wu_Ye_Jiang}, which consists in moving the collocation points from the vertices towards the interior of the mesh elements, so that they are evenly distributed over the discretized geometry and do not coincide with vertices or lie on edges.

\begin{figure}[h]
	\centering
	\begin{minipage}[t]{0.30\textwidth}
		\centering
		\includegraphics[width=3.5cm,height=3.5cm]{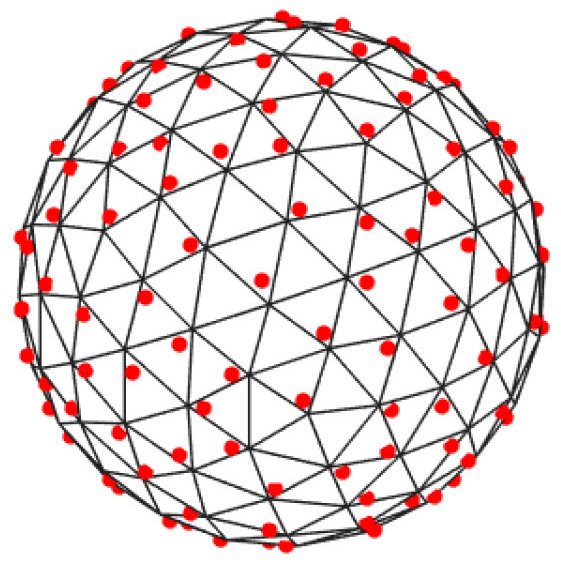}
	\end{minipage}
	\hfill
	\begin{minipage}[t]{0.30\textwidth}
		\centering
		\includegraphics[width=3.8cm,height=3.8cm]{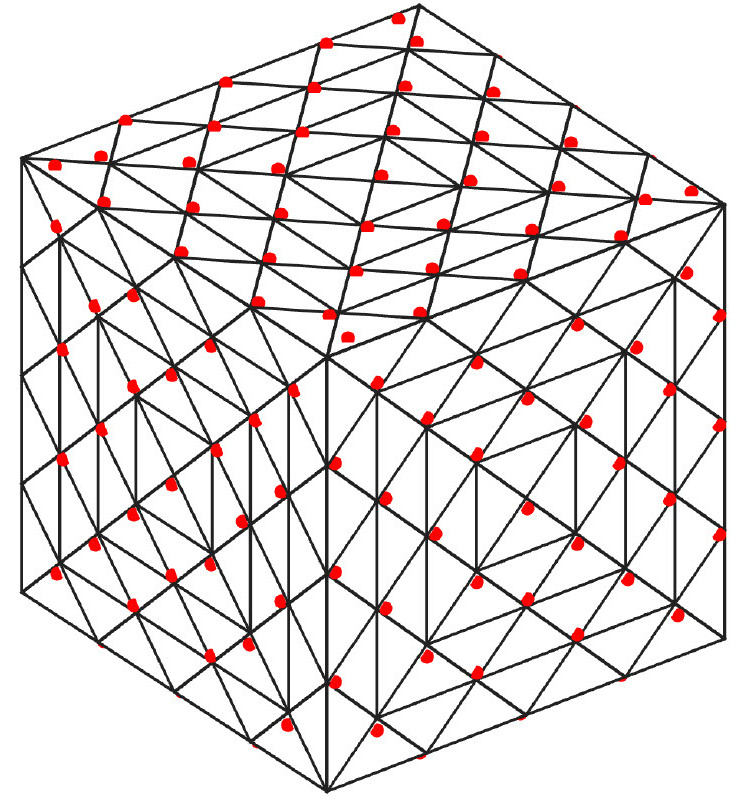}
	\end{minipage}
	\hfill
	\begin{minipage}[t]{0.30\textwidth}
		\centering
		\includegraphics[height=4.2cm]{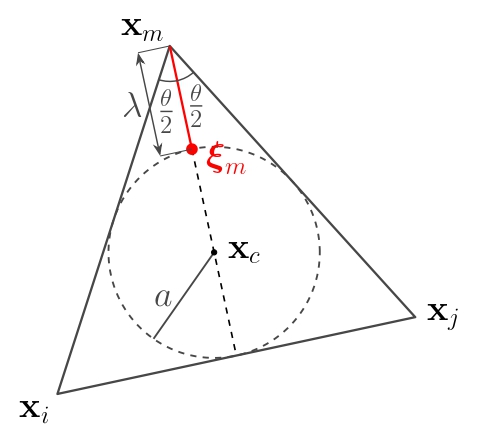}
	\end{minipage}
	\caption{\footnotesize{Triangulations of a sphere and a cube with shifted collocation points, and shifted collocation point placement.}}
	\label{fig:e-coll}
\end{figure}

For convenience of the reader, we report here the above mentioned strategy.
Denoting by $\{\x_m\}_{m=1}^{M}$ the set of all vertices of the triangles in the surface mesh, let $\mathfrak{T}_m$ denote the number of triangles $\{{\T}_{m,\ell}\}_{\ell=1,\ldots,\mathfrak{T}_m}$ sharing the vertex $\x_m$.  We sort the elements \( {\T}_{m,\ell} \) in descending order with respect to a parameter, which is determined by the weight of the vertex angle and the area of the sharing elements, and is defined as
\begin{equation*}
\eta_{m,\ell} = \frac{\angle_{{\T}_{m,\ell}}}{\sum_{j=1}^{\mathfrak{T}_m} \angle_{{\T}_{m,j}}}
+
\frac{A_{{\T}_{m,\ell}}}{\sum_{j=1}^{\mathfrak{T}_m} A_{{\T}_{m,j}}},
\quad \ell = 1, 2, \ldots, \mathfrak{T}_m,
\end{equation*}
where \( \angle_{{\T}_{m,\ell}} \) and \( A_{{\T}_{m,\ell}} \) denote the angle at \( \x_m \) and the area of the element  \( {\T}_{m,\ell} \), respectively. 
The \( m \)-th collocation node $\boldsymbol{\xi}_m$ is shifted, starting from $\x_m$, inside to the element ${\T}_{m,\bar{\ell}}$ corresponding to the largest value of \( \eta_{m,{\ell}} \), along the associated bisector of the angle $\theta$, by a properly  prescribed distance \( \lambda \). Denoting by 
$\x_m, \x_i, \x_j$ the vertices of the triangle ${\T}_{m,\bar{\ell}}$, 
once the distance \( \lambda \) is defined, the position of the collocation point \( \boldsymbol{\xi}_m \) is obtained by solving the following linear system
\begin{equation*}
\begin{bmatrix}
(\mathbf{x}_i - \mathbf{x}_m)^T \\
(\mathbf{x}_j - \mathbf{x}_m)^T \\
\mathbf{n}_{m,\bar{\ell}}^T
\end{bmatrix}
\boldsymbol{\xi}_m =
\begin{bmatrix}
\lambda \left| \mathbf{x}_i - \mathbf{x}_m \right| \cos(\theta/2) + (\mathbf{x}_i - \mathbf{x}_m)^T \mathbf{x}_m \\
\lambda \left| \mathbf{x}_j - \mathbf{x}_m \right| \cos(\theta/2) + (\mathbf{x}_j - \mathbf{x}_m)^T \mathbf{x}_m \\
\mathbf{n}_{m,\bar{\ell}}^T \mathbf{x}_m
\end{bmatrix},
\end{equation*}
\( \mathbf{n}_{m,\bar{\ell}} \) being the normal vector to the triangular element ${\T}_{m,\bar{\ell}}$. The distance $\lambda$ is defined as 
$\lambda = \|\x_c - \x_m \| - a$,
where $\x_c$ and $a$ are the center and radius of the incircle of the triangle ${\T}_{m,\bar{\ell}}$, respectively. In Figure \ref{fig:e-coll} we show the distribution of the shifted collocation points for the two reference domains used in the numerical examples and a triangular element containing the shifted collocation point \(\boldsymbol{\xi}_m\), which, according to the prescribed distance \(\lambda\), lies at the intersection of the angle bisector originating from \(\mathbf{x}_m\) and the incircle of the triangle ${\T}_{m,\bar{\ell}}$. As stated in \cite{Wu_Ye_Jiang}, for uniform meshes (as those we have considered here) the choice of \(\lambda\) there proposed and here adopted ensures a uniform distribution of the shifted points.

We have performed additional numerical tests to investigate the sensitivity of our CQ-collocation method to the choice of \(\lambda\). We observed that, by moving the collocation point \(\boldsymbol{\xi}_m\) along the angle bisector, the numerical accuracy remains essentially unchanged, as long as the shifted point lies close to the original point \({\bf x}_m\) and, in particular, does not move beyond the point \(\x_c\). These results indicate that the method is not particularly sensitive to the specific choice of \(\lambda\).

\subsubsection{Numerical cubature}\label{sec:cub}

In this section we provide details on the numerical computation of the space integrals appearing in the definition of the matrix entries \eqref{eq:matrix_entries}.  To this aim, we preliminarily remark that, since $\operatorname{supp}(\psi_i) = \cup_{\ell\in\mathcal{\T}_i}  {\T}_{i,\ell}$, the prototype integral we focus on is restricted to a triangle. To simplify the notation, we denote by $\T$ the generic integration triangle, and by $\x$ the generic collocation point.
The integrals involved in \eqref{eq:matrix_entries}, associated with the shifted collocation method, can be classified into four types:
	
	\begin{itemize}
		\item[$i)$] regular integrals, corresponding to a collocation point $\x$ located in a triangle that does not share any vertex with $\T$;
		
		\item[$ii)$] near singular integrals, corresponding to a collocation point $\x$ belongs to a triangle that share at least one vertex with $\T$;
		
		\item[$iii$)] weakly singular integrals, arising when the kernel exhibits a singularity of order $\mathcal{O}(1/r)$ for $r\to 0$ and the collocation point $\x$ belongs to $\T$;
		
		\item[$iv)$] strongly singular integrals, arising when the kernel exhibits a singularity of order $\mathcal{O}(1/r^2)$ for $r\to 0$ and the collocation point $\x$ belongs to $\T$.

	\end{itemize}

	 Different numerical strategies are then adopted, depending on the specific case.
	 
	 \
	 
	For case $i)$ we rely on the cubature rules proposed in \cite{SommVia}, which are defined over general polygonal domains and characterized by a prescribed degree of exactness. In particular, these formulas are exact for all bivariate polynomials of total degree up to $\delta$, with a number of nodes not exceeding the dimension of the corresponding polynomial space, namely
	$N_\delta=(\delta+1)(\delta+2)/2$ (see \cite[Table 2]{SommVia}).
	Moreover, these formulas guarantee positive weights and interior nodes, which are desirable properties in terms of stability and geometric robustness.

\

For case $ii)$, when the collocation point $\x$ does not belong to $\T$ but is located in its vicinity, the integrand may exhibit near-singular behaviour, potentially compromising the accuracy of the numerical integration. To mitigate this effect, whenever $\T$ shares at least one vertex with the triangle to which $\x$ belongs, we employ the same cubature rule as in case $i)$, but with a higher degree of exactness.

	\

	For case $iii)$ the triangle $\T$ to which the collocation point $\x$ belongs, is divided into three sub-triangles by connecting its vertices to $\x$. Each sub-triangle is then mapped onto the reference triangle $\wh{\T}$ of vertices $(0,0,0)$, $(1,0,0)$, $(0,1,0)$, polar coordinates are introduced on $\wh{\T}$, and an $n\times n$ Gauss--Legendre product quadrature rule is applied with respect to the radial and angular variables. 
    This procedure is described in detail in case $iv)$, as it is part of the quadrature strategy adopted therein.
	
	\
	
	For case $iv)$ the computation of the Cauchy principal value requires special treatment. Although this is a well-established procedure (see, for instance, \cite{Guiggiani}, \cite{MonegatoHyper}), we provide some details for our setting.
	To simplify the notation, we refer to formula \eqref{eq:Dkappa_tangential}, with $\Gamma$ replaced by $\T$. Recalling the expression \eqref{eq:lapl_transf} of $K$ and using the explicit representation of its gradient, the strongly singular integral, first term in \eqref{eq:Dkappa_tangential}, can be written as
	\begin{equation}\label{eq:CPV0}
	\begin{aligned}
	I &= - \frac{1}{4\pi}\int_\T s\frac{e^{-rs}}{r}\,[\bfn\times\nabla_\bfx r]\cdot [\bfn\times\nabla v(\bfy)]\,\text{d}\T_\bfy\\
	&\hspace{2cm} -\frac{1}{4\pi}\int_\T \frac{e^{-rs}}{r^2}\,[\bfn\times\nabla_\bfx r]\cdot [\bfn\times\nabla v(\bfy)]\,\text{d}\T_\bfy =: I_1 - I_2,
	\end{aligned}
	\end{equation}
	where $r = \|\x-\y\|$ and $\bfn = \bfn_\bfx \equiv \bfn_\bfy$ denotes the common unit normal vector at $\T$ and at $\bfx$. 
	Since $\nabla_{\bfx} r = \mathcal{O}(1)$ as $r\to 0$, it follows that $I_1$ is only weakly singular with a singularity $\mathcal{O}(1/r)$, whereas the second term $I_2$ is strongly singular, exhibiting a singularity $\mathcal{O}(1/r^2)$. Therefore, the integral $I_2$ must be interpreted in the Cauchy Principal Value (CPV) sense, as it is shown next.
	
	By adding and subtracting $1$ in the numerator of the quotient appearing in $I_2$, we rewrite it as follows
	\begin{equation}\label{eq:CPV1}
	\begin{aligned}
	I_2 &= \frac{1}{4\pi}\int_\T \frac{e^{-rs}-1}{r^2}\,[\bfn\times\nabla_\bfx r]\cdot [\bfn\times\nabla v(\bfy)]\,\text{d}\T_\bfy \\[5pt]
	&\hspace{2cm}+ \frac{1}{4\pi}\int_\T \frac{1}{r^2}\,[\bfn\times\nabla_\bfx r]\cdot [\bfn\times\nabla v(\bfy)]\,\text{d}\T_\bfy=: I_{21}+I_{22}.
	\end{aligned}
	\end{equation}
	We note that $I_{21}$ is weakly singular with singularity $\mathcal{O}(1/r)$, so that we  discuss on $I_{22}$ which has to be defined in the CPV sense. Following \cite{MonegatoHyper}, let $\Cex$ denote the disk centred at $\x$ with radius $\epsilon$
	and rewrite
	\begin{equation*}
	I_{22}
	=\frac{1}{4\pi}\lim_{\epsilon\to 0^+}
	\int_{\T\setminus\Cex}
	\frac{1}{r^2}
	\,[\bfn\times\nabla_{\bfx}r]
	\cdot
	[\bfn\times\nabla v(\bfy)]
	\,\text{d}\T_\bfy.
	\end{equation*}
	\begin{figure}[h]
		\centering
		\includegraphics[width=12.0cm,height=5.0cm]{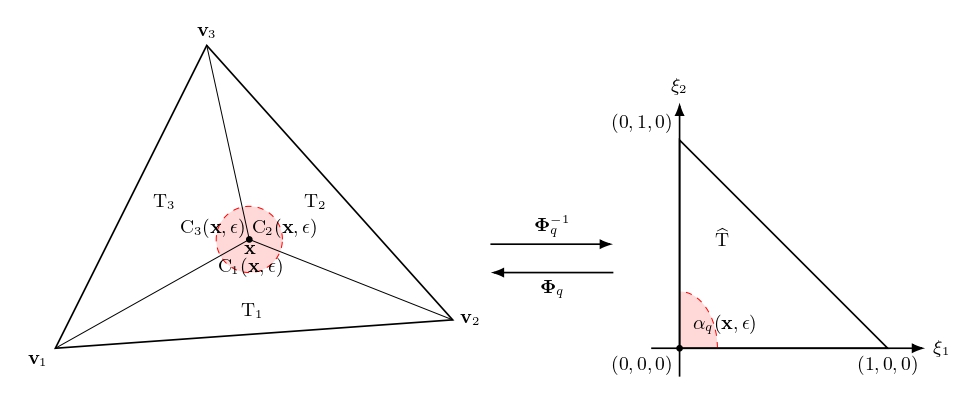}
		\caption{\footnotesize{The integration domain \(\T\setminus\Cex\) with its sub-domains  \(\T_q\setminus\Cexq\), $q = 1,2,3$ (left plot), each mapped to the reference one $\wh{\T}\setminus\alpha_q(\x,\epsilon)$ by \(\bfPhi^{-1}_q\) (right plot).}}\label{fig:Triangle_CPV}
	\end{figure}
	
	Let $\bfv_1, \bfv_2, \bfv_3$ be the vertices of $\T$. By splitting the integration domain $\T\setminus\Cex$ into three subregions determined by its intersection with the three subtriangles $\T_q$ of vertices $\x, \bfv_q, \bfv_{q+1}$, $q = 1,2,3$ (with the convention $\bfv_4\equiv \bfv_1$), we get  
	\begin{equation*}
	I_{22}
	=\frac{1}{4\pi}\lim_{\epsilon\to0^+}
	\sum_{q=1}^3\int_{\T_q\setminus\Cexq}
	\frac{1}{r^2}
	\,[\bfn\times\nabla_{\bfx}r]
	\cdot
	[\bfn\times\nabla v(\bfy)]
	\,\text{d}\T_\bfy,
	\end{equation*}
	where $\Cexq := \Cex\cap \T_q$ (see Figure \ref{fig:Triangle_CPV}).
			
	We map each subtriangle $\mathrm{T}_q$ onto the reference triangle $\widehat{\mathrm{T}}$,
	by means of the inverse of the affine
	transformation $\boldsymbol{\Phi}_q\colon \widehat{\mathrm{T}} \to \mathrm{T}_q$ defined as
	\[
	\boldsymbol{\Phi}_q(\boldsymbol{\xi}) = \mathbf{B}_q \boldsymbol{\xi} + \mathbf{x},
	\qquad \boldsymbol{\xi} = [\xi_1, \xi_2, 0]^T,
	\]
	where
	\[
	\mathbf{B}_q = \big[\, \mathbf{v}_q - \mathbf{x},\ \mathbf{v}_{q+1} - \mathbf{x},\ \mathbf{n} \,\big],
	\qquad q = 1,2,3.
	\]
	Since $\mathbf{x}$ lies in the interior of $\mathrm{T}$, the first two columns of
	$\mathbf{B}_q$ are linearly independent, so that $\mathbf{B}_q$ is invertible.
	Hence, $I_{22}$ can be rewritten as
	\[
	I_{22} = \lim_{\epsilon \to 0^+} \sum_{q=1}^{3} |\det(\mathbf{B}_q)|
	\int_{\widehat{\mathrm{T}} \setminus \alpha_q(\mathbf{x},\epsilon)}
	\frac{\mathrm{dot}_q(\mathbf{x}, \boldsymbol{\xi})}{\|\mathbf{x} - \boldsymbol{\Phi}_q(\boldsymbol{\xi})\|^2}
	\, d\xi_1\, d\xi_2 ,
	\]	
	where we have shortly denoted
	$$
	\operatorname{dot}_q(\bfx,\boldsymbol{\xi}) = \frac{1}{4\pi}\left[\bfn\times\nabla_{\bfx}\|\x-\bfPhi_q(\boldsymbol{\xi})\|\right]\cdot\left[\bfn\times\nabla v(\bfPhi_q(\boldsymbol{\xi})))\right]
	$$
	and
	$$
	\alpha_q(\x,\epsilon)=\bfPhi_q^{-1}(\Cexq).
	$$
	Introducing the polar coordinates $(\rho,\vartheta)$ centred at $(0,0)$ such that $\xi_1 = \rho\cos\vartheta$ and $\xi_2 = \rho\sin\vartheta$, with $0\leq \vartheta\leq \pi/2$ and $0\leq \rho\leq 1 /(\sin\vartheta+\cos\vartheta)$, we have 
	$$
	\begin{aligned}
	\|\bfx-\bfPhi_q(\boldsymbol{\xi})\|^2 &= \|\bfB_q\boldsymbol{\xi}\|^2= \rho^2\sum_{i=1}^{3}\left((\bfB_q)_{i1}\cos\vartheta + (\bfB_q)_{i2}\sin\vartheta\right)^2 \\
	&=: \rho^2 \sum_{i=1}^{3} f^2_{q,i}(\vartheta)=: \rho^2 f^2_q(\vartheta),
	\end{aligned}
	$$
	$$
	\nabla_\bfx \|\bfx-\bfPhi_q(\boldsymbol{\xi})\| = -\frac{1}{f_q(\vartheta)}\left[f_{q,1}(\vartheta),\,f_{q,2}(\vartheta),\,f_{q,3}(\vartheta)\right]^T,
	$$
	and 
	$$
	\alpha_q(\x,\epsilon)= \left\{(\rho,\vartheta)\colon 0\leq\rho\leq \frac{\epsilon}{f_q(\vartheta)},\ 0\leq \vartheta\leq \frac{\pi}{2}\right\}.
	$$
	Recalling that \(v\) is a linear basis function and hence has a constant gradient, and observing that \(\nabla_{\bfx} r\) depends only on \(\vartheta\), we conclude that 
	the function \(\operatorname{dot}_q(\cdot,\cdot)\) does not depend on $\rho$. By abuse of notation, 
	we denote it by $\operatorname{dot}_q(\x,\vartheta)$, so that
	\begin{equation}\label{eq:CPV2}
	\begin{aligned}
	I_{22} &= \lim_{\epsilon\to 0^+}\sum_{q=1}^3|\det(\bfB_q)| \int_0^{\pi/2} \frac{\operatorname{dot}_q(\x,\vartheta)}{f^2_q(\vartheta)} \,\text{d}\vartheta
	\int^{1/(\cos\vartheta+\sin\vartheta)}_{\epsilon/f_q(\vartheta)} \frac{1}{\rho}\,\text{d}\rho \\
	& \hspace{0.5cm} =  \sum_{q=1}^3|\det(\bfB_q)| \int_0^{\pi/2}\frac{\operatorname{dot}_q(\x,\vartheta)}{f^2_q(\vartheta)}\ln\left(\frac{f_q(\vartheta)}{\cos\vartheta+\sin\vartheta}\right)\,\text{d}\vartheta\\
	& \hspace{1cm}  - \lim_{\epsilon\to 0^+} \ln\epsilon\sum_{q=1}^3|\det(\bfB_q)|\int_0^{\pi/2}\frac{\operatorname{dot}_q(\x,\vartheta)}{f^2_q(\vartheta)}\,\text{d}\vartheta. 
	\end{aligned}
	\end{equation}
	It is worth noting that the original singularity $\mathcal{O}(1/r^2)$ in $I_{22}$
	has been weakened to the milder one of order $\mathcal{O}(1/\rho)$ by the Jacobian of the polar coordinate transformation. 
	
	We remark that  (see \cite[Theorem 3.1]{MonegatoHyper}) a necessary and sufficient condition for the existence of the limit in \eqref{eq:CPV2} is
	\begin{equation}\label{eq:esistenza_CPV}
	\sum_{q=1}^3|\det(\bfB_q)|\int_0^{\pi/2}\frac{\operatorname{dot}_q(\x,\vartheta)}{f^2_q(\vartheta)}\,\text{d}\vartheta=0.
	\end{equation}
		Condition \eqref{eq:esistenza_CPV} has been numerically confirmed and, hence, the CPV of $I_{22}$ is 
	\begin{equation}\label{eq:CPV3}
	I_{22}  =  \sum_{q=1}^3|\det(\bfB_q)| \int_0^{\pi/2}\frac{\operatorname{dot}_q(\x,\vartheta)}{f^2_q(\vartheta)}\ln\left(\frac{f_q(\vartheta)}{\cos\vartheta+\sin\vartheta}\right)\,\text{d}\vartheta,
	\end{equation}
and it is evaluated by means of an \(n\)-point Gauss--Legendre quadrature rule.
	
Finally, the weakly singular integrals \(I_1\) in~\eqref{eq:CPV0} and \(I_{21}\) in~\eqref{eq:CPV1} are computed using the same triangle subdivision strategy, without the need to exclude the disc $\Cex$.

\subsection{Stability of the CQ-collocation BEM}\label{sec:stab}

It is well known that, unlike the Galerkin methods, the collocation approaches generally lack a theoretical framework ensuring stability and convergence. In particular, for the problem under consideration, rigorous proofs of stability and convergence are not currently available. 

 To provide insights on the numerical stability of our approach, we preliminarily recall that the accurate computation of the matrix entries of the final linear system, arising from the discretization of the space integrals, is crucial for the accuracy and stability of the numerical procedure, either combined with a collocation or with a Galerkin method in space. Insufficient quadrature accuracy may lead to spurious oscillations in numerical approximation of the density function.
In \cite{Chappell}, theoretical results on stability and convergence are proved for the Galerkin discretization of \eqref{eq:hyper_W2_t_bis}. In particular, the stability analysis relies on the assumption that the space integrals are evaluated either analytically or to machine precision. In practical computations, however, these integrals are approximated with a finite accuracy. Such approximation errors may compromise the stability of the numerical method, depending on the selected time and space discretization parameters.
Consequently, as also stated in \cite{Chappell} and \cite{Lubich2}, excessively small time steps cannot be employed unless the space integrals are evaluated with sufficiently high accuracy. The same phenomenon has been observed for the collocation method proposed here.
 Indeed, extensive numerical experiments indicated that the accurate evaluation of boundary integrals is crucial to preserve numerical stability. In particular, numerical oscillations consistently disappear as the number of cubature nodes is increased. This behaviour has already been observed in our earlier papers (see, e.g., \cite{FMS_Dir}, \cite{FMS_Neu}, \cite{FS_fascette}) and is further supported by the results reported in Example 1 (see Table \ref{table:instab} and Figure \ref{fig:figurastab}). These findings provide strong numerical evidence that, when sufficient quadrature accuracy is employed, the CQ-collocation method provides numerical evidence of unconditional stability.

\section{Numerical results}\label{sec:num_res}

In the numerical experiments, we consider two benchmark spatial domains: the sphere and the cube, representing regular and non-regular geometries, respectively. In both cases, the surface is decomposed by a mesh of $\mathcal{N}_\Gamma$ planar triangular elements. The mesh refinement is performed by dividing each triangle into four sub-triangles, obtained by connecting the midpoints of its edges, effectively halving the mesh diameter $h$ at each refinement level. In Table~\ref{tab:AB}, we report the mesh size $h$, together with the corresponding numbers of triangles $\mathcal{N}_\Gamma$ and vertices $M$, for the refinement levels used in the numerical experiments, referring respectively to the unit sphere and to the cube of side length 1.

Similarly, the time interval $[0,T]$ is uniformly partitioned with time step size $\Delta_t = T/N$, which is also halved at each temporal refinement by doubling the number $N$ of the subintervals. 

\begin{table}[h!]
	\centering
	\footnotesize
		\begin{tabular}{cccccc}
			\hline
			\multicolumn{3}{c}{sphere} & \multicolumn{3}{c}{cube} \\
			\hline
			$h$ & $\mathcal{N}_\Gamma$ & $M$ & $h$ & $\mathcal{N}_\Gamma$ & $M$ \\
			\hline
			$1.05\times10^{0\phantom{-}}$ & 20 & 12 &  --- & --- & --- \\
			$5.43\times10^{-1}$ & 80 & 42 & $5.00\times10^{-1}$ & 96 & 50 \\
			$2.76\times10^{-1}$ & 320 & 162 & $2.50\times10^{-1}$ & 384 & 194 \\
			$1.38\times10^{-1}$ & 1280 & 642 & $1.25\times10^{-1}$ & 1536 & 770 \\
			$6.92\times10^{-2}$ & 5120 & 2562 & $6.25\times10^{-2}$ & 6144 & 3074  \\
			\hline
		\end{tabular}
		\caption{\footnotesize{Mesh-size, number of triangles and vertices for the unit sphere and the cube of side length 1.}}
		\label{tab:AB}
\end{table}

\

To assess the efficiency of the proposed approach, we apply it to some benchmark examples. Recalling that theoretical convergence results for the CQ-collocation method are currently unavailable, we refer to the known ones of the corresponding CQ-Galerkin method proved in~\cite{Chappell}.
For the CQ-collocation approach, we expect a convergence behaviour similar to that of the CQ-Galerkin one. 
To validate this conjecture, recalling that the convolution quadrature is based on a BDF2 scheme, we assume that the problem datum satisfies ${g} \in C^\ell([0, T], C(\mathbb{R}^3))$
for some integer $\ell \geq 4$, and that the following compatibility conditions
$(d^{m}/dt^m){g}(\cdot, t)_{|_{t=0}} = 0$ for $m = 0,\ldots,3$ hold.
In this setting we recall that, when using piecewise linear approximants in space, the convergence rates proved in~\cite{Chappell} for the Galerkin approach are
\begin{equation*}
\| \varphi(\cdot, t_n) - \varphi^{(n)}_{\Delta_t, h} \|_{H^{1/2}(\Gamma)} 
= \mathcal{O}(\Delta_t^2) + \mathcal{O}(h^{3/2}), 
\end{equation*}
\begin{equation*}
\| u(\cdot, t_n) - u^{(n)}_{\Delta_t, h} \|_{H^{1}(\Omegam)} 
= \mathcal{O}(\Delta_t^2) + \mathcal{O}(h^{3/2}), 
\end{equation*}
uniformly with respect to $t_n \in [0, T]$.
Therefore, taking into account known inequalities for Sobolev norms (see, for example, \cite{GRAHAM_HACK_SAUTER}), it appears reasonable to
conjecture, for the $L^2$-norm and the collocation approach, convergence rates of the type
\begin{equation}\label{eq:est_phi}
\| \varphi(\cdot, t_n) - \varphi^{(n)}_{\Delta_t, h} \|_{L^2(\Gamma)} 
= \mathcal{O}(\Delta_t^2) + \mathcal{O}(h^2),
\end{equation}
\begin{equation}\label{eq:est_u}
\| u(\cdot, t_n) - u^{(n)}_{\Delta_t, h} \|_{L^2(\Omegam)} 
= \mathcal{O}(\Delta_t^2) + \mathcal{O}(h^{5/2}).
\end{equation}
These bounds suggest taking in both cases $\Delta_t$ and $h$ such that $\Delta_t \approx h$, which corresponds to the choices of the discretization parameters we have done in the numerical examples.

\

In all the examples, for the approximation of the space integrals described in Section \ref{sec:cub}, unless otherwise specified, the cubature rule in cases $i)$ and $iii)$ has exactness order $\delta = 5$, whereas $n = 8$ Gauss--Legendre quadrature nodes are employed in case $iv)$ and $\delta = 9$ in case $iii)$.

\

All the numerical computations have been performed on a PC with Intel\textsuperscript{\textregistered} Core\textsuperscript{ TM} i7-7700 (3.60 GHz) using parallel Matlab\textsuperscript{\textregistered} codes.

\

\textbf{Example 1.} Let \(\Omega^i\) denote the unit sphere, and consider the Neumann boundary datum on its surface given by 
\[g(\x,t)=\tilde{g}(t), \qquad \text{with} \quad \tilde{g}(t) = t^3\sin(4t)e^{-2t}, \quad  t\in [0,4].\]
In this particular setting, explicit analytical expressions for both \(\varphi\) and \(u\), solutions of equations \eqref{eq:BIE} and~\eqref{pb:esterno}, respectively, are available.
In particular, the density function is defined as (see \cite{Veitphd})
\begin{equation}\label{eq:exactsol_veit}
\begin{aligned}
\varphi(t) &= -2\int_0^t \tilde{g}(t-\tau)\cosh(\tau) \\
&\hskip0.35cm +2\sum_{k=1}^{\lfloor t/2\rfloor}(-1)^{k+1}\sum_{l=1}^{k}\int_{2k}^t c_{k,l}(\tau-2k)^{k-l+1}e^{\tau-2k}\tilde{g}'(t-\tau)\,d\tau
\end{aligned}
\end{equation}
with $c_{k,l}=\binom{k-1}{l-1}\frac{2^{k-l}}{(k-l+1)!}$. Note that this solution depends only on the time variable. The solution $u$ is given by (see Section \ref{sec:anal_sol} of Appendix for details) 
\begin{equation}\label{eq:uesatta_veit}
u(\mathbf{x}, t) = -\frac{1}{r} \int_{r - 1}^t \tilde{g}(t - \tau) \, e^{-(\tau - (r - 1))} \, d\tau, \quad r = \|\x\|
\end{equation}
and has a radial behaviour.
Taking into account the space independence property of $\varphi$, we analyze the error only with respect to the time variable. In particular, since estimates \eqref{eq:est_phi} and \eqref{eq:est_u} are uniform with respect to $t_n \in[0, T]$, the error is calculated using the following mean square formula
\begin{equation}\label{eq:errorephi}
\operatorname{E}^\varphi_{\Delta_t,h}(\bar{\x})= \left(\Delta_t\sum_{n=0}^{N}|\varphi(\bar{\x},t_n) - \varphi^{(n)}_{\Delta_t,h}(\bar{\x})|^2\right)^{1/2}
\end{equation}
for an arbitrarily chosen point $\bar{\x}\in\Gamma$.
The approximate solution $u^{(n)}_{\Delta_t,h}$ is then computed at a point $\mathbf{P}\in\Omegam$ by \eqref{eq:u_rec}, and the corresponding error is analogously calculated using 
\begin{equation}\label{eq:erroreu}
\operatorname{E}^u_{\Delta_t,h}(\mathbf{P}) = \left(\Delta_t\sum_{n=0}^{N}|u(\mathbf{P},t_n) - u^{(n)}_{\Delta_t,h}(\mathbf{P})|^2\right)^{1/2}.
\end{equation}
We apply here the shifted collocation method and in Table \ref{t:table2} we report the errors \eqref{eq:errorephi} and \eqref{eq:erroreu} at $\bar\x=\tilde{\x}/\|\tilde{\x}\|$, with $\tilde{\x}=(-1,\,(1+\sqrt{5})/2,\,0)$ and  $\p = (0,0,1.5)$, respectively.
Moreover, we report the Estimated Convergence Order (EOC) defined as
\[\operatorname{EOC}^* = \log_2\left( \frac{\operatorname{E}^*_{\Delta_t,h}(\bullet)}{\operatorname{E}^*_{\frac{\Delta_t}{2},\frac{h}{2}}(\bullet)}\right), \qquad * = \varphi, u, \quad \bullet = \bar{\x}, \mathbf{P}\]
and the spectral condition number $\kappa_2(\mathbf{A}^{(0)})$ of the system matrix in \eqref{E:SYS}. The results therein reported confirm that halving both discretization parameters recovers the expected quadratic convergence order, with an optimal condition number.

\setlength{\tabcolsep}{3pt}
\begin{table}[h!]	
	\centering
	\footnotesize
		\begin{tabular}{ccccccc}
			\hline
			\multicolumn{1}{c}{$\Delta_t$} & \multicolumn{1}{c}{$h$} & \multicolumn{1}{c}{$\operatorname{E}^u_{\Delta_t,h}(\p)$} & \multicolumn{1}{c}{EOC$^u$} & \multicolumn{1}{c}{$\operatorname{E}^\varphi_{\Delta_t,h}(\bar{\x})$} & \multicolumn{1}{c}{EOC$^\varphi$} & \multicolumn{1}{c}{$\kappa_2(\mathbf{A}^{(0)})$} \\
			\hline
			$1.25\times 10^{-1}$  & $1.05\times 10^{0\phantom{-}}$ & $1.68\times 10^{-2}$ & \multicolumn{1}{c}{---} & $5.79\times 10^{-2}$ & \multicolumn{1}{c}{---} & $2.66$ \\
			$6.25\times 10^{-2}$  & $5.47\times 10^{-1}$ & $5.39\times 10^{-3}$ & $1.64$ & $1.44\times 10^{-2}$ & $2.01$ & $2.87$ \\
			$3.13\times 10^{-2}$  & $2.76\times 10^{-1}$ & $1.37\times 10^{-3}$ & $1.98$ & $3.94\times 10^{-3}$ & $1.87$ & $2.25$ \\
			$1.56\times 10^{-2}$  & $1.38\times 10^{-1}$ & $3.37\times 10^{-4}$ & $2.02$ & $9.22\times 10^{-4}$ & $2.10$ & $2.57$  \\
			$7.81\times 10^{-3}$  & $6.92\times 10^{-2}$ & $8.36\times 10^{-5}$ & $2.01$ & $2.32\times 10^{-4}$ & $1.99$ & $2.85$ \\
			\hline
		\end{tabular}
		\centering
		\caption{\footnotesize{Example 1. Errors, convergence rates, and condition numbers for the collocation method on the unit sphere with $T=4$.}}
		\label{t:table2}
\end{table}

\begin{figure}[h]
	\centering
	\begin{minipage}{0.45\textwidth}
		\centering
		\includegraphics[width=0.90\textwidth]{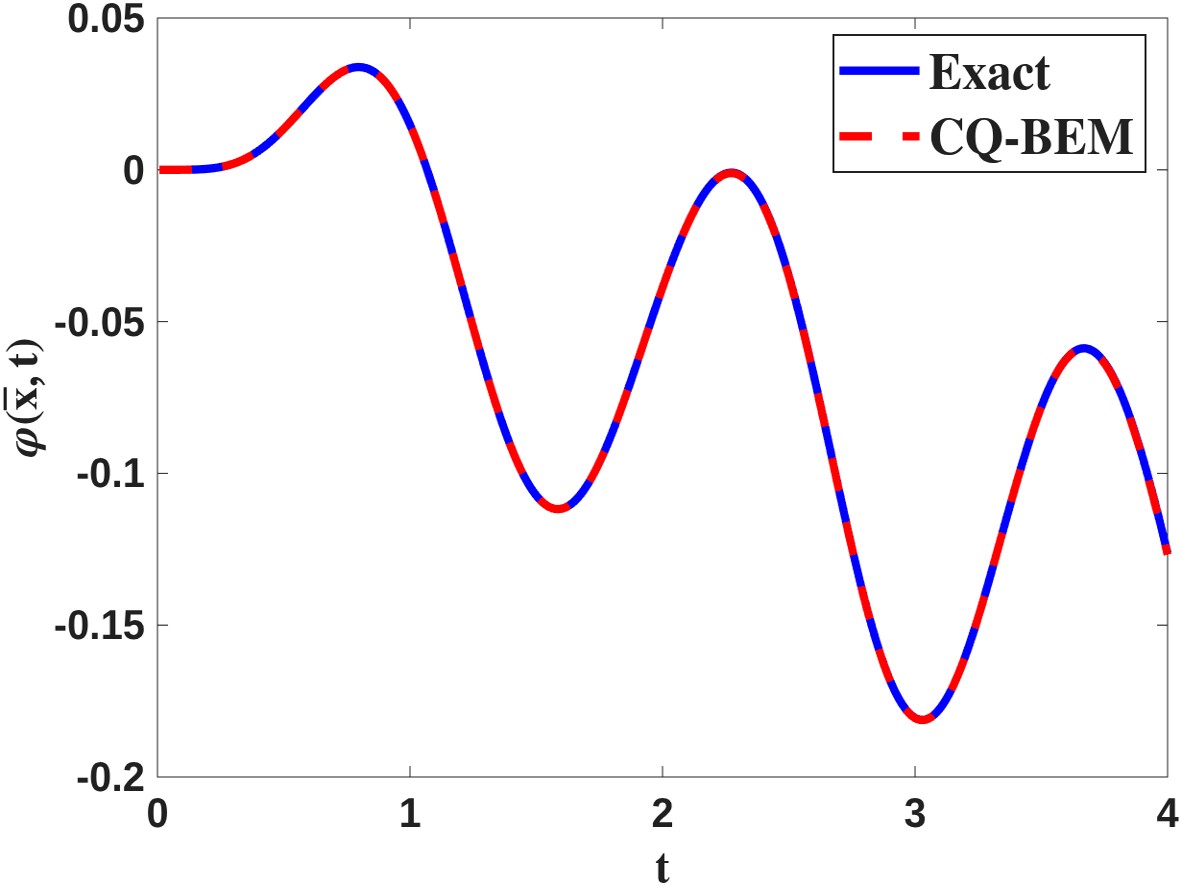}
	\end{minipage}
	\begin{minipage}{0.45\textwidth}
		\centering
		\includegraphics[width=0.90\textwidth]{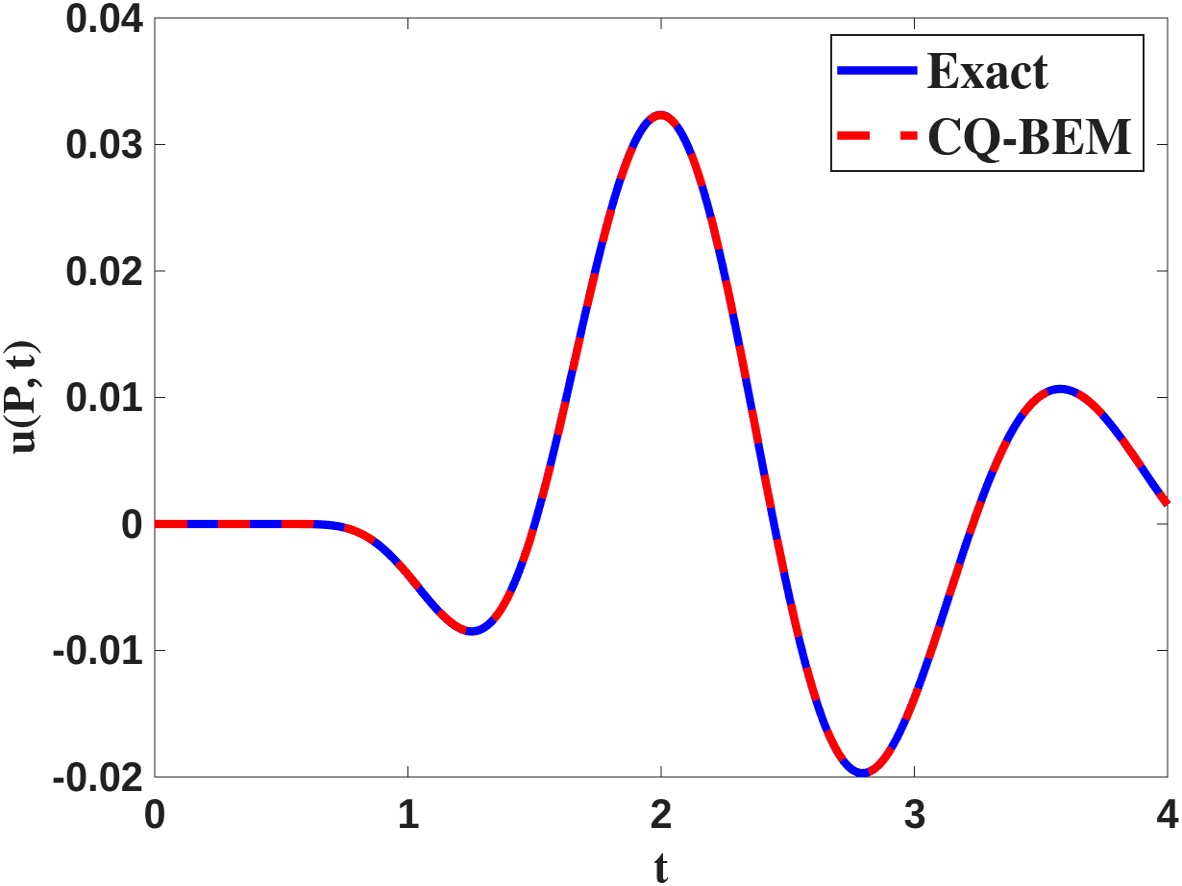}
	\end{minipage}
	\caption{\footnotesize{Example 1. Behaviour of exact and BEM densities (left) and of exact and BEM solutions (right).}}\label{fig:figura1}
\end{figure}

In Figure~\ref{fig:figura1} we compare the behaviour of the exact and BEM-approximated functions for both the density \( \varphi \) and the solution \( u \), which exhibit a very good agreement.

\

\begin{figure}[h]
	\centering
	\begin{minipage}{0.45\textwidth}
		\centering
	\includegraphics[width=0.9\textwidth]{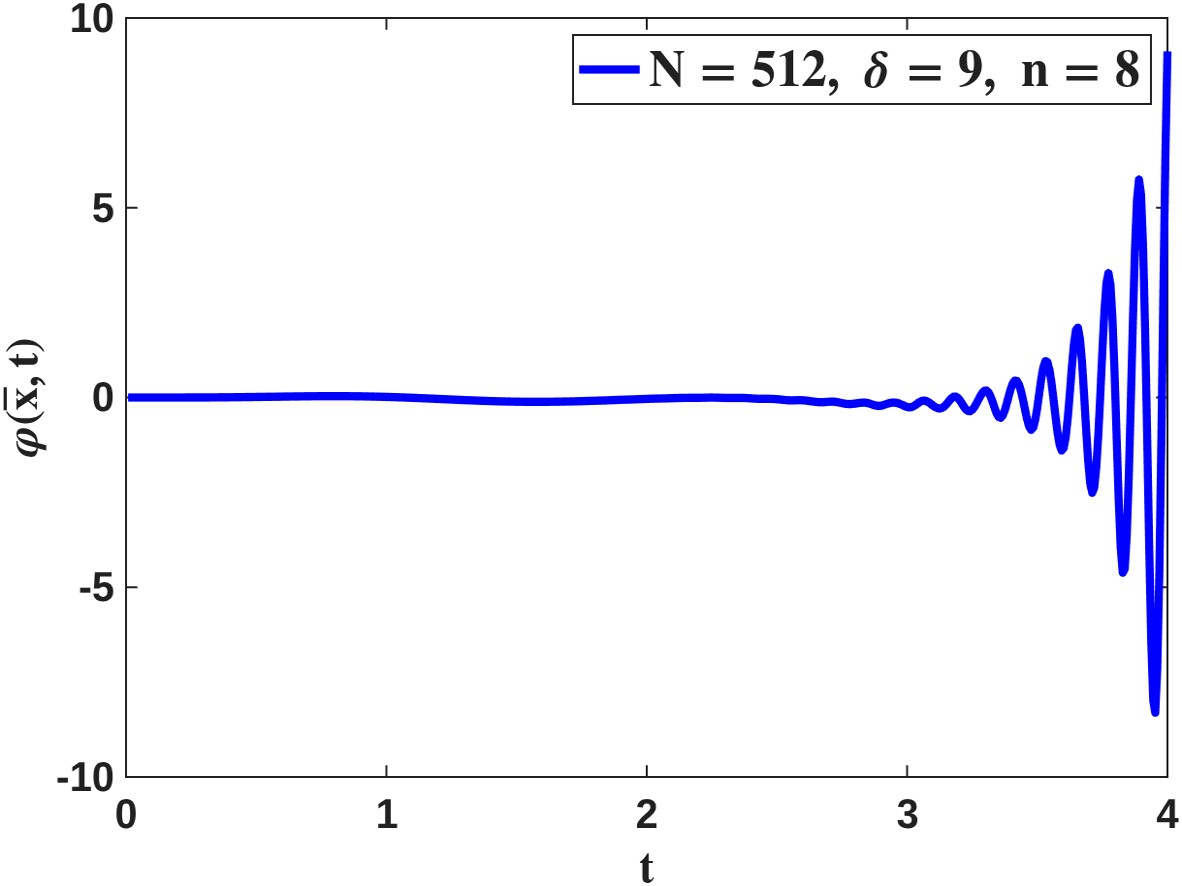}
	\end{minipage}
	\begin{minipage}{0.45\textwidth}
		\centering
		\includegraphics[width=0.9\textwidth]{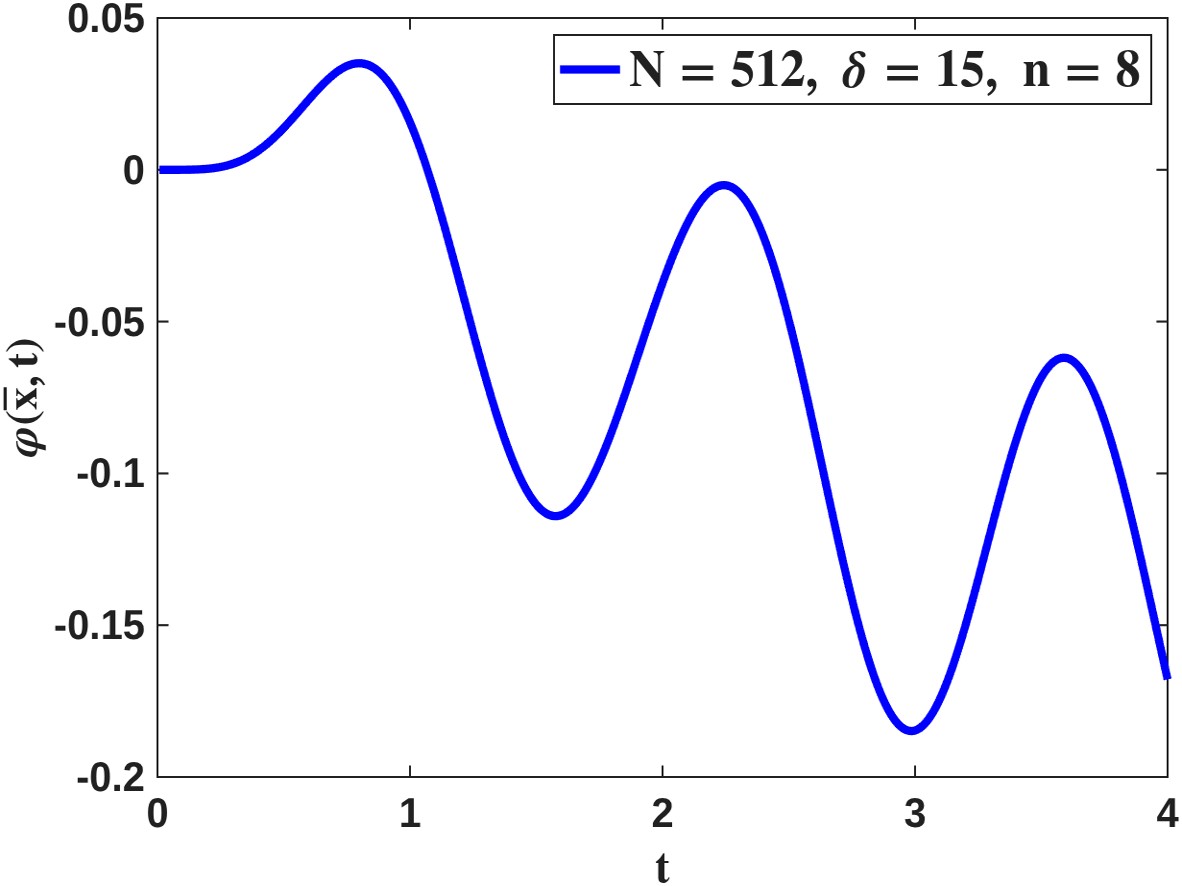}
	\end{minipage}
	\begin{minipage}{0.45\textwidth}
		\centering
		\includegraphics[width=0.9\textwidth]{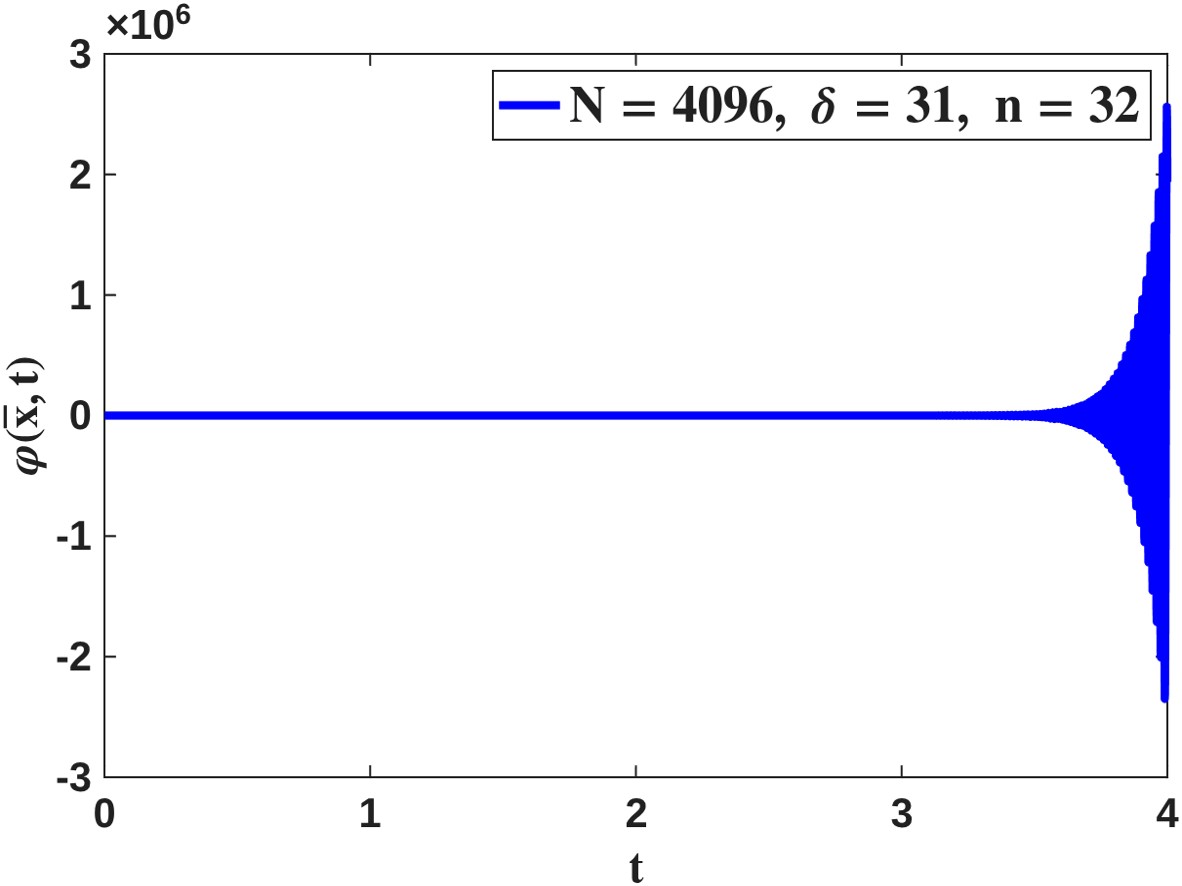}
	\end{minipage}\begin{minipage}{0.45\textwidth}
		\centering
		\includegraphics[width=0.9\textwidth]{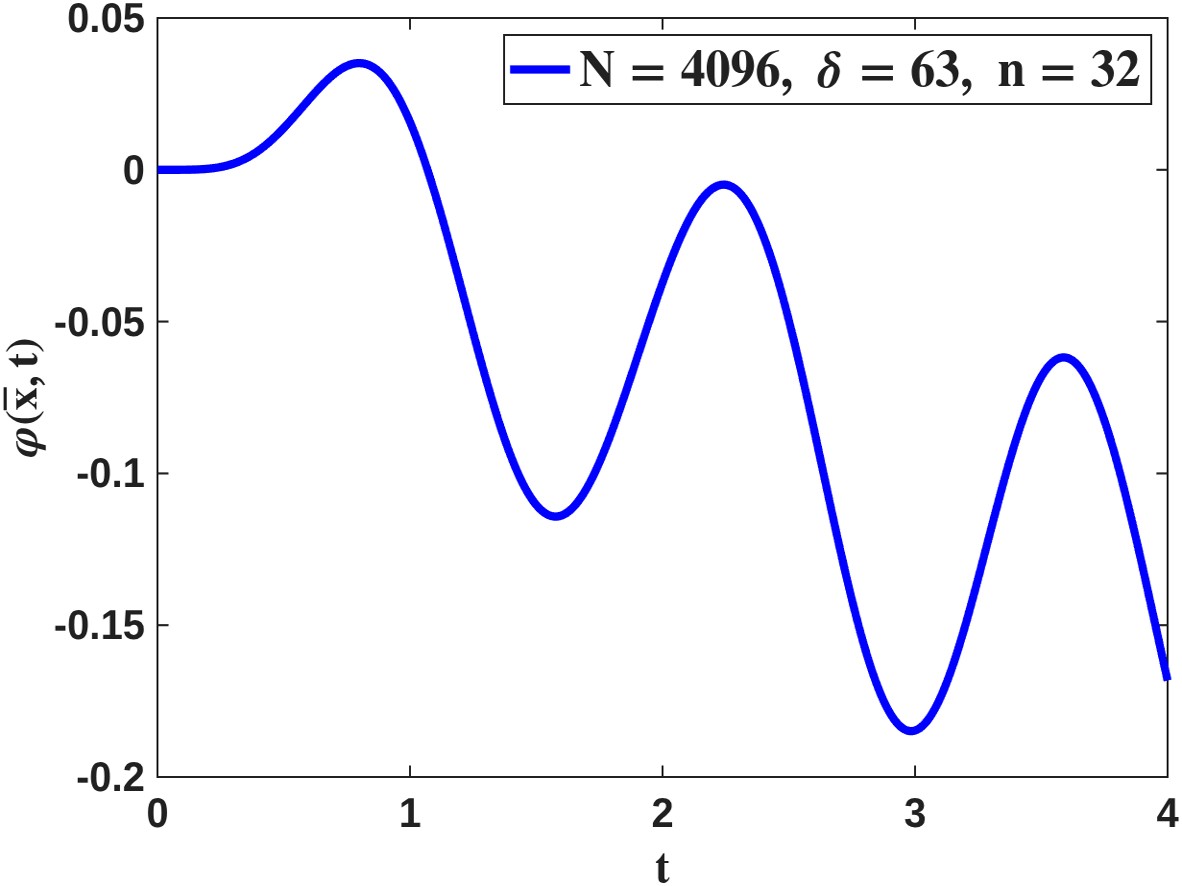}
	\end{minipage}
	\caption{\footnotesize{Example 1. Time behaviour of the BEM density for $h\approx 5.47\times 10^{-1}$ with varying $N$, $\delta$ and $n$.}}\label{fig:figurastab}
\end{figure}

We conclude Example 1 by discussing the numerical stability of the fully discretized scheme. To show that the discretization parameters $h$ and $\Delta_t$ can be chosen independently, and hence that the method is not subject to any CFL-type stability restriction, we fix the spatial mesh size $h\approx 5.47 \times 10^{-1}$, and progressively halve the time step $\Delta_t = T/N$, by doubling $N$. 
	As soon as oscillations or explosions in the numerical approximation appear, we increase the accuracy of the cubature rule. In Table \ref{table:instab} we report the errors of the discrete BEM density function $\operatorname{E}^\varphi_{\Delta_t,h}(\bar{\x})$ defined by \eqref{eq:errorephi}. The results have been obtained by increasing $\delta$ and/or $n$ only in the triangles containing the shifted collocation point and in the triangles sharing at least one vertex with it; otherwise, we have fixed $\delta = 9$. In Figure~\ref{fig:figurastab}, we illustrate the time evolution of the density function. The expected accuracy is recovered and the spurious oscillations disappear once the accuracy of the integral computation is improved. Numerical tests reveal that this effect is due to the localized support of the integrand inside the integration triangle, which becomes smaller as $\Delta_t$ decreases (the same phenomenon has also been observed in \cite{FS_fascette}). Hence, a larger number of quadrature points is required to ensure that sufficiently many nodes lie within the support. This phenomenon does not occur when both $h$ and $\Delta_t$ are reduced.

	In light of the considerations reported in Section~\ref{sec:stab},	provided that the cubature rule used for the evaluation of the spatial integrals is sufficiently accurate,	the numerical results confirm the numerical unconditional stability of the proposed method.

Based on our numerical experiments, no clear relationship between the qua\-dra\-ture parameters $\delta$ and/or $n$ and the time step $\Delta_t$ could be identified. For this reason, whenever a higher integration accuracy was required, we adopted the simple strategy of doubling $\delta$ and/or $n$.
	
	 We remark that a similar numerical behaviour has been observed for the CQ-Galerkin approach, for which the accurate evaluation of the double integrals is equally crucial, in agreement with the analysis in \cite{Chappell}.

\setlength{\tabcolsep}{6pt} 
\begin{table}[h!]
	\centering
	\footnotesize
	\begin{tabular}{ccccc}
		\hline
		\multicolumn{1}{c}{$N$}&\multicolumn{1}{c}{$\Delta_t$} & \multicolumn{1}{c}{$\delta$} &
		\multicolumn{1}{c}{$n$} & \multicolumn{1}{c}{$\operatorname{E}^\varphi_{\Delta_t,h}(\bar{\x})$} \\
		\hline
		$512$ & $7.81\times 10^{-3}$  & $9$ & $8$ & $2.33\times 10^{0\phantom{-}}$ \\
		 &  & $15$ & $8$ & $2.24\times 10^{-2}$ \\
		\hline
		$1024$ & $3.91\times 10^{-3}$ & $15$ & $8$ & $9.50\times 10^{4\phantom{-}}$ \\
		 &  & $15$ & $16$ & $2.25\times 10^{-2}$ \\
		\hline
		$2048$ & $1.95\times 10^{-3}$ & $15$ & $16$ & $9.27\times 10^{6\phantom{-}}$ \\
		 &  & $31$ & $32$ & $2.26\times 10^{-2}$ \\
		\hline
		$4096$ & $9.77\times 10^{-4}$ & $31$ & $32$ & $4.11\times 10^{5\phantom{-}}$ \\
		 &  & $63$ & $32$ & $2.24\times 10^{-2}$ \\
		\hline
		$8192$ & $4.88\times 10^{-4}$ & $63$ & $32$ & $2.72\times 10^{64}$ \\
		 &  & $63$ & $64$ & $2.24\times 10^{-2}$ \\
		\hline
	\end{tabular}
	\centering
	\caption{\footnotesize{Example 1. Errors for fixed $h \approx 5.47\times 10^{-1}$, with varying $\Delta_t$, $\delta$ and $n$.} }
	\label{table:instab}
\end{table}

\

\textbf{Example 2.} 
In this example, we compare the CQ-collocation method proposed here with the CQ-Galerkin approach studied in~\cite{Chappell}. The numerical experiments reported therein are restricted to specific axisymmetric domains, where meshes and quadrature techniques are constructed within a two-dimensional framework rather than a genuinely three-dimensional one, leading to a substantially different numerical treatment. In contrast, we solve the problem in its full three-dimensional setting, without relying on any symmetry. 

As in~\cite{Chappell}, we consider the radiation of a spherically symmetric wave, defined by
\begin{equation}
\label{eq:exactsol}
\begin{aligned}
u(\x,t) = &\frac{H(t+a-r)}{r} \Bigg( \frac{2}{3} 
- \cos\left( \frac{\pi(r - t + 5a)}{3a} \right)
+ \frac{2}{5} \cos\left( \frac{2\pi(r - t + 5a)}{3a} \right) \\
& - \frac{1}{15} \cos\left( \frac{\pi(r - t + 5a)}{a} \right) \Bigg) 
\exp\left( -\left( \frac{r - t + 2a}{3a} \right)^4 \right),
\end{aligned}
\end{equation}
where  $r = \|\mathbf{x}\|$, $\mathbf{x} \in \Omegam$, $t\in [0,T]$, \( H \) denotes the Heaviside step function and, using the same geometric configuration as in Example~1, \( a = 1 \). We observe that the corresponding Neumann boundary datum \( g = \partial u/\partial \mathbf{n} \) inherits the spherical symmetry of \( u \).
We apply the shifted collocation method and, as in \cite{Chappell}, we fix the final time \( T = 16 \) and we compute only the pointwise errors \( \operatorname{E}^u_{\Delta_t,h}(\p) \) at \( \p = (0, 0, 3) \), the exact density $\varphi$ being not analytically known. The corresponding results are reported in Table~\ref{t:table5}. 
To facilitate a meaningful comparison with \cite[Table I]{Chappell}, we adopt the same time step values \( \Delta_t \), and vary the level refinement of the sphere such that \( h \approx \Delta_t \). Although the parameter \( h \) is not directly comparable to the mesh size associated with the one-dimensional generating curve employed in \cite{Chappell}, under our configuration the resulting errors of the CQ-collocation method are very similar to those obtained therein for the CQ-Galerkin one, and the observed convergence rates confirm the expected optimal quadratic behaviour.

\begin{table}[h!]
	\centering	
	\footnotesize
	\begin{tabular}{cccccc}
		\hline
		\multicolumn{1}{c}{$\Delta_t$} & \multicolumn{1}{c}{$h$} & \multicolumn{1}{c}{$\operatorname{E}^u_{\Delta_t,h}(\p)$} & \multicolumn{1}{c}{EOC$^u$} & \multicolumn{1}{c}{$\kappa_2(\mathbf{A}^{(0)})$} &
		\multicolumn{1}{c}{CPU(s)} \\
		\hline
		$5.00\times 10^{-1}$  & $1.05\times 10^{0\phantom{-}}$     & $2.62\times 10^{-1}$ & \multicolumn{1}{c}{---} & $2.44$ & $1.83\times 10^{0\phantom{-}}$ \\
		$2.50\times 10^{-1}$ & $5.47\times 10^{-1}$ & $1.02\times 10^{-1}$ & $1.36$ & $2.49$ & $3.90\times 10^{0\phantom{-}}$ \\
		$1.25\times 10^{-1}$ & $2.76\times 10^{-1}$ & $2.93\times 10^{-2}$ & $1.80$ & $1.88$ & $1.54\times 10^{1\phantom{-}}$ \\
		$6.25\times 10^{-2}$ & $1.38\times 10^{-1}$ & $7.48\times 10^{-3}$ & $1.97$ & $2.23$ & $3.15\times 10^{2\phantom{-}}$ \\
		$3.13\times 10^{-2}$ & $6.92\times 10^{-2}$ & $1.88\times 10^{-3}$ & $1.99$ & $2.52$ & $5.95\times 10^{3\phantom{-}}$ \\
		\hline
	\end{tabular}
	\caption{\footnotesize{Example 2. Errors, convergence rates, condition numbers and CPU times for the collocation method on the unit sphere with $T=16$.}}
	\label{t:table5}
\end{table}

\begin{figure}[h!]
	\centering
	\begin{minipage}{0.48\textwidth}
		\centering
		\includegraphics[width=5.0cm,height=5.0cm]{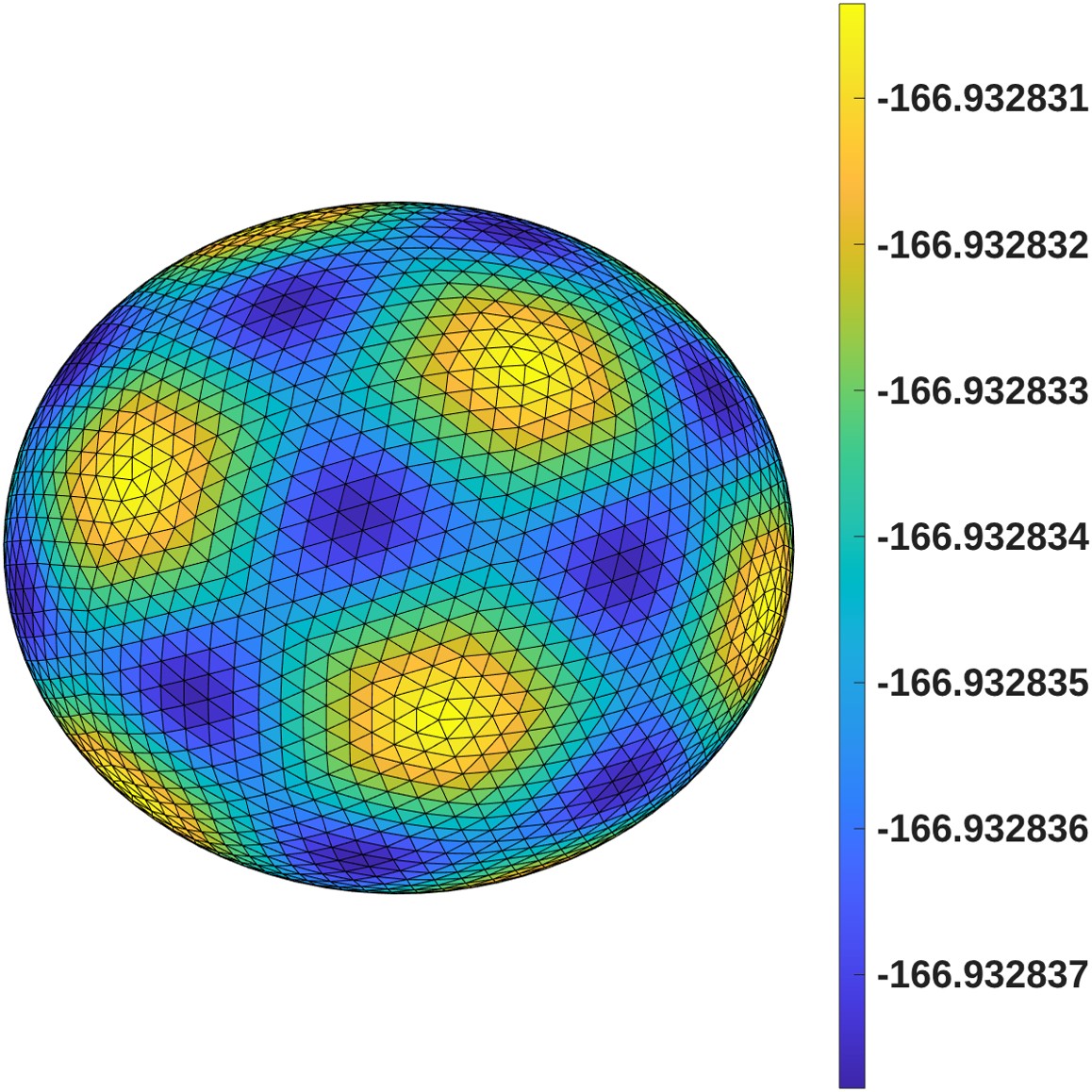}
	\end{minipage}\hskip0.8cm
	\begin{minipage}{0.45\textwidth}
		\centering
		\includegraphics[width=6cm]{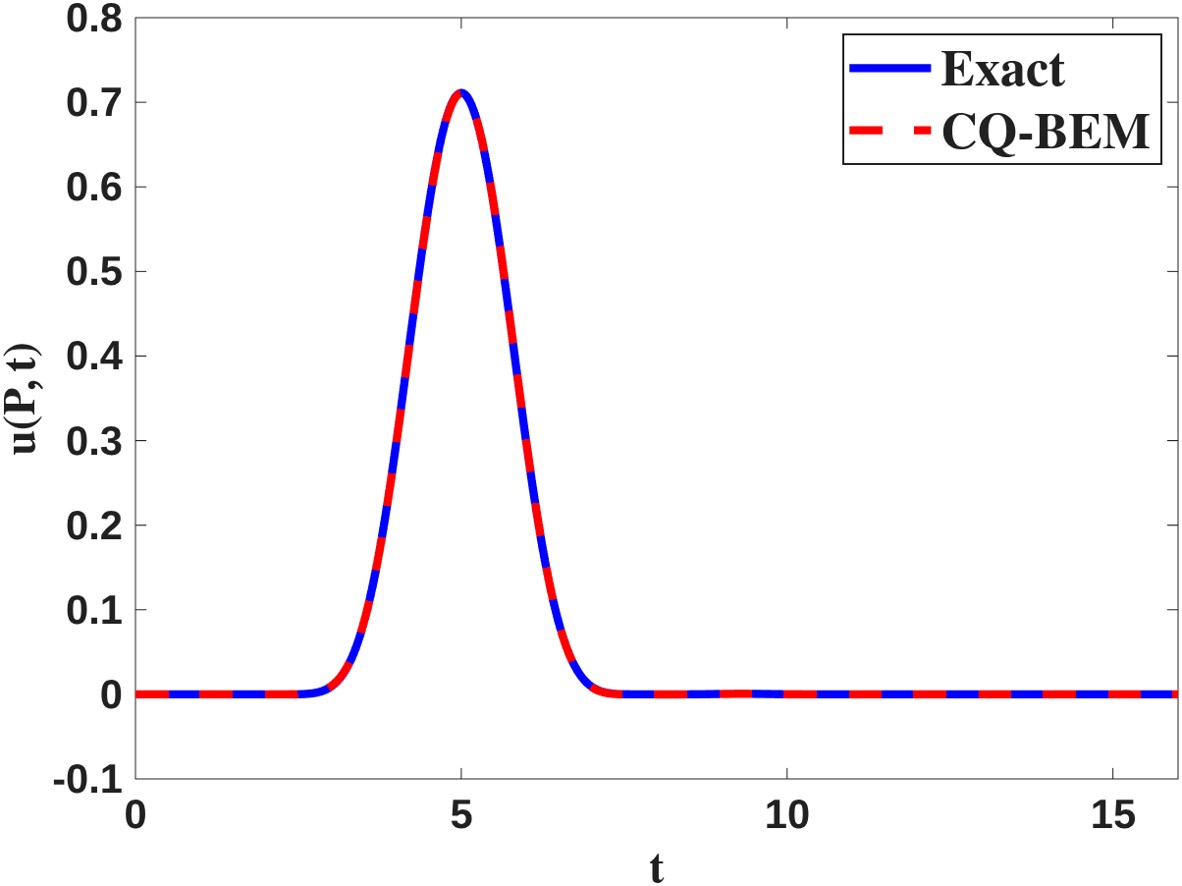}
	\end{minipage}
	\caption{\footnotesize{Example 2. Behaviour of the density $\varphi_{\Delta_t,h}^{(N)}$ on $\Gamma$ (left) and time history of $u_{\Delta_t,h}(\p,t)$ (right).}}\label{fig:sfera}
\end{figure}

\
The left plot of Figure~\ref{fig:sfera} illustrates the computed density \( \varphi_{\Delta_t,h}^{(N)} \) over the surface \( \Gamma_\Delta \) at the final time $t_N$; the right one depicts the numerical approximation \( u_{\Delta_t,h}(\p,t) \) in comparison with the exact solution \( u(\p,t) \) at the fixed spatial point \( \p  \), across the full time range \( t \in [0,16] \).

For completeness, we also implement the CQ-Galerkin method in the three-dimensional setting. In this case, we consider the weak formulation of the boundary integral equation \eqref{eq:hyper_W2}, given by
\begin{align}
&\int_\Gamma \frac{\partial}{\partial \n_\x} 
\left[ \int_\Gamma \frac{\partial }{\partial \n_\y}  \left( \int_0^t G(\x-\y, t-\tau)\,\varphi(\y,\tau)\, \mathrm{d}\tau \right) \mathrm{d}\Gamma_\y \right] w(\x)\, \mathrm{d}\Gamma_\x 
\nonumber\\
& \qquad\qquad =\int_\Gamma g(\x,t)\, w(\x)\, \mathrm{d}\Gamma_\x, 
\label{eq:hyper_W2_weak}
\end{align}
where \( w \) denotes a generic test function.
Following a procedure analogous to the collocation method, we apply the Lubich convolution quadrature, and adopt piecewise linear basis functions both for the trial and test spaces. This leads to a fully discrete linear system of the same form as in \eqref{E:SYS}.
In assembling the system matrices, the integral singularity is weakened by applying integration by parts twice (see \cite[p. 257]{Steinbach2013}). Hence, we consider the following representation, alternative to the left side of \eqref{eq:hyper_W2_weak},
\begin{equation}\label{eq:W_galerk}
\begin{aligned}
&\int_\Gamma \frac{\partial}{\partial \n_\x} \Bigg( 
\int_\Gamma \frac{\partial}{\partial \n_\y} K(r, s)\,v(\y)\, \mathrm{d}\Gamma_\y \Bigg) 
w(\x)\, \mathrm{d}\Gamma_\x \\
& \qquad\qquad = - \int_\Gamma \int_\Gamma 
K(r,s) \left( \n_\y \times \nabla v(\y) \right) \cdot 
\left( \n_\x \times \nabla w(\x) \right) \, \mathrm{d}\Gamma_\y \, \mathrm{d}\Gamma_\x  \\
&\qquad\qquad\quad -\, s^2 \int_\Gamma \int_\Gamma 
\left( \n_\y \cdot \n_\x \right) K(r,s)\,v(\y)\,w(\x) \, \mathrm{d}\Gamma_\y \, \mathrm{d}\Gamma_\x,
\end{aligned}
\end{equation}
which involves only weakly singular integrals, since the integrand functions have a weak singularity \( \mathcal{O}(1/r)\).
The results reported in Table~\ref{t:table6} indicate that the collocation and Galerkin methods exhibit very comparable performance in terms of both accuracy and conditioning. We remark, incidentally, that the errors reported in Table~\ref{t:table6} are very similar to those presented in \cite[Table~I]{Chappell}. 
However, it is worth emphasizing that the collocation method achieves comparable accuracy at a significantly lower computational cost, making it a more efficient and competitive approach in practice, despite the presence of strongly singular integrals. This is confirmed by the CPU times (in seconds) reported in Tables~\ref{t:table5} and~\ref{t:table6}, obtained using the same programming environment to ensure a fair comparison of the computational performance of the two methods. The CPU time ratios between the Galerkin and collocation methods, corresponding to the space-time refinements considered in Tables~\ref{t:table5} and~\ref{t:table6}, are approximately \(1.3\), \(2.4\), \(16.7\), \(12.6\), and \(19.9\), respectively. These results clearly demonstrate the superior computational efficiency of the shifted collocation method, especially for finer discretizations.

\setlength{\tabcolsep}{6pt} 
\begin{table}[h]
	\centering	
	\footnotesize
		\begin{tabular}{cccccc}
			\hline
			\multicolumn{1}{c}{$\Delta_t$} & \multicolumn{1}{c}{$h$} & \multicolumn{1}{c}{$\operatorname{E}^u_{\Delta_t,h}(\p)$} & \multicolumn{1}{c}{EOC$^u$} & \multicolumn{1}{c}{$\kappa_2(\mathbf{A}^{(0)}$)} & \multicolumn{1}{c}{CPU(s)}\\
			\hline
			
			$5.00\times 10^{-1}$  & $1.05\times 10^{0\phantom{-}}$ & $2.53\times 10^{-1}$ & \multicolumn{1}{c}{---}  & $2.08$ & $2.33\times 10^{0\phantom-}$ \\
			$2.50\times 10^{-1}$  & $5.47\times 10^{-1}$  & $9.89\times 10^{-2}$ & $1.36$ & $2.46$ & $9.50\times 10^{0\phantom{-}}$ \\
			$1.25\times 10^{-1}$ & $2.76\times 10^{-1}$  & $2.78\times 10^{-2}$ & $1.83$ & $2.68$ & $2.57\times 10^{2\phantom{-}}$ \\
			$6.25\times 10^{-2}$ & $1.38\times 10^{-1}$ & $7.10\times 10^{-3}$ & $1.97$ & $2.83$ & $3.96\times 10^{3\phantom{-}}$ \\
			$3.13\times 10^{-2}$ & $6.92\times 10^{-2}$ & $1.91\times 10^{-3}$ & $1.89$ & $3.06$ & $1.18\times 10^{5\phantom{-}}$ \\
			\hline
		\end{tabular}
		\centering
		\caption{\footnotesize{Example 2. Errors, convergence rates, condition numbers and CPU times for the Galerkin method on the unit sphere with $T=16$.}}
		\label{t:table6}
\end{table}

We conclude this example with the additional investigation on the performance of the collocation method over longer time intervals. 
 Specifically, we apply the shifted collocation method to the same problem with a final time \( T = 50 \). The numerical results reported in Table~\ref{t:table7} are obtained with the same choices of the discretization parameters \( \Delta_t \) (hence with a larger $N$) and $h$ of Table \ref{t:table5}. The expected convergence behaviour is maintained, demonstrating the robustness of the scheme even over extended time intervals.

\begin{table}[h!]
	\centering	
	\footnotesize
		\begin{tabular}{ccccc}
			\hline
			\multicolumn{1}{c}{$\Delta_t$} & \multicolumn{1}{c}{$h$} & \multicolumn{1}{c}{$\operatorname{E}^u_{\Delta_t,h}(\p)$} & \multicolumn{1}{c}{EOC$^u$} & \multicolumn{1}{c}{$\kappa_2(\mathbf{A}^{(0)})$} \\
			\hline
			$5.00\times 10^{-1}$  & $1.05\times 10^{0\phantom{-}}$     & $2.62\times 10^{-1}$ & \multicolumn{1}{c}{---} & $2.44$  \\
			$2.50\times 10^{-1}$ & $5.47\times 10^{-1}$ & $1.02\times 10^{-1}$ & $1.36$ & $2.49$  \\
			$1.25\times 10^{-1}$ & $2.76\times 10^{-1}$ & $2.93\times 10^{-2}$ & $1.80$ & $1.88$ \\
			$6.25\times 10^{-2}$ & $1.38\times 10^{-1}$ & $7.48\times 10^{-3}$ & $1.97$ & $2.23$  \\
			$3.13\times 10^{-2}$ & $6.92\times 10^{-2}$ & $1.88\times 10^{-3}$ & $1.99$ & $2.52$  \\
			\hline
		\end{tabular}
		\caption{\footnotesize{Example 2. Errors, convergence rates and condition numbers for the collocation method on the unit sphere with $T=50$.}}	
		\label{t:table7}
\end{table}

\

	\textbf{Example 3.} In this example, we again consider the unit sphere and investigate the following solution of Problem~\eqref{pb:esterno}:
	\[
	u(\x,t)=-\frac{x_3}{r}\left(\frac{f'(t-r+1)}{r}+\frac{f(t-r+1)}{r^2}\right), \qquad r = \|\x\|
	\]
	where $\x=(x_1,x_2,x_3)\in\Omega$, and $t\in[0,8]$. The causal impulse $f$ is defined as
	\[
	f(\tau)=
	\begin{cases}
	\tau^4 e^{-2\tau}, & \text{if } \tau>0,\\
	0, & \text{otherwise}.
	\end{cases}
	\]
	Differently from the previous examples, the solution considered here exhibits a non-trivial dependence on the spatial variable, by the factor $x_3/r$. Hence, although the solution still represents a radiating wave, its angular distribution is no longer uniform over the boundary. 
		
	 In Table~\ref{table:9} we report the error  $\operatorname{E}^u_{\Delta_t,h}(\mathbf{P})$ defined in \eqref{eq:erroreu}, for $\p=(0,0,3)$. The results confirm the expected convergence behaviour of the method, with an initial decrease of the EOC due to pre-asymptotic effects on the relatively coarse mesh. As the discretization is refined, the asymptotic convergence regime is reached, leading to the observed increase of the EOC.
	 
\begin{figure}[h!]
	\centering
	\begin{minipage}{0.48\textwidth}
		\centering
		\includegraphics[width=5.5cm,height=5.5cm]{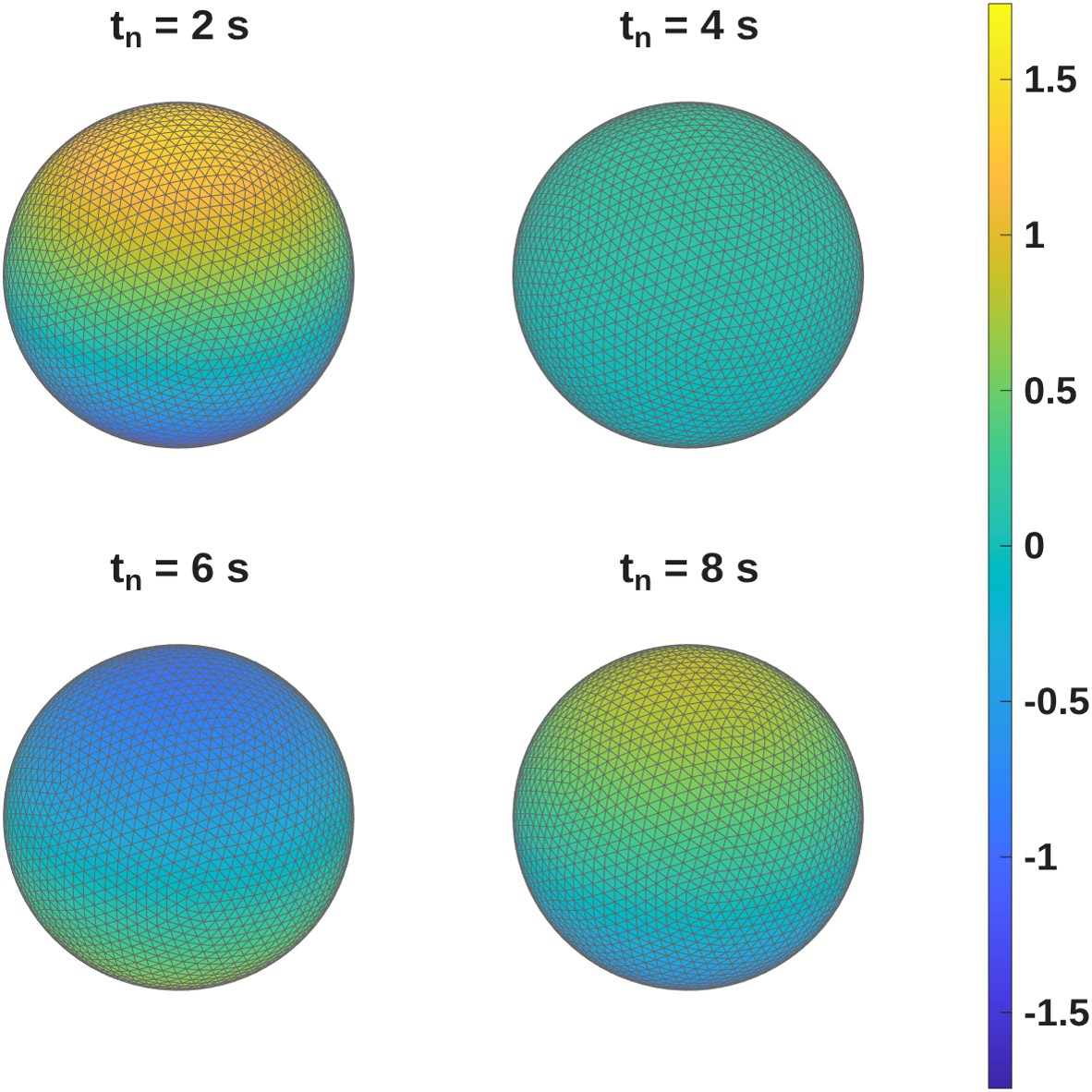}
	\end{minipage}\hskip0.3cm
	\begin{minipage}{0.45\textwidth}
		\centering
		\includegraphics[width=5.8cm]{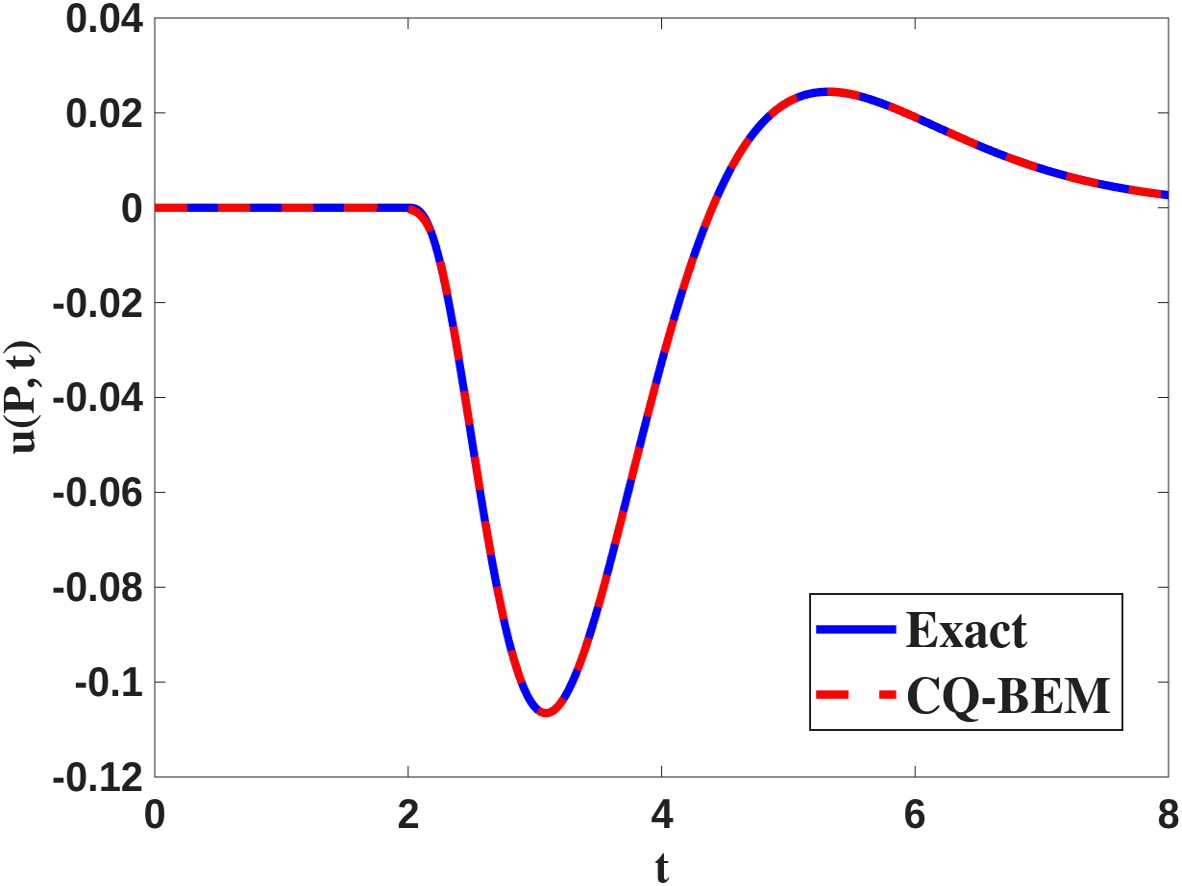}
	\end{minipage}
	\caption{\footnotesize{Example 3. Behaviour of the density $\varphi_{\Delta_t,h}^{(n)}$ on $\Gamma$ at different time instants $t_n$ (left) and time history of $u_{\Delta_t,h}(\p,t)$ (right).}}\label{fig:sfera_ex3}
\end{figure}

	In Figure \ref{fig:sfera_ex3} we show the evolution of the density function $\varphi_{\Delta_t,h}^{(n)}$ on $\Gamma$ at different time instants $t_n$ (left plot), together with the corresponding time history of the solution $u_{\Delta_t,h}(\mathbf{P},t)$ (right plot). For the latter we note that the numerical solution matches the exact one and, in agreement with the causality property of the impulse $f$, reproduces the vanishing response before the arrival of the wave front, as it remains identically zero for $t<2$.
	
	\setlength{\tabcolsep}{6pt} 
	\begin{table}[h]
		\centering	
		\footnotesize
			\begin{tabular}{ccccc}
				\hline
				\multicolumn{1}{c}{$\Delta_t$} & \multicolumn{1}{c}{$h$} & \multicolumn{1}{c}{$\operatorname{E}^u_{\Delta_t,h}(\p)$} & \multicolumn{1}{c}{EOC$^u$} & \multicolumn{1}{c}{$\kappa_2(\mathbf{A}^{(0)}$)} \\
				\hline
				$2.50\times 10^{-1}$& $1.05\times 10^{0\phantom{-}}$ & $5.34\times 10^{-2}$ & \multicolumn{1}{c}{---}  & 2.64\\
				$1.25\times 10^{-1}$& $5.47\times 10^{-1}$ & $2.38\times 10^{-2}$ & $1.17$  & 2.82\\
				$6.25\times 10^{-2}$& $2.76\times 10^{-1}$  & $4.33\times 10^{-3}$ & $2.46$ & $2.17$\\
				$3.13\times 10^{-2}$& $1.38\times 10^{-1}$  & $1.43\times 10^{-3}$ & $1.60$ & $2.51$ \\
				$1.56\times 10^{-2}$& $6.92\times 10^{-2}$ & $3.20\times 10^{-4}$ & $2.16$ & $2.70$ \\
				\hline
			\end{tabular}
			\centering
			\caption{\footnotesize{Example 3. Errors, convergence rates and condition numbers for the collocation method on the unit sphere with $T = 8$.}}
			\label{table:9}
	\end{table}

\

\begin{figure}[h]
	\centering
	\begin{minipage}{0.45\textwidth}
		\centering
		\includegraphics[width=6.0cm,height=5.0cm]{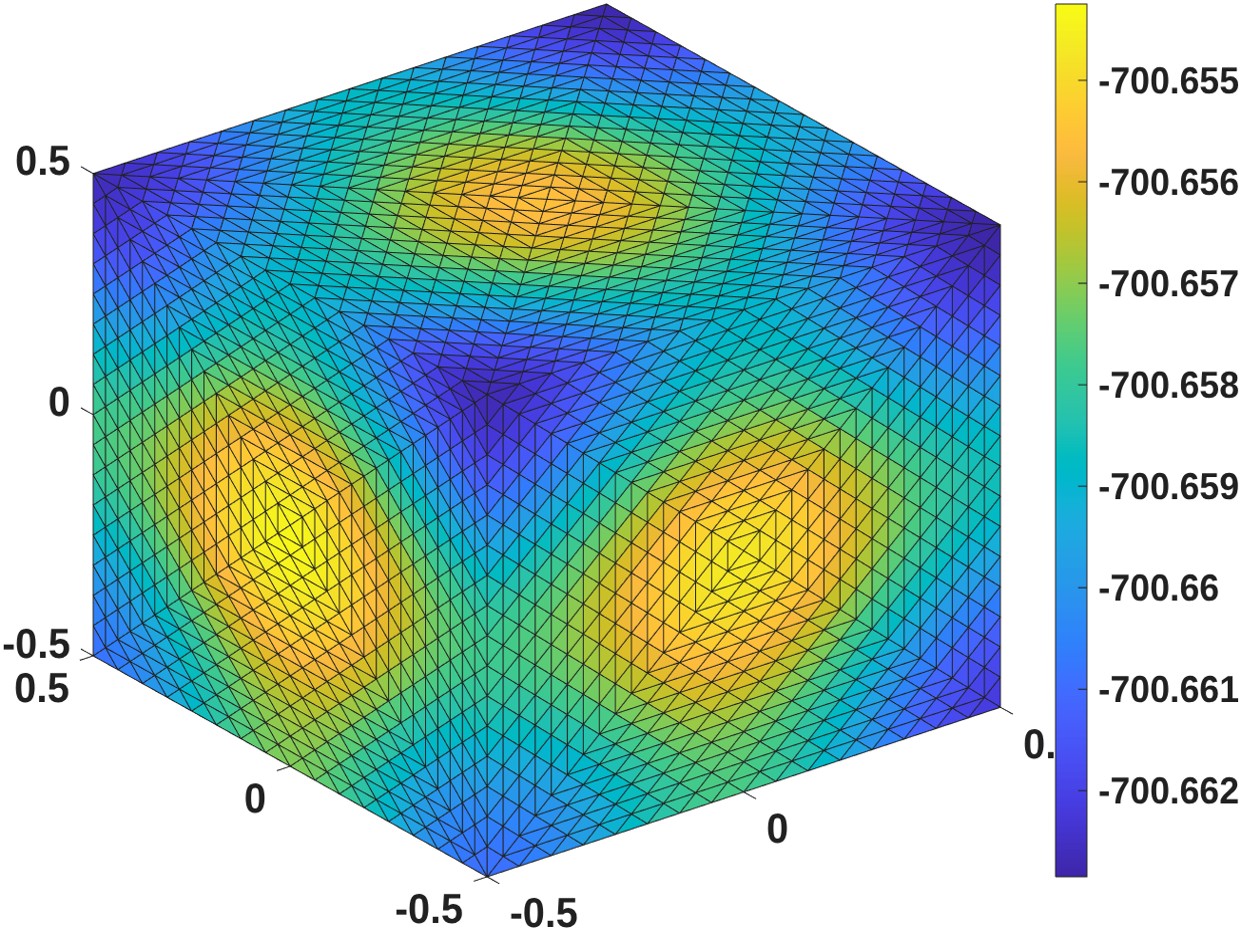}
	\end{minipage}\hskip0.8cm
	\begin{minipage}{0.45\textwidth}
		\centering
		\includegraphics[width=6.0cm]{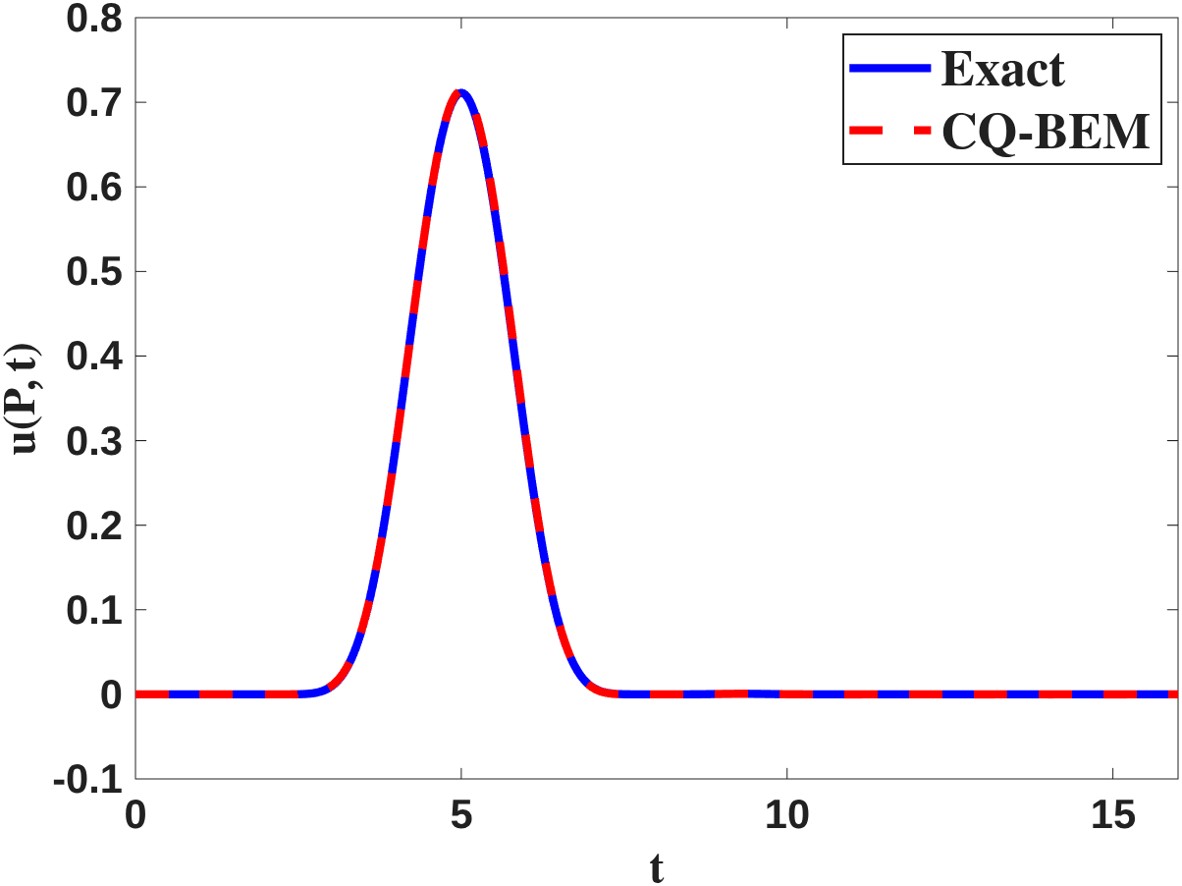}
	\end{minipage}
	\caption{\footnotesize{Example 4. Behaviour of the density $\varphi_{\Delta_t,h}^{(N)}$ on $\Gamma$ (left) and time history of $u_{\Delta_t,h}(\p,t)$ (right).}}\label{fig:cubo}
\end{figure}

\textbf{Example 4.} We consider the same wave defined by~\eqref{eq:exactsol} in Example~2. In this case, however, the wave impinges on a non-smooth domain, the cube \( \Omegap = (-0.5,\,0.5)^3 \) of  side length $1$. 
In Table \ref{table:8} we report the errors of the solution of Problem \eqref{pb:esterno} computed at $\p = (0,0,3)$.  
In Figure~\ref{fig:cubo}, we present two plots: the left one shows the density solution \( \varphi_{\Delta_t,h}^{(N)} \) on the surface \( \Gamma \) at the final time instant, while the right one compares the exact solution \( u(\p,t) \) with its numerical approximation \( u_{\Delta_t,h}(\p,t) \) at the fixed point \( \p \), for \( t \in [0,16] \). These results confirm the quadratic convergence order; furthermore, they are consistent with those obtained for regular domains, as well as with the results presented in \cite[Table~II]{Chappell}, which refer to the CQ-Galerkin method applied to a cylindrical obstacle generated by the rotation of a rectangular profile curve.

\setlength{\tabcolsep}{6pt}
\begin{table}[h]	
	\centering
	\footnotesize
		\begin{tabular}{ccccc}
			\hline
			\multicolumn{1}{c}{$\Delta_t$} & \multicolumn{1}{c}{$h$} & \multicolumn{1}{c}{$\operatorname{E}^u_{\Delta_t,h}(\p)$} & \multicolumn{1}{c}{EOC$^u$} & \multicolumn{1}{c}{$\kappa_2(\mathbf{A}^{(0)})$} \\
			\hline
			
			$2.50\times 10^{-1}$  & $5.00\times 10^{-1}$ & $1.61\times 10^{-1}$ & \multicolumn{1}{c}{---} & $3.28$ \\
			$1.25\times 10^{-1}$  & $2.50\times 10^{-1}$ & $5.31\times 10^{-2}$ & $1.60$ & $3.18$ \\
			$6.25\times 10^{-2}$ & $1.25\times 10^{-1}$ & $1.51\times 10^{-2}$ & $1.81$ & $3.21$ \\
			$3.13\times 10^{-2}$ & $6.25\times 10^{-2}$ & $3.72\times 10^{-3}$ & $2.02$ & $3.25$ \\
			\hline
		\end{tabular}
		\caption{\footnotesize{Example 4. Errors, convergence rates and condition numbers for the collocation method on the cube with $T = 16$.}}
			\label{table:8}
	\vspace{0.5\baselineskip}
\end{table}

\

\begin{figure}[h]
	\hskip-0.8cm\begin{minipage}[t]{0.32\textwidth}
		\centering
		\includegraphics[width=\linewidth]{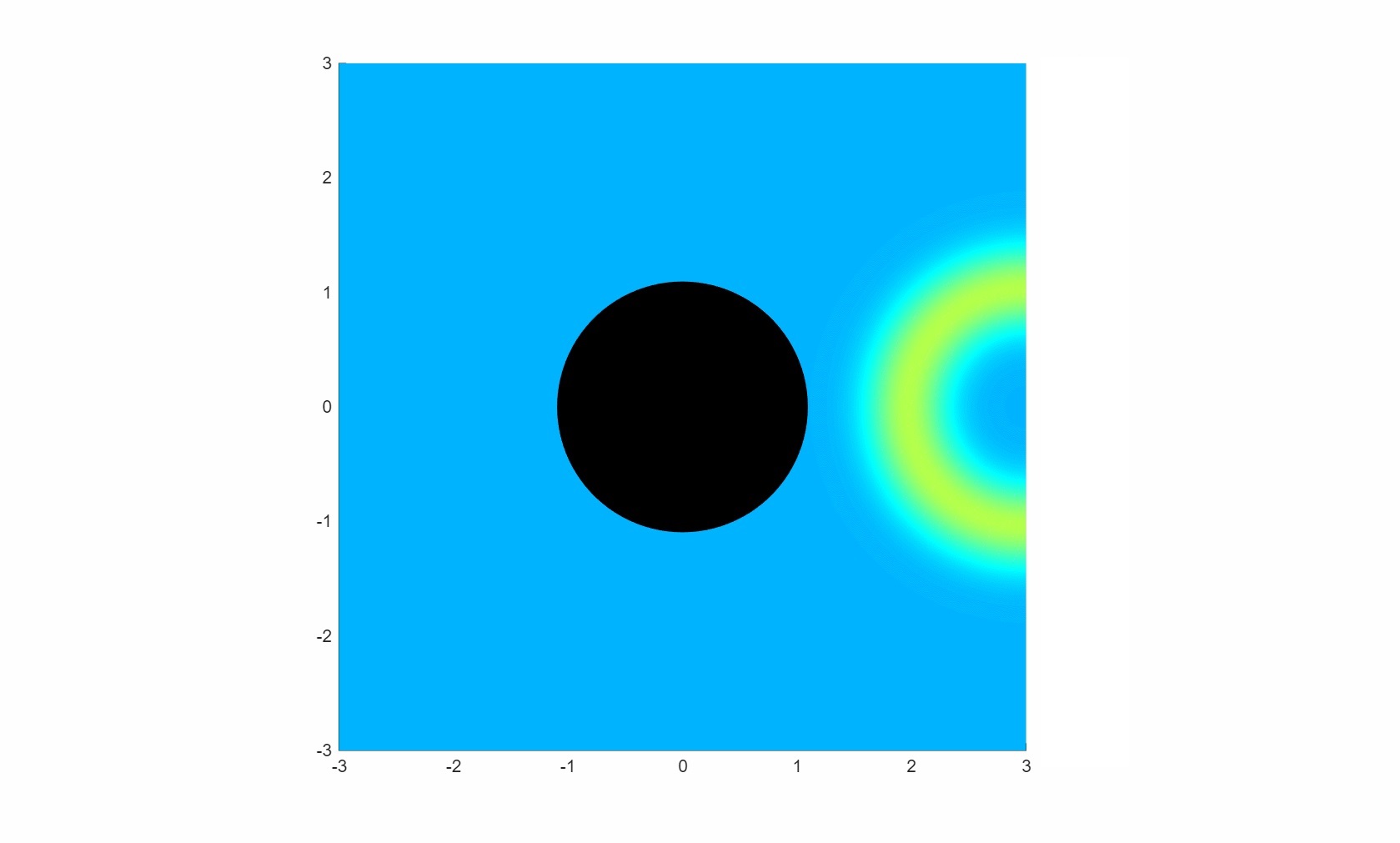}\par\vspace{2pt}
		
		\includegraphics[width=\linewidth]{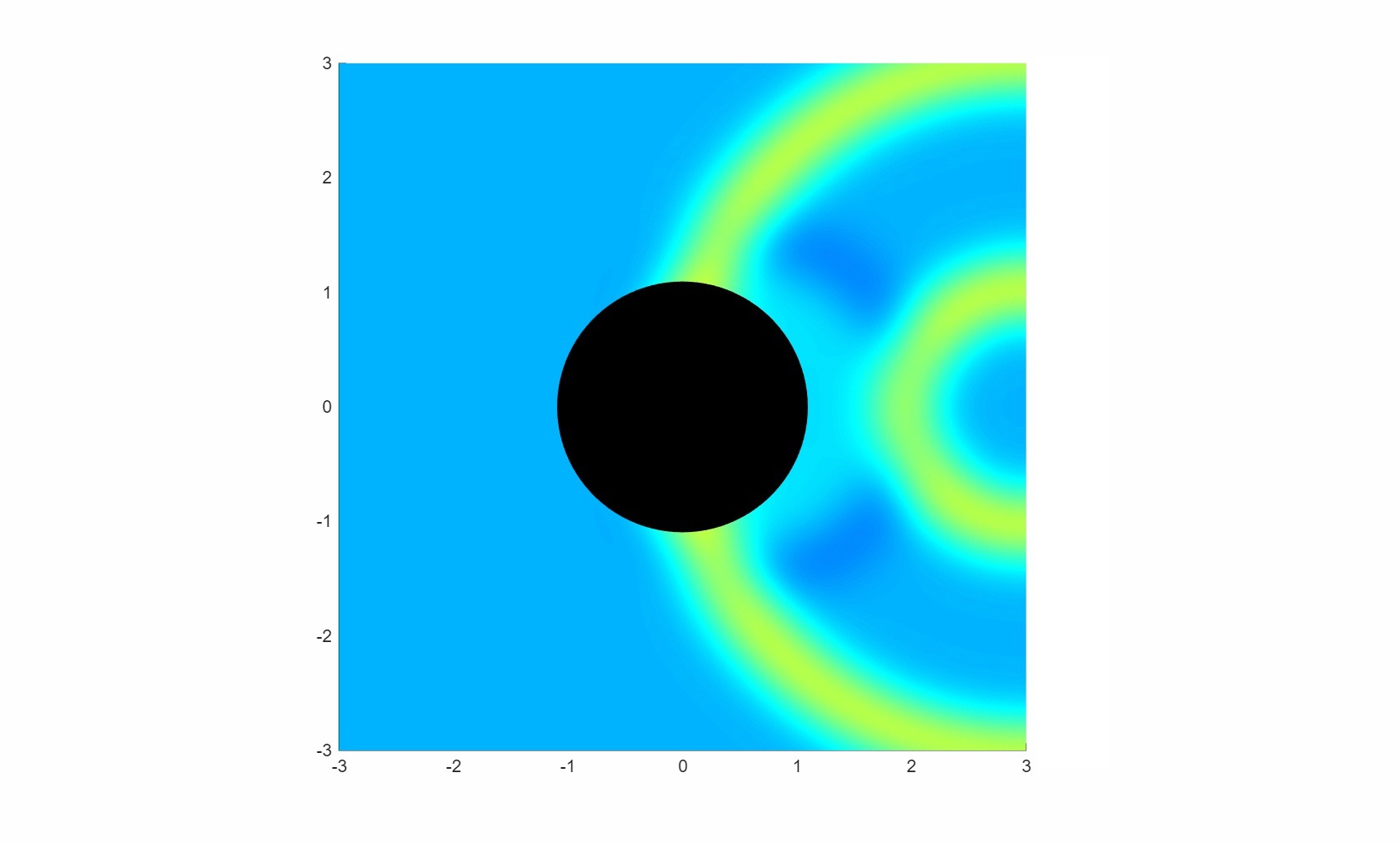}\par\vspace{2pt}
		
		\includegraphics[width=\linewidth]{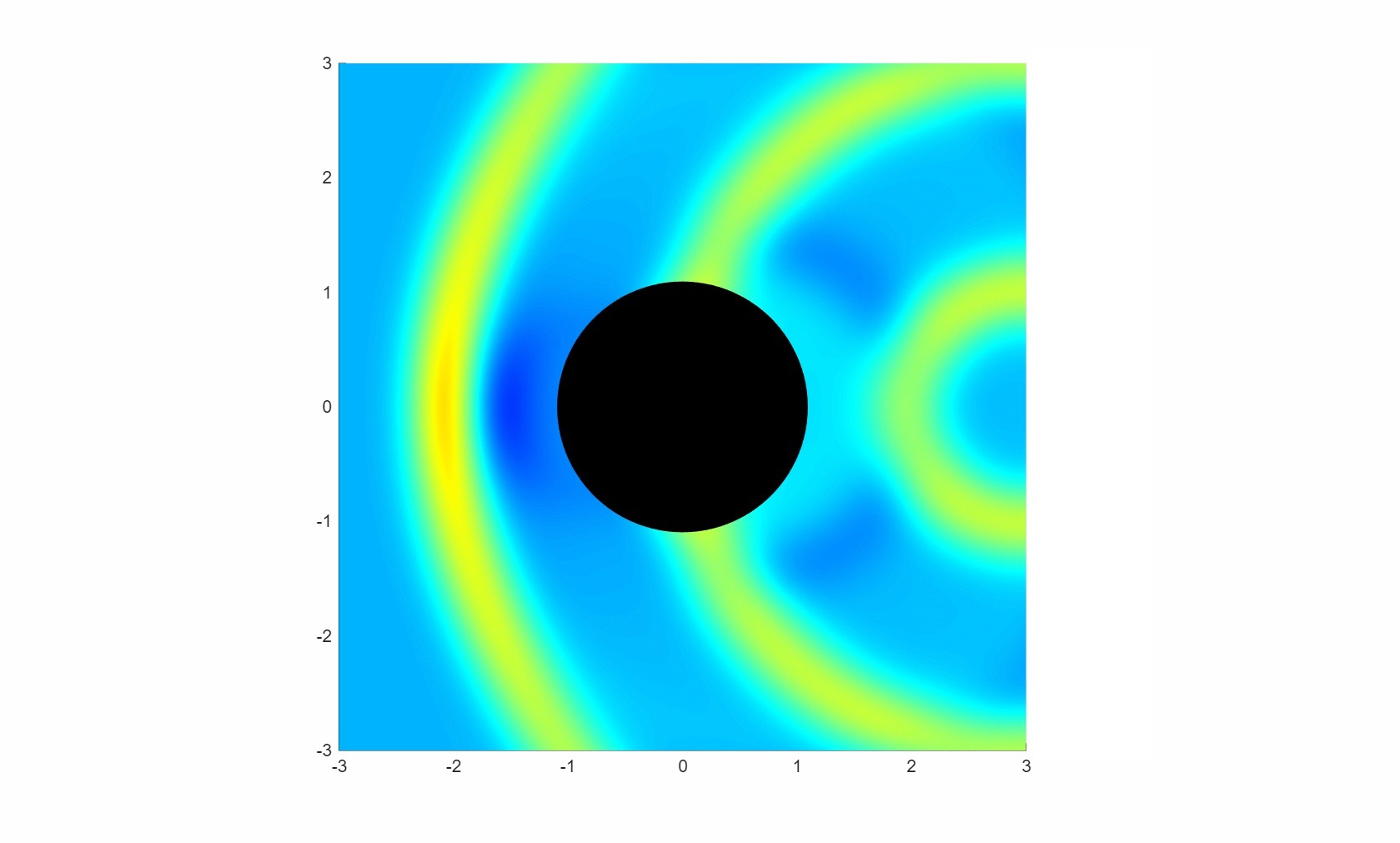}\par\vspace{2pt}
		
		\includegraphics[width=\linewidth]{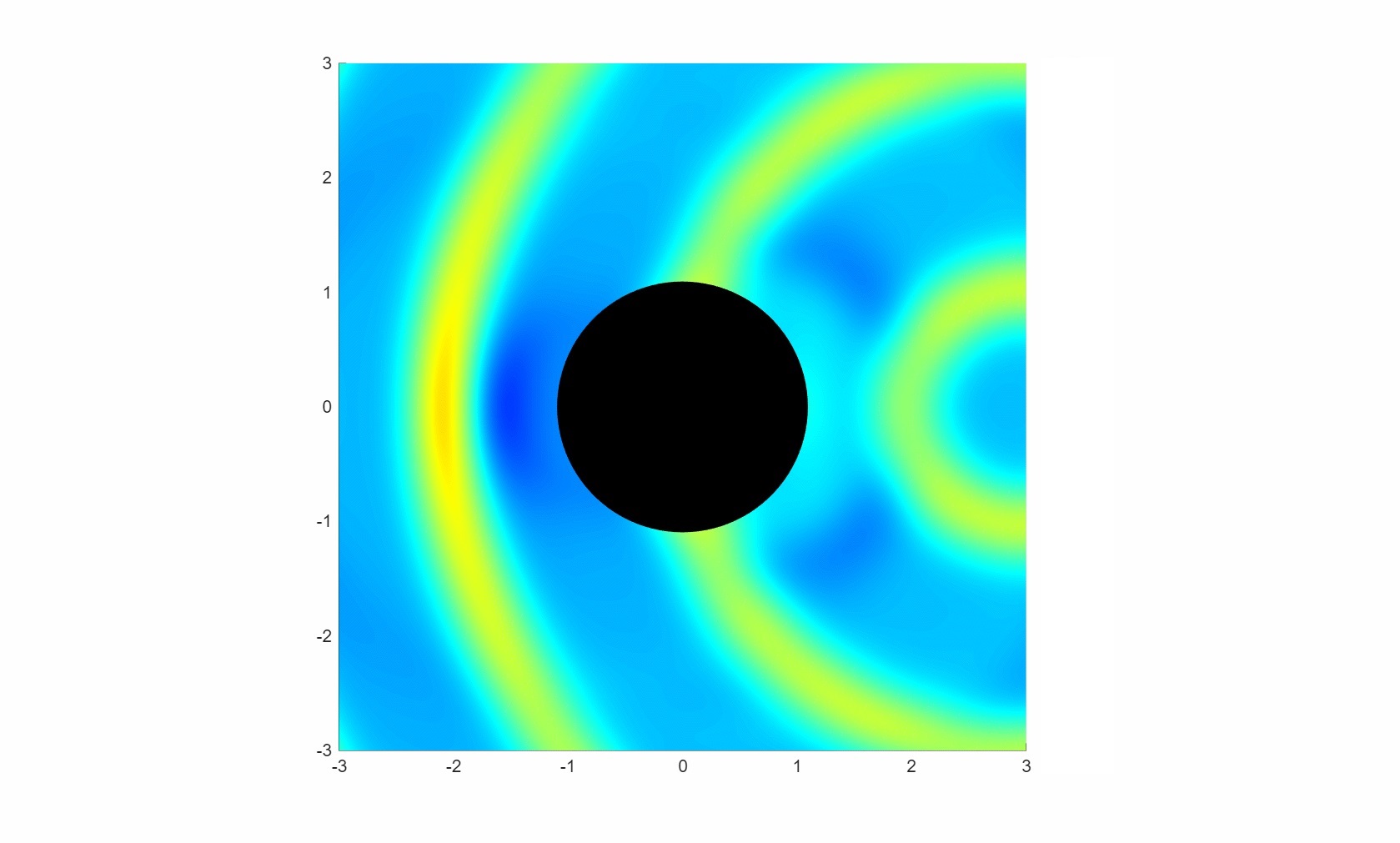}\par\vspace{2pt}
		
		\includegraphics[width=\linewidth]{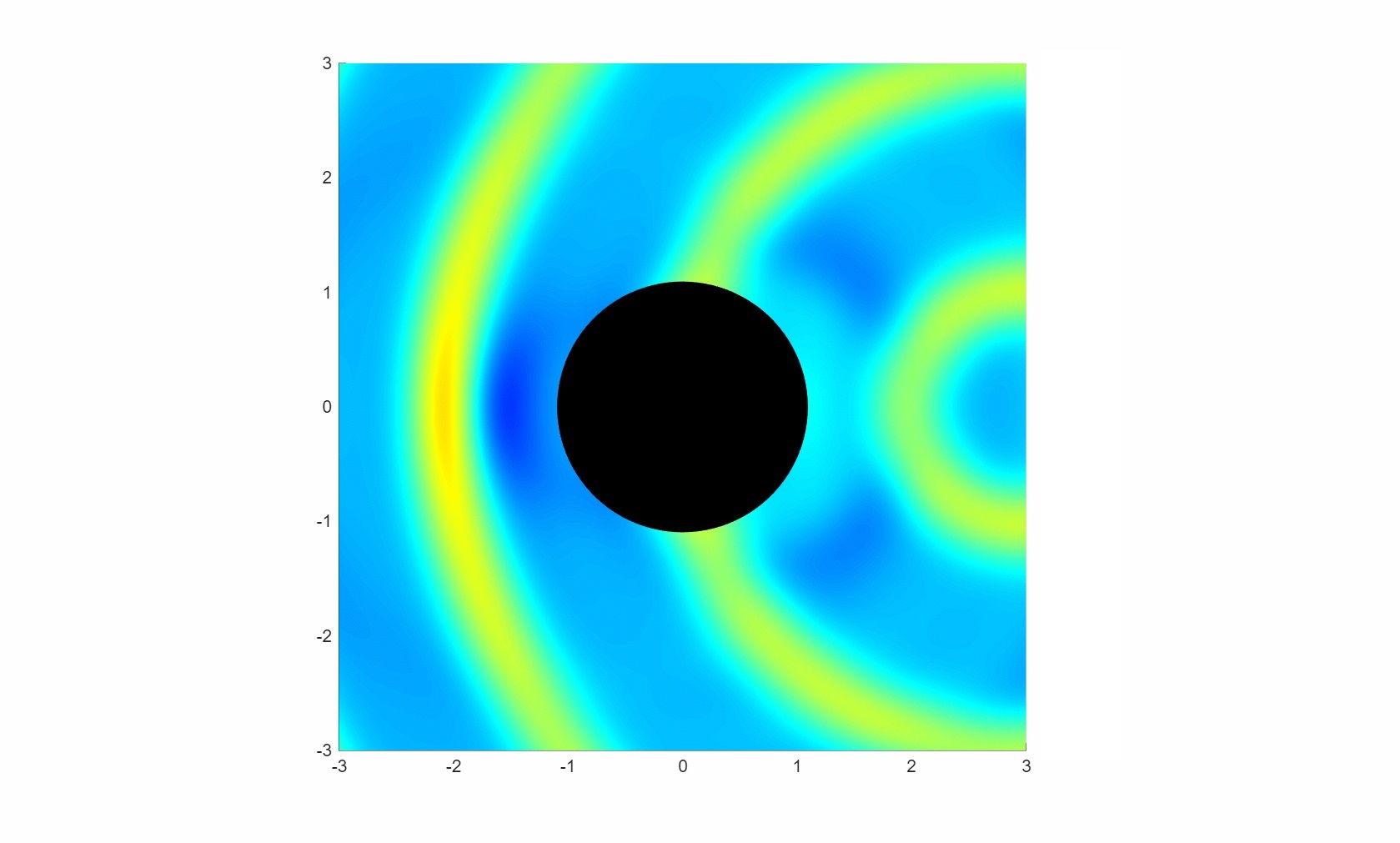}\par\vspace{2pt}
	\end{minipage}\hskip-1cm
	\begin{minipage}[t]{0.32\textwidth}
		\includegraphics[width=\linewidth]{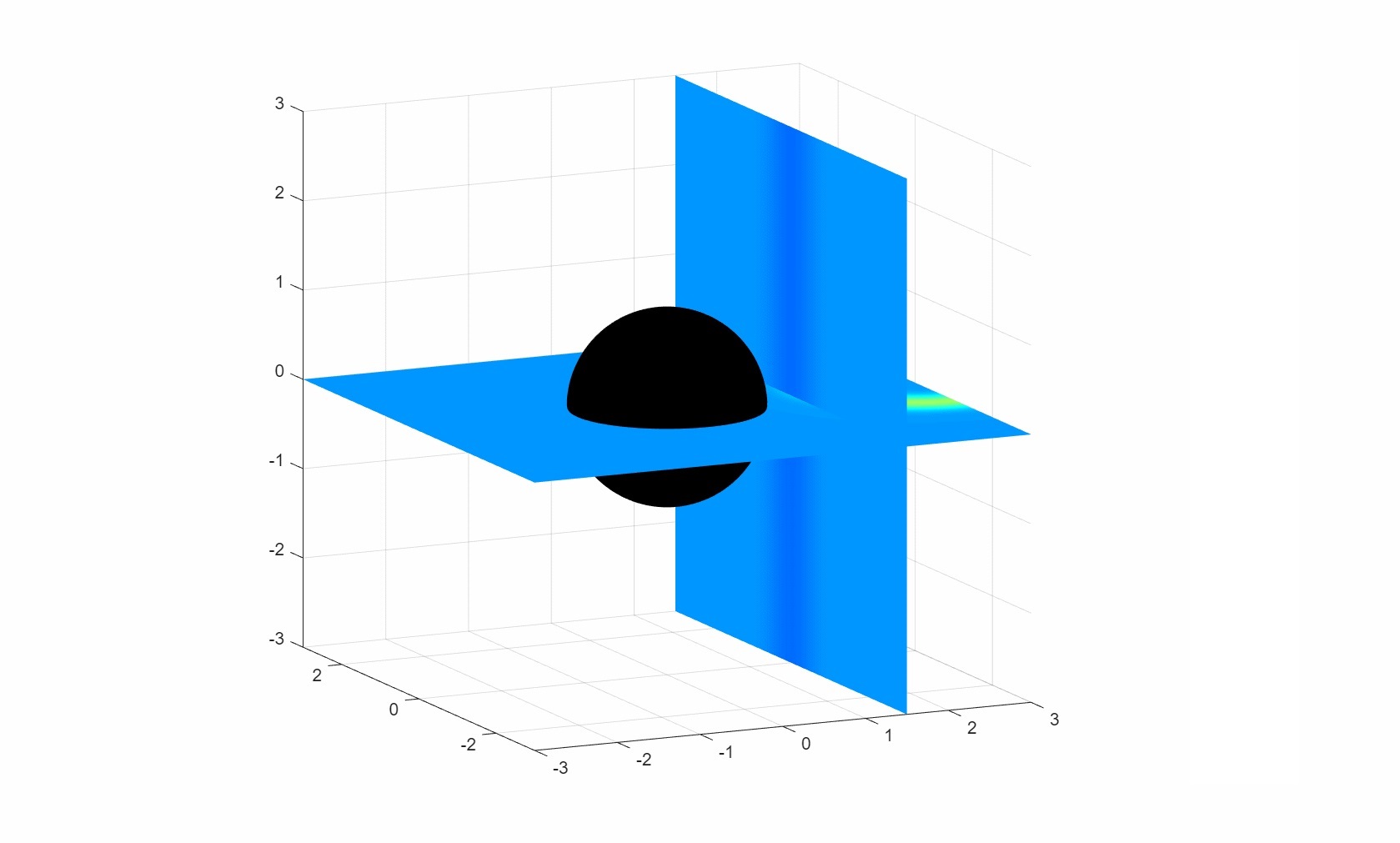}\par\vspace{2pt}
		
		\includegraphics[width=\linewidth]{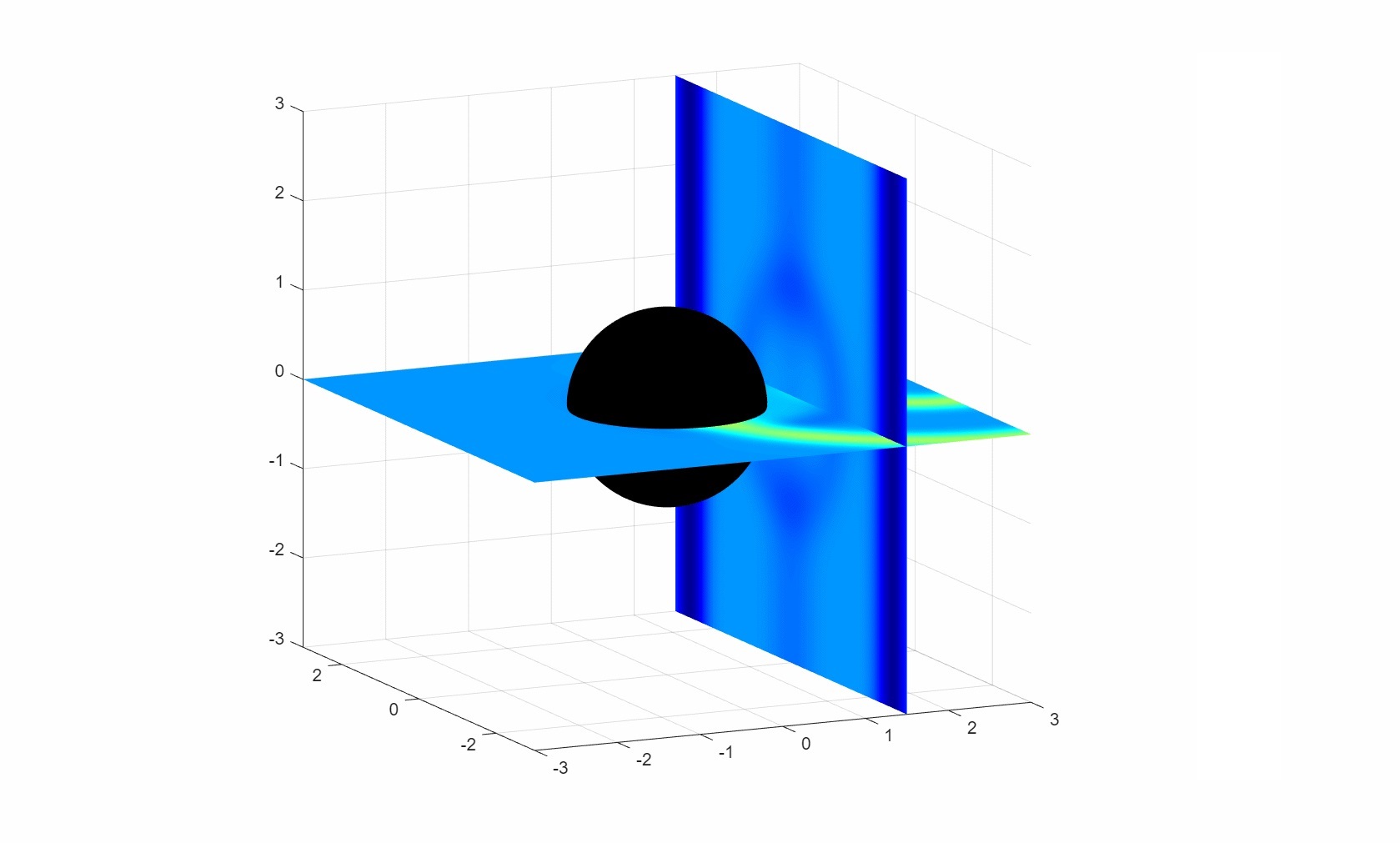}\par\vspace{2pt}
		
		\includegraphics[width=\linewidth]{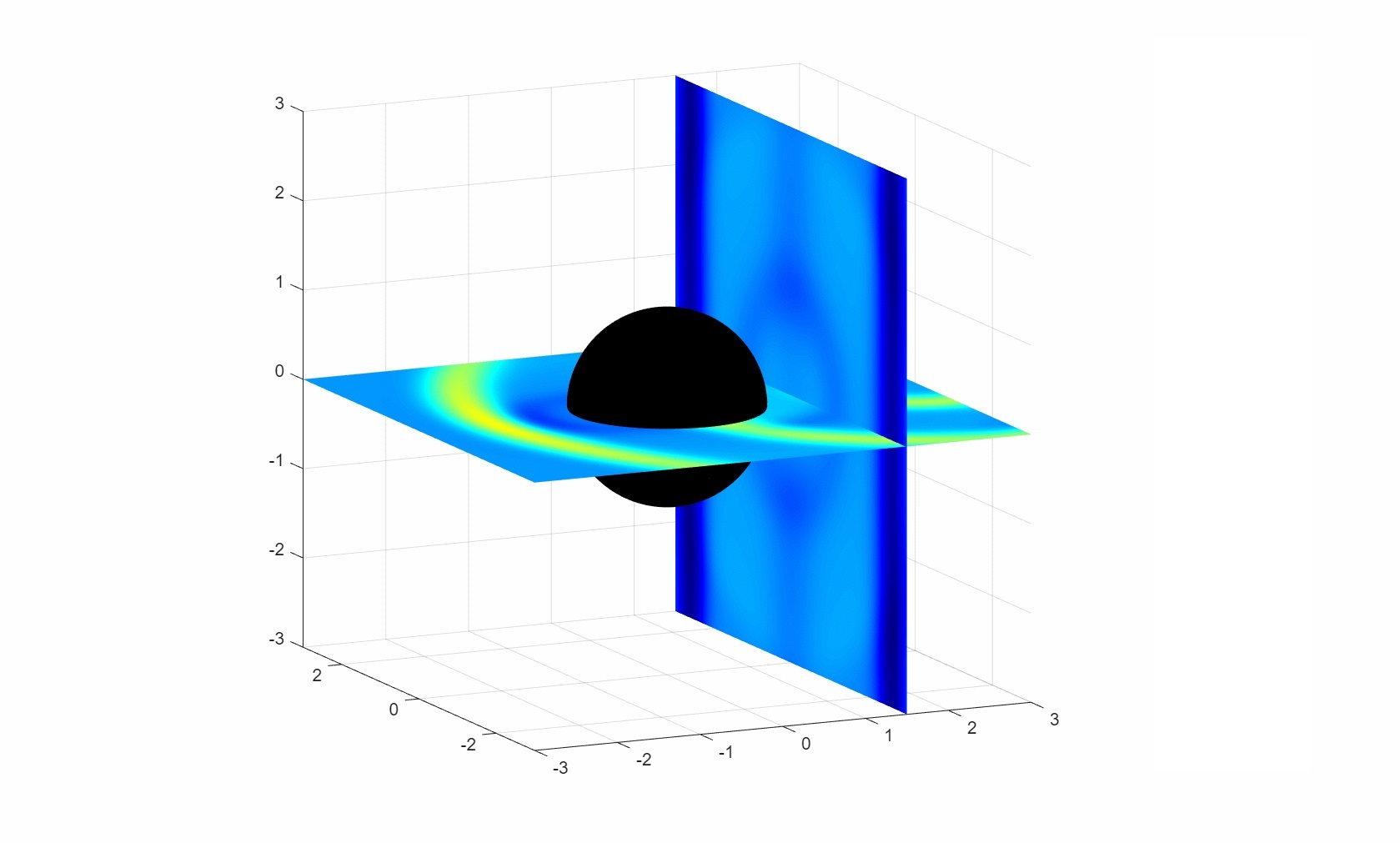}\par\vspace{2pt}
		
		\includegraphics[width=\linewidth]{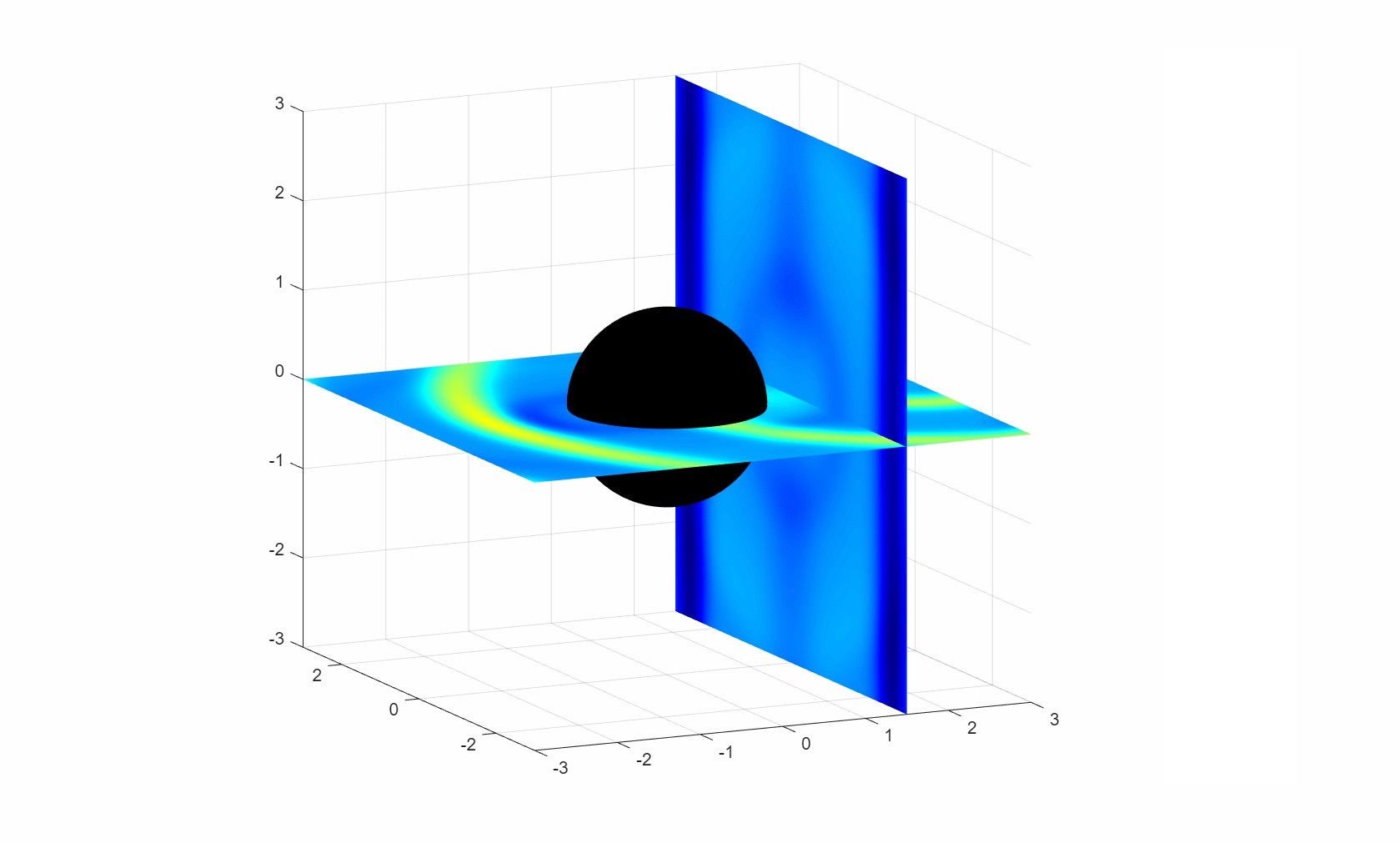}\par\vspace{2pt}    
		
		\includegraphics[width=\linewidth]{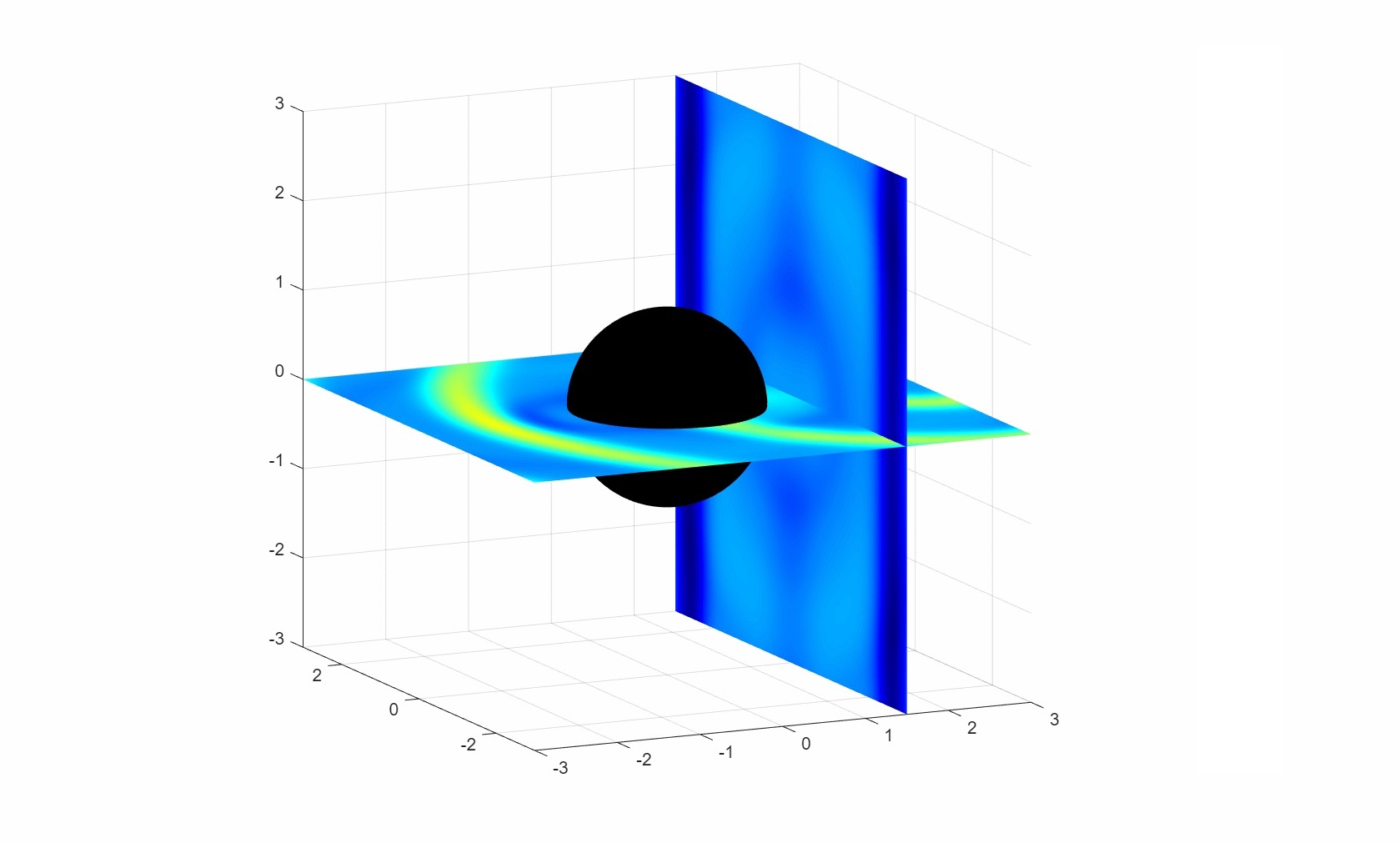}\par\vspace{2pt}    
	\end{minipage}\hskip-0.8cm
	\begin{minipage}[t]{0.32\textwidth}
		\includegraphics[width=\linewidth]{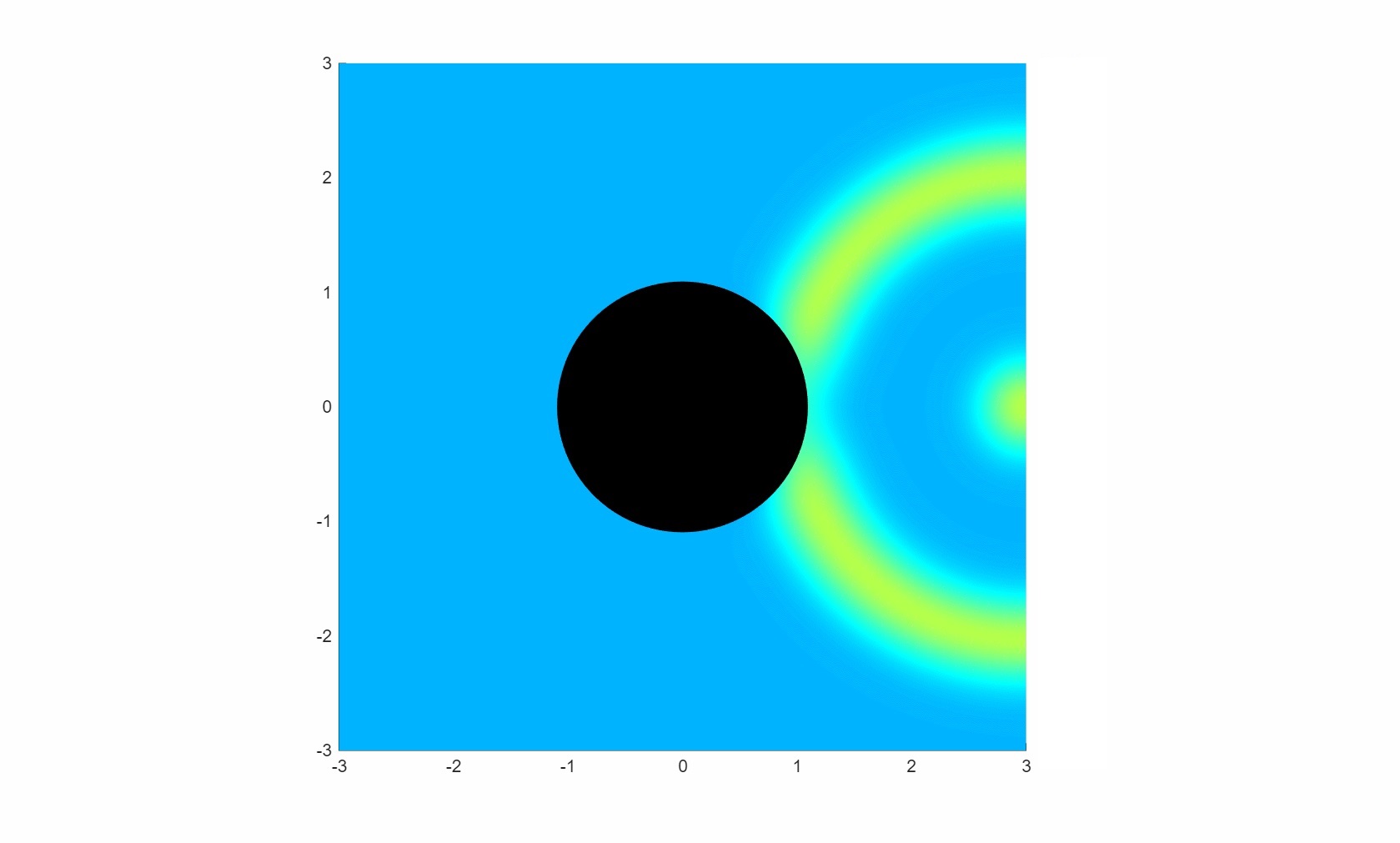}\par\vspace{2pt}
		
		\includegraphics[width=\linewidth]{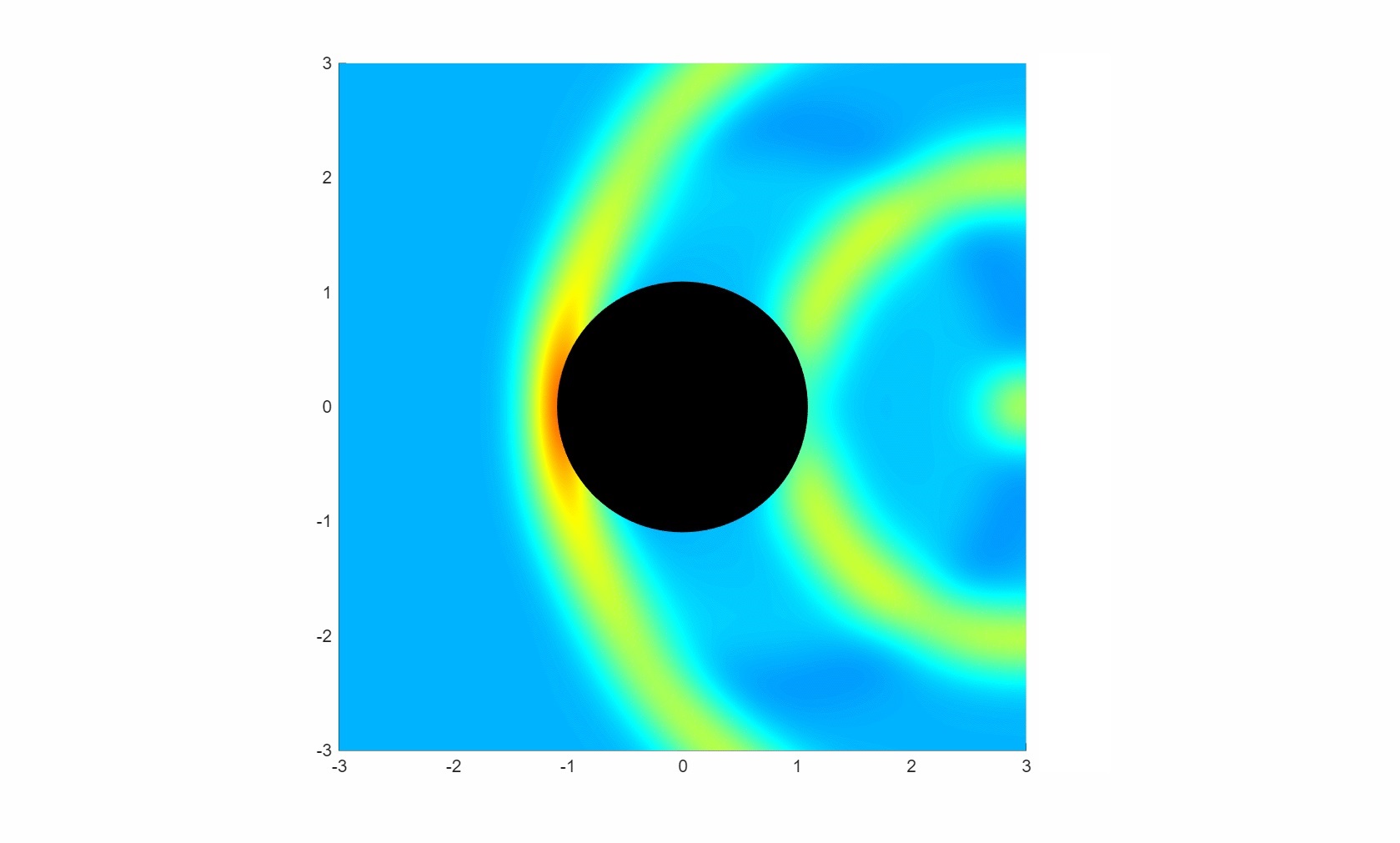}\par\vspace{2pt}
		
		\includegraphics[width=\linewidth]{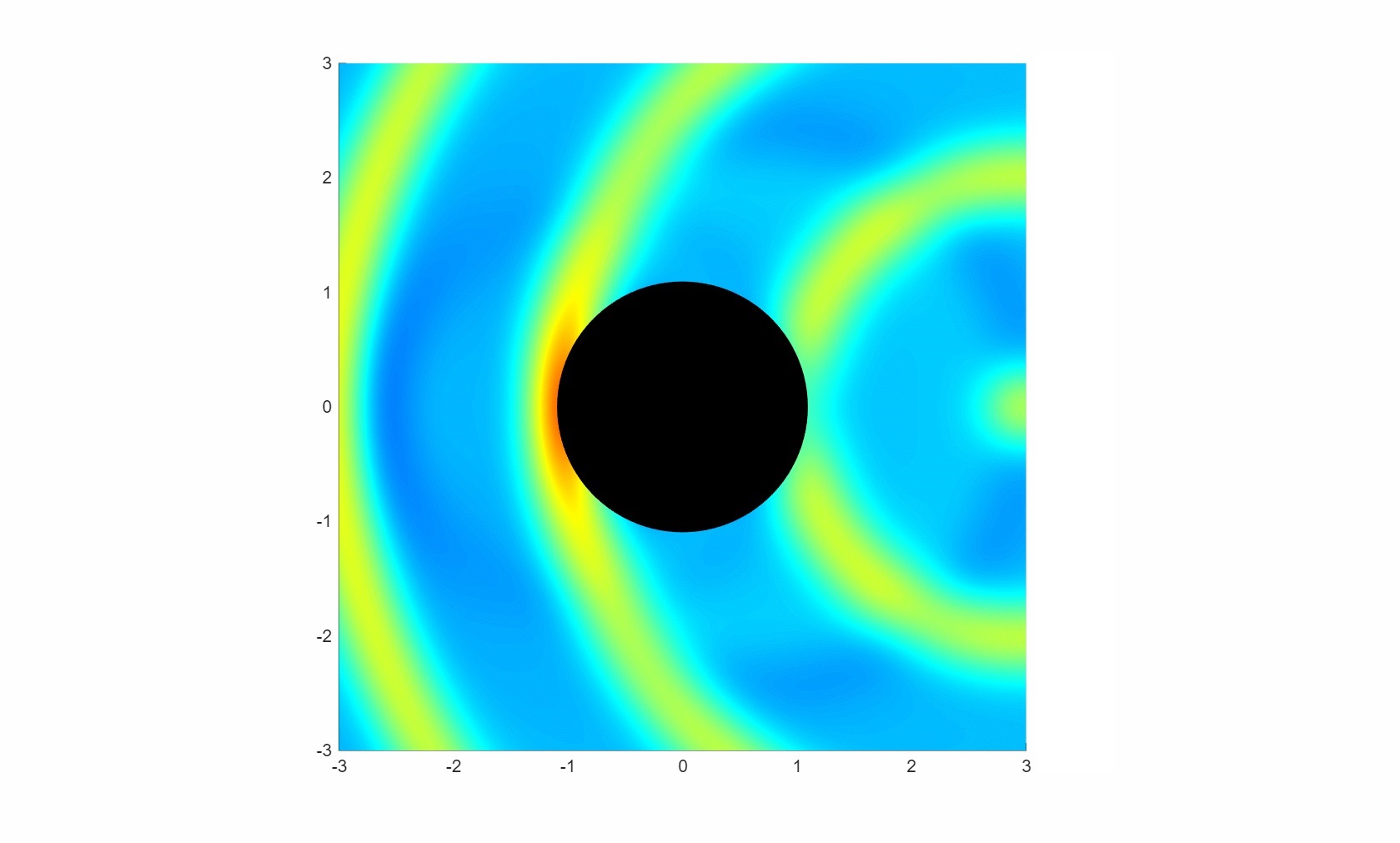}\par\vspace{2pt}
		
		\includegraphics[width=\linewidth]{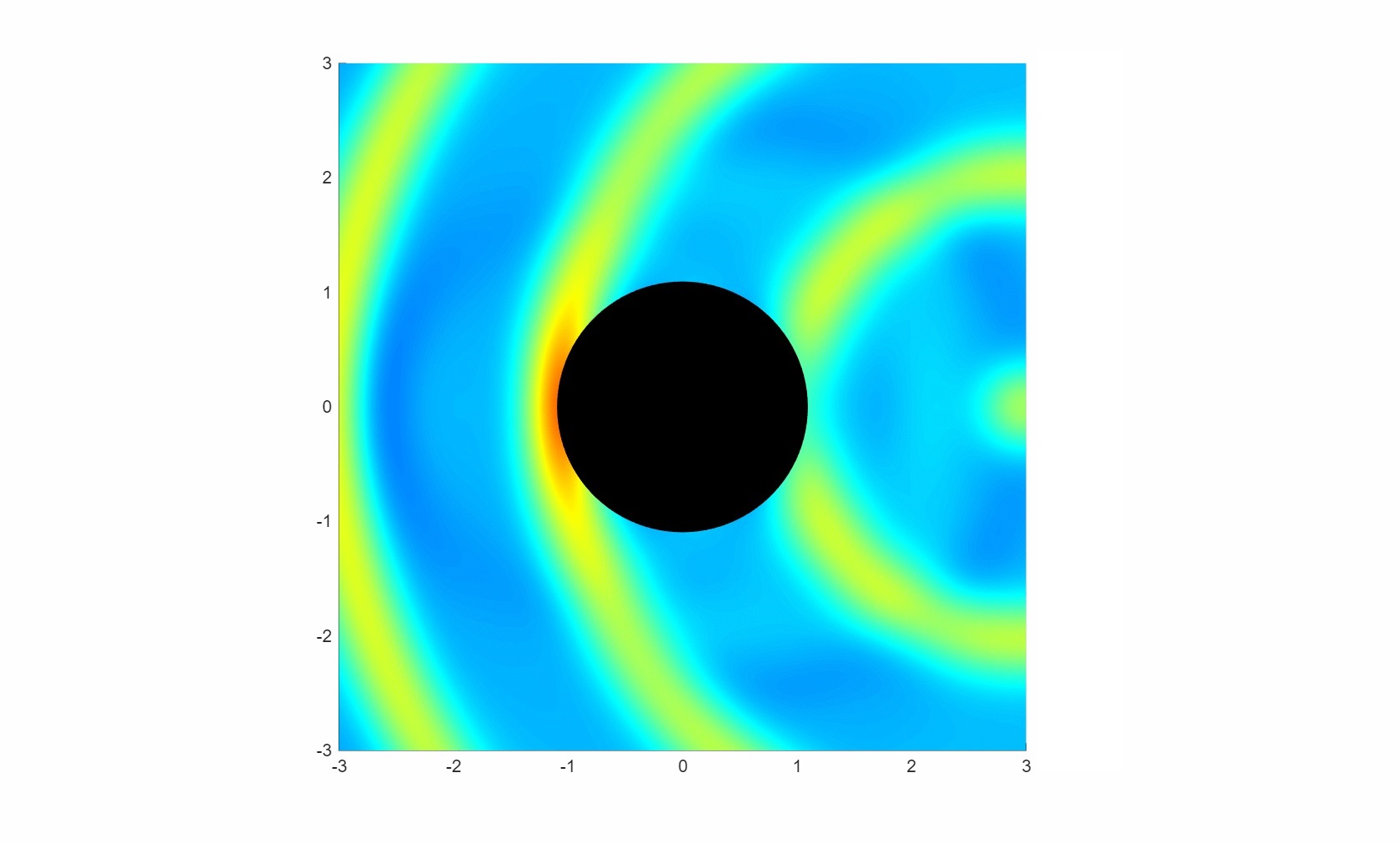}\par\vspace{2pt}    
		
		\includegraphics[width=\linewidth]{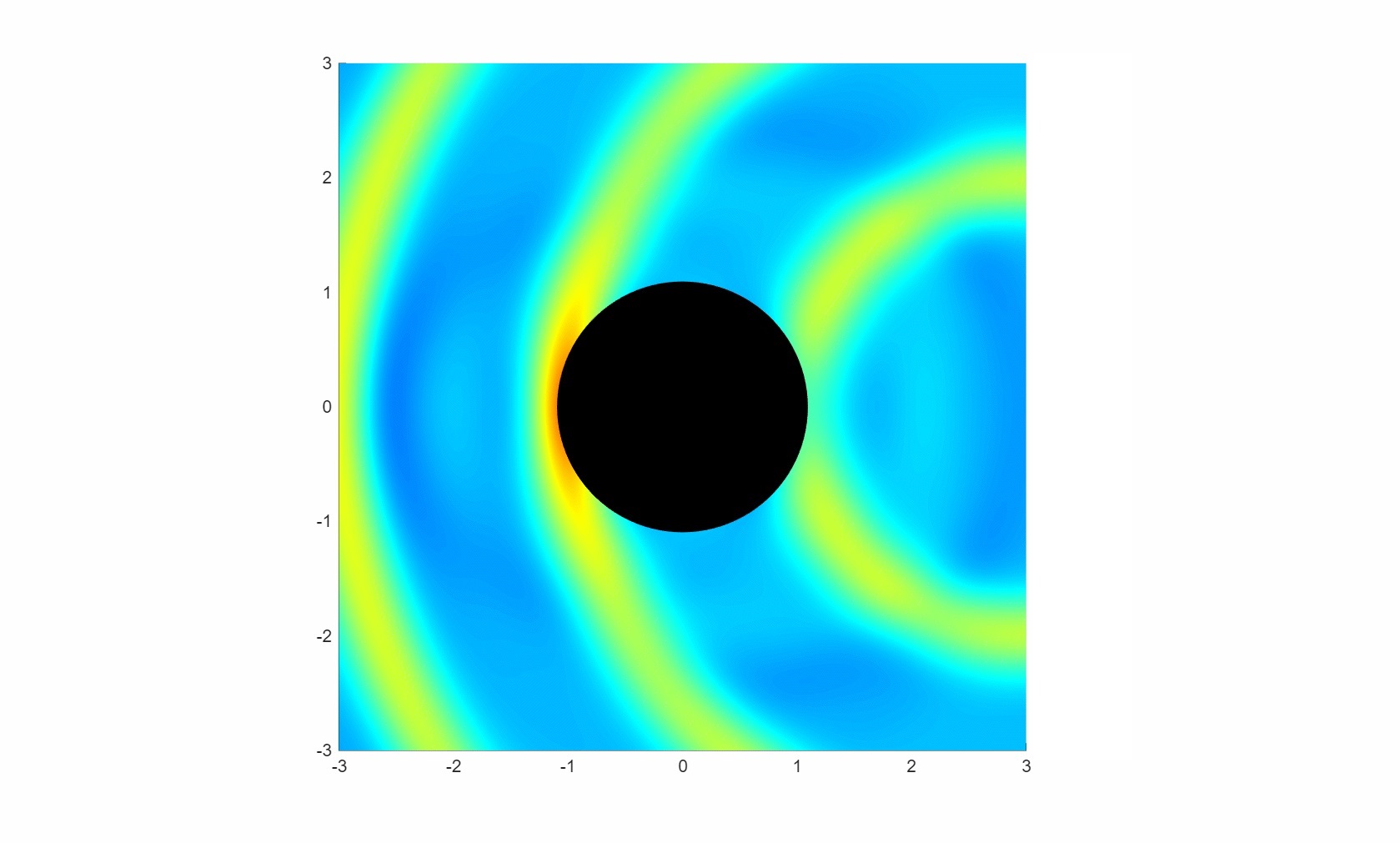}\par\vspace{2pt}    	
	\end{minipage}\hskip-1cm
	\begin{minipage}[t]{0.32\textwidth}
		\includegraphics[width=\linewidth]{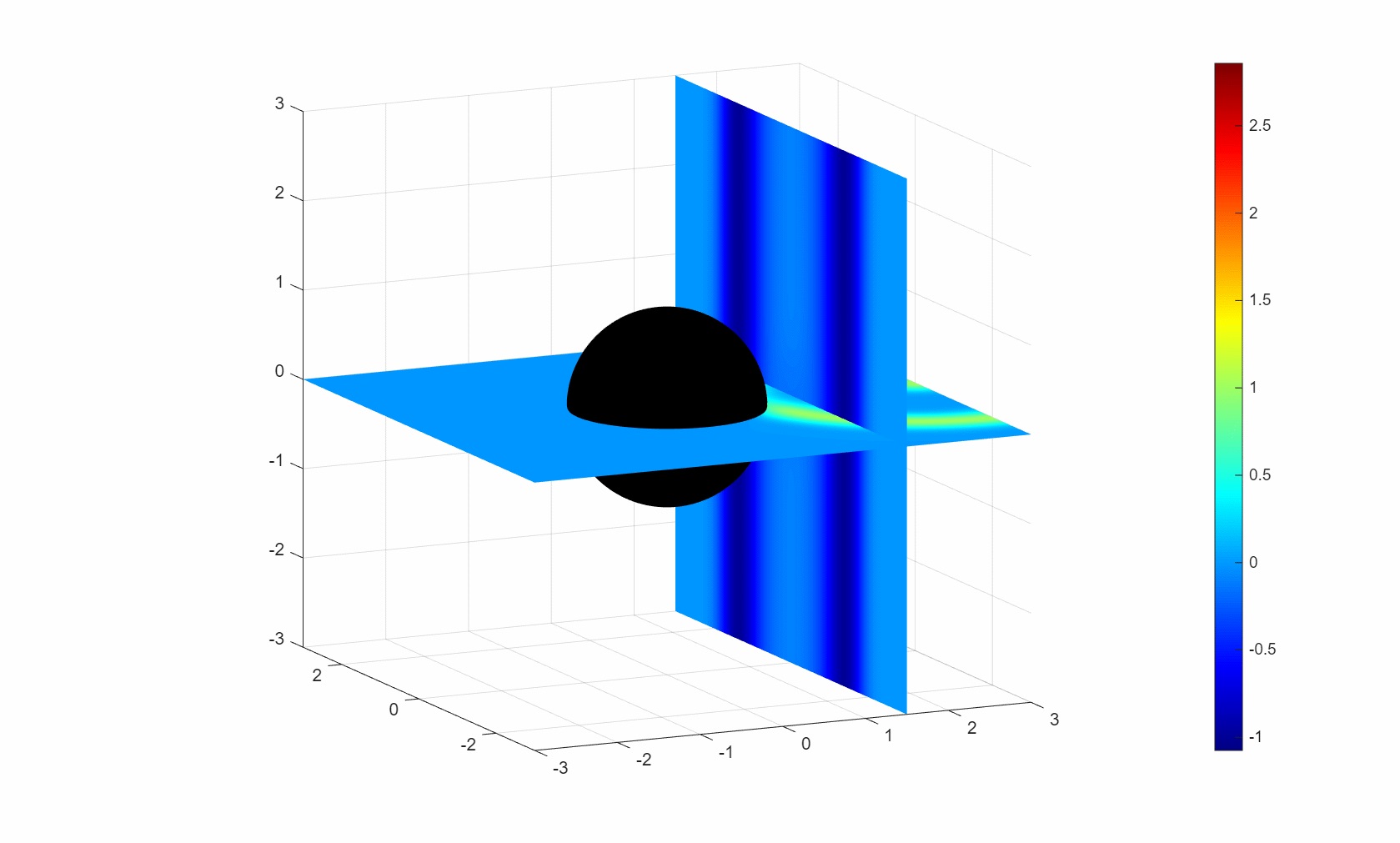}\par\vspace{2pt}
		
		\includegraphics[width=\linewidth]{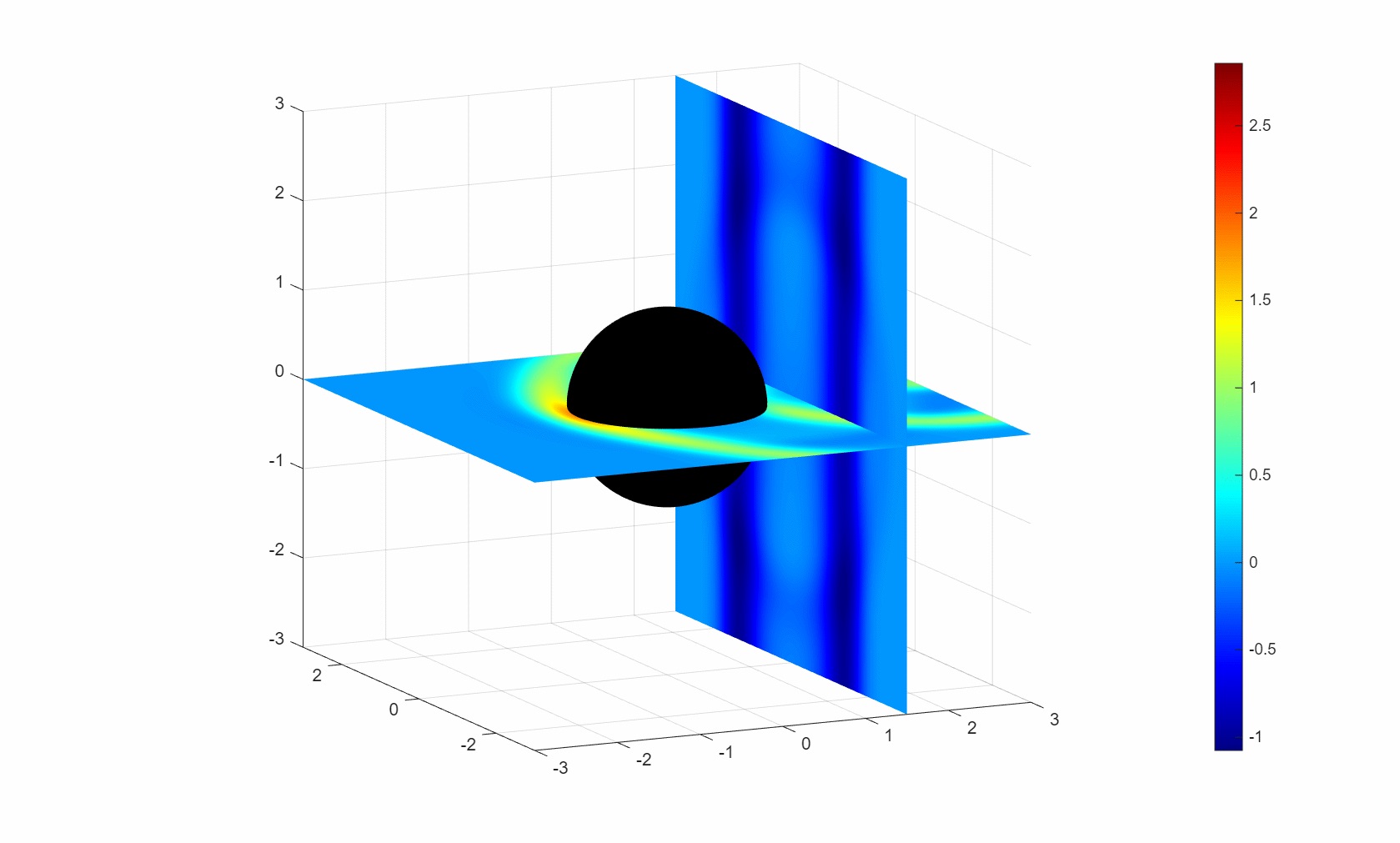}\par\vspace{2pt}
		
		\includegraphics[width=\linewidth]{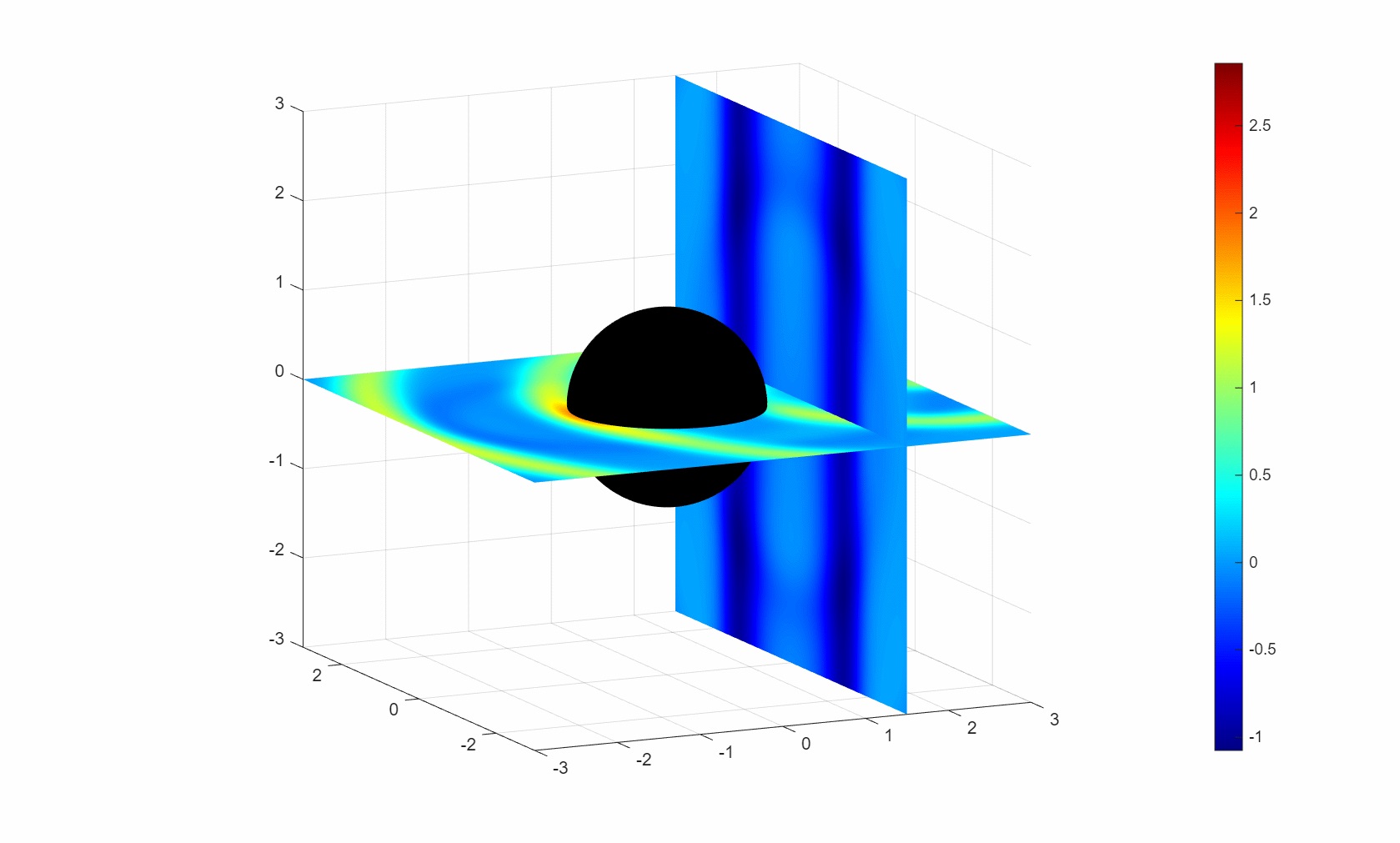}\par\vspace{2pt}
		
		\includegraphics[width=\linewidth]{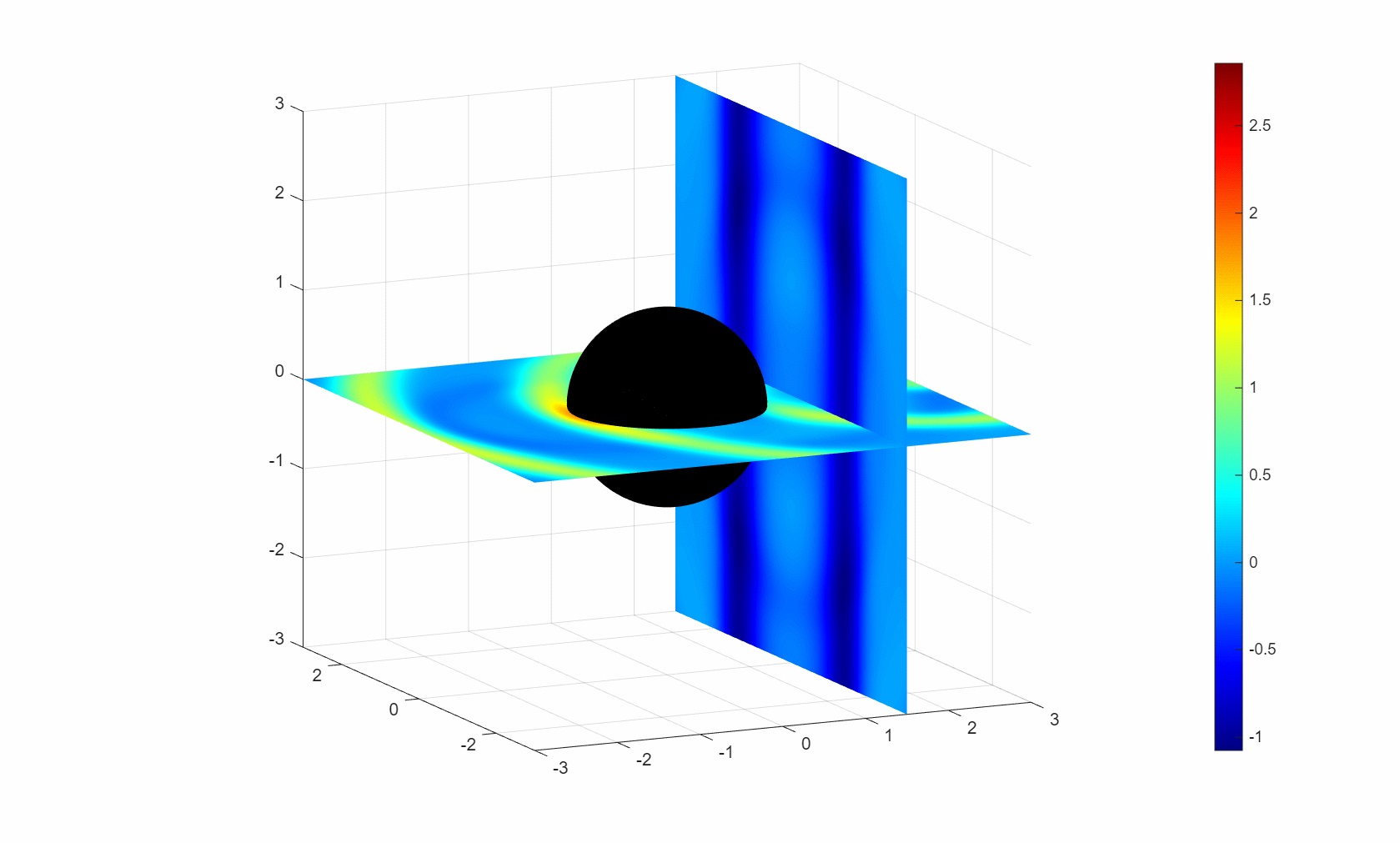}\par\vspace{2pt}    
		
		\includegraphics[width=\linewidth]{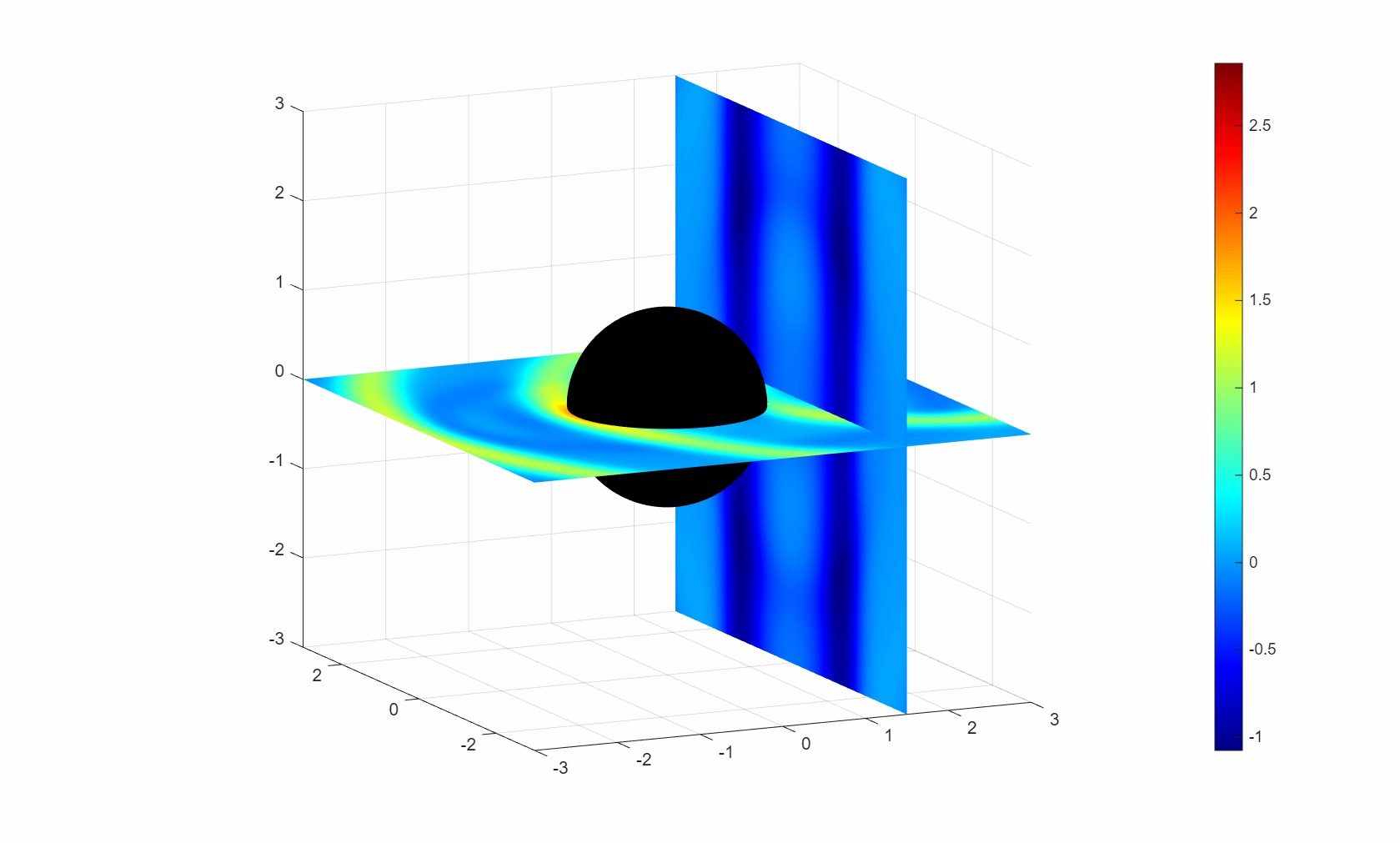}\par\vspace{2pt}    
	\end{minipage}
	\caption{\footnotesize{Example 5. 2D and 3D snapshots of the solution at selected instants.}}\label{fig:snaps2D3D}
\end{figure}

\textbf{Example 5.} As a final example, we present a physically meaningful scattering problem. Specifically, we consider the propagation of a superposition of five cylindrical wave packets
\begin{align*}
u^{\text{inc}}(\x,t) &= \sum_{n=0}^4 e^{-8\left(\sqrt{(x_1-3)^2+x_2^2}-(t-2n)\right)^2}, \quad t\in [0,10]
\end{align*}
each propagating in the \( (x_1, x_2) \)-plane, invariant along the \( x_3 \)-axis and impinging on the unit sphere. 
The wave is centred at $\x = (3,0,0)$ and  
consists of Gaussian pulses emitted at equidistant time intervals \( t = 0, 2, 4, 6, 8 \). Each packet travels at unit speed along circular wavefronts. 
The total wave field $u^{\text{tot}}$ is given by the superposition of the incident wave $u^{\text{inc}}$ and the scattered one $u$, where the latter satisfies Problem~\eqref{pb:esterno} with  Neumann boundary condition on the surface $\Gamma$, given by
$g(\mathbf{x}, t) = -\partial u^{\text{inc}}(\mathbf{x}, t)/\partial \mathbf{n}_\x$.
For the numerical solution of this problem, we apply the Lubich shifted collocation method, decomposing the boundary $\Gamma$ into $\mathcal{N}_\Gamma = 320$ triangular elements and partitioning the time interval into $N = 512$ uniform subintervals. The scattered field is then evaluated within two planar observation domains defined as follows
\[\text{\Pisymbol{psy}{112}}_1 = \left\{(x_1,x_2,x_3): (x_1,x_2) \in [-3,3]^2, \ x_3 = 0 \right\}\cap \Omegam\]
and
\[\text{\Pisymbol{psy}{112}}_2 = \left\{(x_1,x_2,x_3): x_1 = 1.5, \ (x_2,x_3) \in [-3,3]^2\right\}\cap \Omegam.\]
In Figure \ref{fig:snaps2D3D}, we illustrate the time evolution of the total wave field through a sequence of snapshots at selected time instants $t_n \approx 1, 2, \ldots, 10$, displayed in lexicographic order. The two-dimensional plots show the field restricted to the horizontal plane $\text{\Pisymbol{psy}{112}}_1$, while the three-dimensional ones provide a combined view of the wave propagation on both $\text{\Pisymbol{psy}{112}}_1$ and $\text{\Pisymbol{psy}{112}}_2$. These results show that the proposed Lubich shifted-collocation approach provides an accurate and stable tool for the simulation of time-domain scattering phenomena, even in the presence of multiple interacting wave fronts.

\section{Conclusions}\label{sec:conclusions}
In this paper, we have investigated the application of a Lubich convolution quadrature combined with a collocation boundary element method for time-domain wave propagation problems with non-homogeneous Neumann boundary conditions. A shifted collocation strategy has been employed, with particular emphasis on implementation, efficiency, stability, and convergence of the resulting scheme.

Although a formal stability and convergence analysis is lacking, the numerical experiments provide evidence of unconditional stability of the proposed method and exhibit second-order convergence, provided that the boundary integrals are evaluated with sufficient accuracy. 

The applicability and robustness of the approach have been demonstrated through benchmark cases, including a realistic scattering problem. The results confirm that the CQ-collocation boundary element method is an effective and flexible tool for the time-domain simulation of wave phenomena.

\section*{Acknowledgements}
The authors are sincerely grateful to the anonymous referees for their insightful comments and constructive suggestions, which have significantly improved the manuscript.

\section{Appendix}\label{sec:appendix}
\subsection{Analytical solution on the unit sphere}\label{sec:anal_sol}
Let \(\Omega^i\) denote the unit sphere, and consider the Neumann boundary datum on its surface given by \(g(\mathbf{x}, t) =  Y_n^m(\mathbf{x})\tilde{g}(t)\), for a suitable $\tilde{g}$ and \(Y_n^m\) spherical harmonics of degree \(n\) and order \(m\), with \(n = 0, 1, \ldots\) and \(-n \leq m \leq n\). 

To derive an explicit expression for the solution \( u \) of  the exterior Neumann problem \eqref{pb:esterno}, we proceed similarly to the approach presented in~\cite{SilviaOnde2015}, where a Dirichlet wave problem was investigated. Specifically, we observe that in this setting the solution can be written as
\begin{equation}\label{eq:udecomposition}
u(\mathbf{x}, t) = \tilde{u}(t) Y_n^m(\boldsymbol{\xi}), \quad \text{with } \mathbf{x} = r \boldsymbol{\xi} \in \mathbb{R}^3,\ \boldsymbol{\xi} \in \Gamma, 
\end{equation}
where, for notational simplicity, we denote by $\tilde{u}(t)$ a function which depends also on $r$. We consider the Laplace transform of \( u(\mathbf{x}, t) \), denoted by \( \hat{u}(\mathbf{x}, s) \), which satisfies the wave equation in the frequency domain.
Using known properties of the Helmholtz equation, it follows that \( \hat{u} \) satisfies the relation
\begin{equation}\label{eq:uhat}
\hat{u}(s) = \widehat{\tilde{g}}(s) \, \frac{h_n^{(1)}(\imath sr)}{\imath s \, (h_n^{(1)})'(\imath s)},
\end{equation}
where \( h_n^{(1)} \) denotes the spherical Hankel function of the first kind, and \( (h_n^{(1)})' \) its derivative (see~\cite{Nedelec2001}). By applying the inverse Laplace transform to~\eqref{eq:uhat}, we obtain the representation~\eqref{eq:udecomposition}, with
\begin{equation*}
\tilde{u}(t) = \mathcal{L}^{-1} \left( \widehat{\tilde{g}}(s) \frac{h_n^{(1)}(\imath r s)}{\imath s \, (h_n^{(1)})'(\imath s)} \right)(t)
= \tilde{g} * \mathcal{L}^{-1} \left( \frac{h_n^{(1)}(\imath r \cdot)}{(\imath \cdot) \, (h_n^{(1)})'(\imath \cdot)} \right)(t),
\end{equation*}
where \( * \) denotes the time convolution.
Using known properties of spherical Hankel functions (see in particular~\cite[Sec.~10.1.1]{AbraStegun} and~\cite[Sec.~8.466]{GradRyz}), and following ~\cite{SauterVeit2013}, we obtain
\begin{align*}
&h_n^{(1)}(s) = \imath^{-(n+1)} \frac{1}{s} e^{\imath s} y_n\left(-\frac{1}{\imath s}\right), \\
&(h_n^{(1)})'(s) = \imath^{-(n+1)} \frac{e^{\imath s}}{s^3} \left[ s (\imath s - 1) y_n\left(-\frac{1}{\imath s}\right) - \imath y_n'\left(-\frac{1}{\imath s}\right) \right],
\end{align*}
where the polynomials \( y_n \) and their derivatives are given by
\begin{align*}
y_n(s) &:= \sum_{k=0}^n (n, k) \, s^k, \quad 
y_n'(s) := \sum_{k=1}^n (n, k) \, k \, s^{k-1}, \quad 
(n, k) := \frac{(n + k)!}{2^k \, k! \, (n - k)!}.
\end{align*}
In Example 1, we have considered the case \( n = 0 \), for which we have \( y_0(s) = 1 \) and \( y_0'(s) = 0 \). This leads to the simplified expression
\begin{equation*}
\frac{h_0^{(1)}(\imath sr)}{\imath s \, (h_0^{(1)})'(\imath s)} = -\frac{e^{-s(r-1)}}{r(s+1)}.
\end{equation*}
Applying well known inverse Laplace transform properties, we obtain
\begin{equation*}
\mathcal{L}^{-1} \left( \frac{h_0^{(1)}(\imath r \cdot)}{(\imath \cdot)(h_0^{(1)})'(\imath \cdot)} \right)(t)
= -\frac{1}{r} e^{-(t - (r - 1))} H(t - (r - 1)),
\end{equation*}
where \( H \) denotes the Heaviside step function. Therefore, the time-dependent part of the solution becomes
\begin{align}\label{eq:uesatta}
\tilde{u}(t) 
&= -\frac{1}{r} \int_0^t \tilde{g}(t - \tau) \, e^{-(\tau - (r - 1))} H(\tau - (r - 1)) \, d\tau. 
\end{align}
Finally, recalling that \( Y_0^0(\mathbf{x}) = 1/(2\sqrt{\pi}) \), the solution of Problem \eqref{pb:esterno} is
\begin{equation*}
u(\mathbf{x}, t) = \frac{1}{2\sqrt{\pi}} \tilde{u}(t),
\end{equation*}
with \( \tilde{u}(t) \) given by~\eqref{eq:uesatta}.

\end{document}